%% file: main.tex
\documentclass[journal]{new-aiaa}
\input{Preamble}
\title{Collision Avoidance Maneuvers Optimization in the Presence of Multiple Encounters}

\author{Zeno Pavanello\footnote{PhD candidate, Te P\=unaha \=Atea - Space Institute, 20 Symonds Street, Auckland Central, Auckland 1010, New Zealand; zpav176@aucklanduni.ac.nz (Corresponding Author).}, Laura Pirovano\footnote{Research Fellow, Te P\=unaha \=Atea - Space Institute, 20 Symonds Street, Auckland Central, Auckland 1010, New Zealand; laura.pirovano@auckland.ac.nz.} and Roberto Armellin\footnote{Professor, Te P\=unaha \=Atea - Space Institute, 20 Symonds Street, Auckland Central, Auckland 1010, New Zealand; roberto.armellin@auckland.ac.nz. Member AIAA.}}
\affil{Te P\=unaha \=Atea - The Space Institue, The University of Auckland, Auckland 1010, New Zealand}
\author{Andrea De Vittori \footnote{Ph.D. candidate, Department of Aerospace Science and Technology, Via La Masa 34, Milan 20156. Italy.} and 
Pierluigi di Lizia \footnote{Associate Professor, Department of Aerospace Science and Technology, Via La Masa 34, Milan 20156. Italy.}}
\affil{Polytechnic University of Milan, Milan 20156, Italy}

\begin{document}

\maketitle

\begin{abstract}
The optimization of fuel-optimal low-thrust \glspl{cam} in scenarios involving multiple encounters between spacecraft is addressed.
The optimization's objective is the minimization of the total fuel consumption while respecting constraints on the \acrlong{tpoc}. 
The solution methodology combines sequential convex programming, second-order cone programming, and differential algebra to approximate the non-convex optimal control problem progressively. 
A \acrlong{gmm} method is used to propagate the initial covariance matrix of the secondary spacecraft, allowing us to split it into multiple mixands that can be treated as different objects. This leads to an accurate propagation of the uncertainties.
No theoretical guarantee is given for the convergence of the method to the global optimum of the original \acrlong{ocp}. Nonetheless, good performance is demonstrated through case studies involving multiple short- and long-term encounters, showcasing the generation of fuel-efficient \glspl{cam} while respecting operational constraints.\blfootnote{A part of this work was presented at the 2024 AIAA Science and Technology Forum and Exposition, held in Orlando, Florida, January 8-12 2024. Paper AIAA-2024-0845 was entitled "A Convex Optimization Method for Multiple Encounters Collision Avoidance Maneuvers".}
\end{abstract}

\section{Introduction}
Spacecraft operations are encountering growing challenges due to the proliferation of space debris and the saturation of valuable orbital regimes, notably the \gls{geo} and the \gls{leo} ones, where large constellations reside \cite{ESA2023}. Ensuring the safety of spacecraft has become a paramount concern. In this regard, the \acrfull{cam} has emerged as a critical requirement for spacecraft operators, aiming to autonomously and efficiently navigate through potentially hazardous encounters. While the existing body of research extensively addresses single short-term encounters \cite{Bombardelli2015,Gonzalo2019,Hernando-Ayuso2020,Armellin2021}, recent attention has shifted toward scenarios involving multiple close encounters \cite{Kim2012,Lee2019,Frigm2020,Sanchez2021,Masson2022}. This paper presents an innovative approach to address the complex problem of fuel-efficient \glspl{cam} in the context of multiple short- and long-term encounters. In short-term encounters, the problem can be simplified, and each conjunction can be framed in its B-plane, whereas this is not possible in the long-term case because the path of relative motion inside the combined covariance ellipsoid is far from rectilinear, and the encounter is not instantaneous \cite{Patera2003Satellite}.

Conventional strategies for \gls{cam} optimization primarily focus on mitigating a single short-term encounter through the reduction of the \gls{poc} or the increase of the miss distance \cite{Bombardelli2015,Hernando-Ayuso2020,DeVittori2022Geo,DeVittori2023,DeVittori2023a,Dutta2022,Armellin2021}. Very few strategies have been developed for long-term encounters \cite{Mueller2009,Serra2015,Pavanello2023Long}. However, the increasing density of space debris and the prospect of multiple close encounters demand a more comprehensive and sophisticated approach. 
Lately, NASA is taking steps in this direction, with the introduction of the concept of \gls{tpoc}  \cite{Frigm2020} to quantify the contribution of multiple conjunctions to the overall collision probability: they recognize that in an environment that is densely populated, treating each conjunction as a separate event will become unsustainable for operators. The challenge lies in developing methodologies that can efficiently compute optimal trajectories while considering both onboard propulsion's practical constraints and the encounters' combined probabilistic nature. Regarding long-term encounters, few publications have addressed the problem of computing a \gls{cam} \cite{Mueller2009,Serra2015,Pavanello2023Long,Pavanello2023}. In this circumstance, instead of minimizing \gls{poc}, the alternative \gls{ipc} is controlled and kept below a certain threshold for the whole window of interest.

The multiple short-term problem has been faced in previous literature. Kim \textit{et al.} \cite{Kim2012} developed a genetic algorithm approach that deals with up to four conjunctions but must perform a separate $\dv$ impulse for each encounter. Their method does not include a check on \gls{tpoc}; instead, they minimize \gls{poc} of each conjunction separately. This results in a higher final \gls{tpoc}, and thus the method underestimates the risk of a collision. The authors do not provide a run-time for their method, but genetic algorithms are notoriously slower to converge with respect to other methods. 
Sánchez and Vasile \cite{Sanchez2021} developed a comprehensive tool that is based on a linear approximation of the dynamics with impulsive maneuvers and finds the optimal \gls{cam} using a min-max optimization. Their method combines machine learning and multi-criteria decision-making techniques to assess the collision risk and rank possible optimal \glspl{cam}. They propose two main strategies: either computing a single maneuver or treating each conjunction as a separate event. In the first case, only the worst-case scenario is considered; therefore, the computed impulse may not be optimal. In the second case, a separate maneuver is computed for each one of the encounters: this means that after the maneuver for the first conjunction has been computed, \gls{poc} must be re-calculated for each of the following conjunctions. This process must be repeated as many times as the number of the conjunctions. Sánchez and Vasile's approach can be very beneficial for operators on the ground to assist in the selection of the preferable \gls{cam} for a specific scenario. Still, it cannot be used to compute autonomous maneuvers due to the high computational load and the necessary human interaction that it requires.
Massòn \textit{et al.} \cite{Masson2022} cast the problem as a \gls{lp} with non-convex quadratic constraints, and they solve it by preallocating the maneuvering opportunities. They linearize the dynamics using the Yamanaka-Ankersen \gls{stm}, and they impose linear \gls{sk} constraints. Since they employ the $l^1$ norm of the control to keep the problem linear, they could recover a solution with higher propellant consumption if compared with the $l^2$ norm.
To summarize, the state-of-the-art still lacks a method that can deal with an arbitrary number of conjunctions comprehensively (using \gls{tpoc}) and does not need \textit{a-priori} information on the maneuvering time. Moreover, all of the above-referenced works considered Keplerian dynamics, disregarding the contributions of orbital perturbations. Finally, to the best of the authors' knowledge, the multiple long-term conjunction \gls{cam} optimization problem has never been addressed.

An aspect that is closely related to this topic is the nonlinear propagation of the uncertainty in the state of the objects involved in the conjunctions. When the time span of the scenario involves multiple orbits and the uncertainties of the secondary spacecraft are discretely large, the initial Gaussian approximation can lose validity \cite{Luo2017, Yanez2019}. We assume that accurate orbit determination is performed on the primary so that its uncertainty can be linearly propagated and remain Gaussian in the considered time window. Different methods can be used to recover an accurate description of the uncertainty of the secondary after a long propagation. Among these, the Monte Carlo analysis is the most accurate but most computationally expensive; unscented transforms reduce the computational load but assume Gaussian distribution, \gls{gmm} provide a good compromise between accuracy and computational load \cite{Vittaldev2015}. The research papers previously mentioned either do not consider the evolution of the uncertainty (because they use the instantaneous value of covariance given by a \gls{cdm}), or they only consider a linear propagation. In this regard, a first attempt at including the nonlinear propagation of the covariance matrix in the \gls{cam} optimization was taken by Dutta and Misra, who recently developed a pseudo-\gls{smd} approach to treat a single conjunction with \gls{gmm} \cite{Dutta2023}. This work will show how the approach developed for multiple conjunctions defined by \glspl{cdm} can be adapted to treat multiple conjunctions where the initial covariance is propagated nonlinearly using a \gls{gmm}. This allows for a more accurate depiction of the \gls{koz} when defining the convexified \gls{poc} constraint. 

The solution methodology is based on convex optimization methods and is tailored to address the intricacies of multiple encounters and nonlinear uncertainty propagation. This direct optimization approach is chosen for its proven capability in solving complex aerospace engineering problems efficiently \cite{Malyuta2021Advances,Li2022,Pirovano2024}.
These findings highlight the convenience of employing a \gls{socp} to formulate the convexified problem, as it allows for the minimization of the $l^2$ norm of the control history, which is a crucial requirement in spacecraft trajectory optimization problems. Turning the problem into an \gls{socp} involves convexifying all non-convex aspects of the original \gls{ocp}, notably the fuel-optimal objective function, the nonlinear dynamics constraint, and the \gls{poc} constraint \cite{Boyd2}. In the proposed approach, lossless relaxation is used to linearize the objective function. Furthermore, the nonlinear dynamics are automatically linearized using \gls{da}, which gives absolute freedom on the dynamics model considered. The convexification approach proposed in reference \cite{Armellin2021} is used to linearize the \gls{poc} constraint in the short-term problem, and the one proposed in reference \cite{Pavanello2023Long} is employed for the \gls{ipc} constraint in the long-term problem.

To enhance convergence, a trust region constraint, drawing from references \cite{Malyuta2021Tutorial,Losacco2023,Bernardini2023}, is employed to confine the solution space. However, due to the linearization of the dynamics, the optimizer cannot attain the optimal solution in a single \gls{socp} run. Hence, a \gls{scp} is constructed as a sequence of \glspl{socp}, where the solution from the previous iteration serves as the linearization point for the new iteration. This iterative process may necessitate a discrete number of iterations before converging to the optimal solution. Moreover, to account for the combined nature of \gls{tpoc} between the multiple encounters, a third outer iteration process is used to find the optimal maneuver. This process automatically adapts the single \gls{poc} limits associated with the conjunctions. Alternatively, a single \gls{scp} run is performed where, in place of the method proposed in \cite{Armellin2021} and \cite{Pavanello2023Long}, a direct linearization of the \gls{poc} or \gls{ipc} constraint is used.

In the subsequent sections, the paper will present the proposed methodology's mathematical framework and elaborate on the \gls{scp} approach for approximating the optimal control problem. Simulated case studies will demonstrate the approach's effectiveness in generating fuel-efficient \glspl{cam} for multiple short- and long-term encounters. By addressing this increasingly pertinent challenge, this paper aims to offer a valuable contribution to the scope of spacecraft \gls{ca} in cluttered orbital environments.

\section{Formulation of the Collision Avoidance Problem}
We consider multiple conjunctions between a primary spacecraft and $n_{conj}\in\mathbb{N}$ secondary objects. The primary is subject to a close encounter with each secondary at a certain \gls{tca}, indicated with $t_{CA,s}$, where $s\in\{1, ... \hspace{2pt},n_{conj}\}$. The primary and the secondary objects are characterized by a mean state and a covariance matrix at the starting time of the simulation $t_0$ or the conjunction's \gls{tca}, $t_{CA,s}$. The states are represented using any arbitrary set of elements, e.g. Cartesian, Keplerian, or generalized equinoctial elements

\begin{equation}
\begin{aligned}
      \quad &  \vec{x}_p(t_0)\sim \mathcal{N}(\vec{\zeta}_p(t_0),\boldsymbol{C}_p(t_0)), \quad &  \vec{x}_s(t_{CA,s})\sim \mathcal{N}(\vec{\zeta}_s(t_{CA,s}),\boldsymbol{C}_s(t_{CA,s})) \quad & \text{for } s \in \{1,... \hspace{2pt} ,n_{conj}\}
\end{aligned}
\end{equation}
where $\vec{x}_p(t_0)$ and $\vec{x}_s(t_{CA,s})\in\mathbb{R}^6$ are Gaussian \gls{mrv} variables, $\vec{\zeta}_p(t_0)$ and $\vec{\zeta}_s(t_{CA,s})\in\mathbb{R}^6$ are the mean values, $\boldsymbol{C}_p(t_0)$ and $\boldsymbol{C}_s(t_{CA,s})\in\mathbb{R}^{6\times6}$ are the covariance matrices.

The propagation of the state of the primary and the secondaries can be performed under a generic dynamical model
\begin{equation}
\begin{aligned}
    \quad & \dot{\vec{x}}_p(t) = \vec{f}_p(t,\vec{x}_p(t),\vec{u}(t),\vec{p}_p), \quad & \dot{\vec{x}}_s(t) =  \vec{f}_s(t,\vec{x}_s(t),\vec{p}_s) \quad & \text{for } s \in \{1,... \hspace{2pt} ,n_{conj}\},
\end{aligned}
\label{eq:diffeq}
\end{equation}
where $t\in\mathbb{R}_{[t_0,t_f]}$ is the independent time variable, $\vec{x}_p(t),\vec{x}_s(t),\dot{\vec{x}}_p(t)$ and $\dot{\vec{x}}_s(t)\in\mathbb{R}^6$ are the time-dependent states of the primary and the secondaries and their derivatives, $\vec{u}(t)\in\mathbb{R}^3$ is the control action, $\vec{p}_p\in\mathbb{R}^{m_p}$ and $\vec{p}_s\in\mathbb{R}^{m_s}$ are vectors of parameters, $\vec{f}_p(\cdot): \mathbb{R}_{[t_0, t_f]}\times\mathbb{R}^6\times\mathbb{R}^3\times\mathbb{R}^{m_p}\rightarrow\mathbb{R}^6$ and $\vec{f}_s(\cdot): \mathbb{R}_{[t_0, t_f]}\times\mathbb{R}^6\times\mathbb{R}^{m_s}\rightarrow\mathbb{R}^6$ are continuous functions.. 

To precisely establish the \gls{poc} constraint, it is crucial to explicitly express the relative position of the primary object w.r.t. to each secondary. If the state of the two objects is expressed in Cartesian \gls{eci} coordinates ($\Vec{x}_p(t)=[\vec{r}_p(t); \hspace{4pt} \Vec{v}_p(t)]$), as it is the case in \cref{sec:results} the relative position is the subtraction of the first three elements of the two normally distributed \gls{mrv}
\begin{equation}
  \begin{aligned}
  \quad & \Delta\vec{r}_s(t) = \vec{r}_p(t)-\vec{r}_s(t)
  \quad & s\in\{1, ... \hspace{2pt},n_{conj}\},
  \end{aligned}
\end{equation}
Given that the subtraction is a linear transformation, the relative position is also normally distributed
\begin{equation}
\begin{aligned}
      \quad & \Delta\vec{r}_s(t)\sim \mathcal{N}\big(\vec{\mu}_p(t)- \vec{\mu}_s(t),\boldsymbol{P}_p(t)+\boldsymbol{P}_s(t)\big) \quad & s\in\{1, ... \hspace{2pt},n_{conj}\},
    \label{eq:rel}
\end{aligned}
\end{equation}

The standard \gls{ocp} in the continuous domain for the short-term conjunction  \gls{cam} problem is stated as follows
\begin{subequations}
\begin{align}
\min_{u} \quad & J = \int_{t_0}^{t_f} u(t)\mathrm{d}t, \label{eq:ocpObj}\\
\mbox{s.t.} 
\quad &  \dot{\vec{x}}_p = \boldsymbol{f}(t,\vec{x}_p(t),\vec{u}(t)), \label{eq:ocpDyn}\\
\quad & \tpc(t) \leq \tpclim, \label{eq:ocpPc}\\
\quad & \vec{x}(t_0) = \vec{x}_0, \label{eq:ocpInit}\\
\quad & u(t) = \sqrt{u_1(t)^2+u_2(t)^2+u_3(t)^2}, \label{eq:ocpUEq} \\
\quad & u(t) \leq u_{max}, \label{eq:ocpUmax}
\end{align}
\label{eq:ocp}
\end{subequations}
where $\tpc\in\mathbb{R}(\cdot):\mathbb{R}\rightarrow\mathbb{R}$ is the \gls{tpoc}, which is monotonically increasing with time, and $\tpclim\in\mathbb{R}$ is the limit value imposed to the \gls{tpoc}.
In Problem \eqref{eq:ocp}, \cref{eq:ocpObj} is the fuel minimization objective function, \cref{eq:ocpDyn} is the dynamics constraint, \cref{eq:ocpPc} is the \gls{poc} constraint, \cref{eq:ocpInit} is the initial state bound, \cref{eq:ocpUEq} is an auxiliary constraint to define the norm of the control, and \cref{eq:ocpUmax} is the bound on the maximum value of the control action. The collision probability can include $n_{conj}$ conjunctions. Since \glspl{cam} typically involve a $\dv$ on the order of magnitude of [\si{mm/s}] or [\si{cm/s}], the mass loss due to the maneuver is not considered in the equations of motion and the optimization problem \cite{Armellin2021}.

\subsection{Multiple Short-Term Conjunctions}
Conjunctions can be categorized into two types: short-term and long-term. Short-term conjunctions, being more frequent, have received a more comprehensive investigation. In such instances, the relative velocity is notably high, resulting in near-instantaneous collisions. Consequently, the dynamics of the conjunction can be effectively approximated as linear without sacrificing accuracy, and the event is typically studied on the two-dimensional B-plane \cite{Alfriend2000,Armellin2021,Hernando-Ayuso2020,DeVittori2022}.

In a scenario where $n_{conj}$ consecutive conjunctions take place, the B-plane of a conjunction $s$ is centered on the secondary object, the $\eta$ axis is aligned with the direction of the relative velocity of the primary w.r.t. the secondary, and the $\xi\zeta$ plane is perpendicular to the relative velocity axis. Since \gls{tca} is the time when the miss distance is lowest, it follows that  $\Delta\vec{r}_s(t_{CA,s})\cdot\Delta\vec{v}_s(t_{CA,s}) = 0$ and thus the relative position lies on the $\xi\zeta$ plane. The following discussion assumes that all the uncertainty is concentrated around the secondary object, and all the mass is concentrated around the primary \cite{Alfriend2000,Zhang2020}. In this way, the primary's state, which has to be optimized, can be treated as a deterministic variable.

The \gls{poc} of a single encounter is given by the integration of the \gls{pdf} of the relative position over the circle $\mathbb{C}_\mathrm{HBR}$ defined by the combined \gls{hbr} of the two spacecraft and centered in the primary
\begin{equation}
P_{C,s}=\frac{1}{(2 \pi)\sqrt{\mathrm{det}\big(\boldsymbol{P}_\mathcal{B}(t_{CA,s})\big)}} \iint_{\mathbb{C}_\mathrm{HBR}} \mathrm{exp}\left(-\frac{\Delta\vec{r}_\mathcal{B}\transp(t_{CA,s})\vec{P}_\mathcal{B}^{-1}(t_{CA,s})\Delta\vec{r}_\mathcal{B}(t_{CA,s})}{2}\right) \mathrm{d} A,
\label{eq:pceq}
\end{equation}
where $\Delta\vec{r}_\mathcal{B}(t_{CA,s})\in\mathbb{R}^2$ $\vec{P}_\mathcal{B}(t_{CA,s})\in\mathbb{R}^{2\times2}$ are the projection of the relative position mean value and covariance on the encounter B-plane at \gls{tca}. The argument of the exponential function is the \gls{smd} at \gls{tca}, $d_{m,s}^2$, divided by $2$. Various approaches have been proposed to numerically solve this integral over the years; in this work, we will use Chan's method \cite{Chan2004International}. By employing a numerical inversion of Chan's formula \cite{DeVittori2022}, it is possible to get the \gls{smd} limit from a \gls{poc} limit for the conjunction ($\smdlim = f(\bar{P}_C)$). Chan's \gls{poc} formula is generally preferred over others because of its accuracy and computational efficienty \cite{Alfano2007Review}.

Assuming that the conjunctions are uncorrelated, the \gls{tpoc} of the multiple encounter scenario is computed as the complement of the probability of not colliding with any object
\begin{subequations}
    \begin{equation}
    \tpc = 1-\prod_{s=1}^{n_{conj}} \big(1-P_{C,s}\big),
    \label{eq:pcTot}
    \end{equation}
    \begin{equation}
    \tpc \approx\sum_{s=1}^{n_{conj}}P_{C,s}.
    \label{eq:pcTotApprox}
    \end{equation}
\end{subequations}

In the hypothesis of small probabilities, it can be shown that the product in \cref{eq:pcTot} can be replaced by the summation in \cref{eq:pcTotApprox}. We will always use \cref{eq:pcTot} when computing the \gls{tpoc}, but we will leverage the approximated formula \cref{eq:pcTotApprox} to compute first guesses for the limits values of the \gls{poc} of the single conjunctions.
One might assume that the hypothesis of uncorrelation is a strong one because the trajectory of the primary would be affected by a collision, thus changing the \gls{poc} of the subsequent conjunctions. Nonetheless, as Frigm \textit{et al.} highlight \cite{Frigm2020}, from an operation point of view, "there is no reason to suppose that a close approach with a particular secondary is likely to promote a close approach with some other secondary."

\subsection{Multiple Long-Term Conjunctions}
Patera \cite{Patera2003Satellite} demonstrates that when the relative trajectory of the primary spacecraft diverges significantly from a straight line within the combined covariance ellipsoid, the encounter qualifies as long-term, requiring the use of a distinct model for the computation of \gls{poc}. Despite numerous recent efforts to develop efficient methods for computing \gls{poc} in long-term encounters (references \cite{Coppola2012,Alfano2014,Xu2011Research,Wen2022}), the presence of nonlinear relative trajectories within the combined uncertainty ellipsoid makes the risk quantification and, consequently, the design of long-term \glspl{cam} very complex.
Quantitative thresholds to distinguish between the two types of conjunctions have been proposed. Among these, Chan's \cite{Chan2004Short} states that conjunctions are short-term if the path of relative motion may be considered a straight line over a distance of $8-25$ \si{km} with a deviation lower than $2$ \si{m}; Dolado \textit{et al.}'s \cite{Dolado2012Satellite}, instead, identify conjunctions as short-term if the relative velocity at \gls{tca} is higher than $10$ \si{m/s}. 
In recent literature, it has been shown that \gls{ipc} is a valuable alternative to \gls{poc} when addressing long term encounters \cite{NunezGarzon2022,Pavanello2023,Jones2013Satellite,Adurthi2015Conjugate,Zhang2020}: in general, an increase in \gls{poc} is accompanied by a high value of \gls{ipc}, so controlling the evolution of IPoC allows to indirectly control the growth of PoC. For this reason, in agreement with a good portion of the state-of-the-art, we use \gls{ipc} as a risk metric rather than the more complex \gls{poc} when dealing with long-term encounters.

The integral of the \gls{ipc} at any time $t\in\mathbb{R}_{[t_0,t_f]}$ must be performed over the the three spatial dimensions 
\begin{equation}
P_{IC,s}(t)=\frac{1}{(2 \pi)^{3/2}\mathrm{det}(\boldsymbol{P}(t))^{1/2}} \iiint_{\mathbb{S}_\mathrm{HBR}} \mathrm{exp}\left(-\frac{\Delta\vec{r}_s\transp(t)\vec{P}_{rel,s}^{-1}(t)\Delta\vec{r}_s(t)}{2}\right) \mathrm{d} V,
\label{eq:ipceq}
\end{equation}
where both the relative position mean value $\Delta\vec{r}_s(t)\in\mathbb{R}^3$ and its covariance $\vec{P}_{rel,s}(t)\in\mathbb{R}^{3\times3}$ evolve in time. Multiple approaches have been proposed to solve this integral. We will consider one analogous to Alfriend's \cite{Alfriend2000} for long-term encounters (its derivation can be found in \cite{Pavanello2023Long}). Under the assumption of this method, i.e., the \gls{pdf} of the relative position is considered constant inside the hard body sphere, \gls{ipc} has a direct dependence on \gls{smd}

\begin{equation}
   P_{IC,s}(t) = \sqrt{\frac{2}{\pi \mathrm{det}\big(\boldsymbol{P}(t)\big)}}\frac{R^3}{3}\mathrm{exp}\left(-\frac{d_{m,s}^2(t)}{2}\right).
   \label{eq:constipc}
\end{equation} 
In this case, the \gls{smd} is computed using the three-dimensional relative position variable.
Analogously to \cref{eq:pcTot}, the \gls{tipc} in a multiple long-term encounters scenario is given by the combination of the single uncorrelated conjunctions 

\begin{equation}
    \tipc(t) = 1-\prod_{s=1}^{n_{conj}} \big(1-P_{IC,s}(t)\big) \approx \sum_{s=1}^{n_{conj}} P_{IC,s}(t),
    \label{eq:ipcTot}
\end{equation}
Where the approximated equation is obtained using the same approximation of \cref{eq:pcTotApprox}.

\subsection{Gaussian Mixture Model for Uncertainty Propagation}
When dealing with multiple encounters with the same secondary or with a long-term conjunction, the time frame of interest can be very extended, spanning in the order of magnitude of the orbital periods. In the former case, the B-plane configuration depends on the propagated covariance, whereas in the latter, it is necessary to know the values of the elements of the covariance matrix of the relative position at every time to compute the \gls{ipc} from \cref{eq:constipc}. If the propagation period is sufficiently long, the nonlinearities of the dynamics propagation become predominant, and the Gaussianity of the state is lost.

An effective way to recover a precise propagation of the uncertainty, which avails in the computation of the collision metrics, is the use of a \gls{gmm}. Through this technique, the \gls{pdf} of the secondary spacecraft is split into $n_{mix}\in2\mathbb{N}+1$ mixands. Each mixand is characterized by a mean value and a covariance matrix derived from the original distribution. The propagations of the mixands happen independently, and they are handled as different secondaries.
The initial \gls{pdf} split of the secondary's state yields 
\begin{equation}
p[\vec{x}(t)]=\sum_{c=1}^{n_{mix}} \gamma_c p_g\left[\vec{x}(t) ; \vec{\mu}_c(t), \vec{P}_c(t)\right],
\end{equation}
where $c\in\{1, ... \hspace{2pt},n_{mix}\}$ indicates the single mixand, $p$ is the \gls{pdf} of the secondary object, $p_g$ indicates the Gaussian \gls{pdf}, and $\gamma_c\in\mathbb{R}$ are a set of weights such that $\sum_{c=1}^{n_{mix}} \gamma_c = 1$. It is assumed that most of the uncertainty is concentrated in the state of the secondary, so there is no need to split the state of the primary.

Vittaldev's uni-variate library method is used to perform the split \cite{Vittaldev2015}. Given the original \gls{pdf} $p[\vec{x}(t)]$, the uni-variate splitting library is applied along the direction $\hat{a}$
\begin{subequations}
    \begin{equation}
     \vec{\mu}_c=\vec{\mu}+\mu_c \boldsymbol{S} \hat{a}^\star,
    \end{equation}
    \begin{equation}
    \boldsymbol{P}_c=\boldsymbol{S}\left[\boldsymbol{I}+\left(\sigma_c^2-1\right) \hat{a}^\star \hat{a}^{\star T}\right] \boldsymbol{S}^T,
    \end{equation}
\end{subequations}
where the matrix $\boldsymbol{S}$ is the Cholesky decomposition of $\vec{P}$, $\mu_c$ and $\sigma_c$ are from the uni-variate splitting library and $\hat{a}^{\star}$ is

\begin{equation}
    \hat{a}^{\star}=\frac{S^{-1} \hat{a}}{\left\|S^{-1} \hat{a}\right\|_2}.
\end{equation}

Once the split has been performed, $n_{mix}$ fictitious secondary objects are defined at \gls{tca}, which can be treated independently in the \gls{ca} framework. Each of these objects encounters the primary spacecraft each of the conjunction times $n_{conj}$. The \gls{tpoc} or \gls{tipc}, then, is derived from \cref{eq:pcTot} or \cref{eq:ipcTot} with the introduction of a nested product function and the mixand weights

\begin{subequations}
\begin{equation}
   \tpc = 1- \prod_{s=1}^{n_{conj}} \prod_{c=1}^{n_{mix}}\left(1-\gamma_c P_{C,cs}\right))\approx\sum_{s=1}^{n_{conj}} \sum_{c=1}^{n_{mix}}\left(\gamma_c P_{C,cs}\right),
    \label{eq:pcGmm} \end{equation}
\begin{equation}
   \tipc(t) = 1- \prod_{s=1}^{n_{conj}} \prod_{c=1}^{n_{mix}} \left(1-\gamma_c P_{IC,cs}(t)\right)\approx\sum_{s=1}^{n_{conj}} \sum_{c=1}^{n_{mix}} \left(\gamma_c P_{IC,cs}(t)\right),
    \label{eq:ipcGmm} 
\end{equation}
\end{subequations}
where $P_{C,cs}$ and $P_{IC,cs}\in\mathbb{R}$ are the \gls{poc} and \gls{ipc} of the $c$-th mixand in the $s$-th conjunction. The approximated equations are obtained using the same approximation of \cref{eq:pcTotApprox}.
The splitting direction $\hat{a}$ is arbitrary, but it determines how good an approximation the \gls{gmm} provides w.r.t. the real distribution. Following works like Losacco's \cite{Losacco2023} and Vittaldev's \cite{Vittaldev2016Space}, a good selection of the splitting direction depends both on the \gls{nli} of the propagation and on the components of the initial covariance matrix. The \gls{nli}, indicated with $\vec{\nu}$, gives a quantitative indication of the effect of the nonlinearities in a certain direction; in reference \cite{Losacco2023}, how to compute it is effectively shown. A vector $\vec{\phi}\in\mathbb{R}^6$ made of the 2-norm of each column of the initial covariance matrix is used to quantify the entity of the uncertainty in each direction

\begin{equation}
    \begin{aligned}
        \phi_k = \sqrt{\sum_{l=1}^6 P_{lk}^2} \quad & k \in [1,6],
    \end{aligned}
\end{equation}
The splitting direction, then, considers the information from both the \gls{nli} and the covariance

\begin{equation}
    \hat{a} = \frac{\Vec{\nu}\circ\vec{\phi}}{||\Vec{\nu}||\cdot||\vec{\phi}||},
\end{equation}
where $\circ$ indicates the Hadamard product.

\section{Convex Formulation}

The \gls{ocp} Problem \ref{eq:ocp} is transformed into a \gls{scp} to be solved iteratively. The methods applied are described in detail in \cite{Pavanello2023}: the \gls{scp} iteratively approximates the global solution of the \gls{ocp} via a sequence of \glspl{socp}. Here, we recall the constraints and the objective function and adapt it to the multiple encounters problem.
Three sources of non-convexities are present in the original \gls{ocp} as formulated in Problem \ref{eq:ocp}:
\begin{enumerate}
    \item  orbital dynamics in \cref{eq:ocpDyn} are highly nonlinear, thus non-convex.
    \item \cref{eq:ocpPc} is a nonlinear constraint because it involves the use of \cref{eq:pcTot} and the inversion of Chan's formula for short-term encounters or \cref{eq:constipc} for long-term ones.
    \item \cref{eq:ocpObj} and \cref{eq:ocpUEq} yield a nonlinear objective function.
\end{enumerate}
These issues are addressed in the following sections, where the dynamics constraint, the \gls{poc} constraint, and the objective function are convexified using different techniques. A more detailed explanation of the convexification approach can be found in \cite{Armellin2021} and \cite{Pavanello2023}.

\subsection{Convexification of the dynamics}
\label{sec:convexDynamics}
The dynamics of the problem are first discretized and then automatically linearized using \gls{da}.

\subsubsection{Discretization of the Dynamics}
The continuous time variable $t\in\mathbb{R}_{[t_0, t_f]}$ is substituted by the discrete time variable $t_i\in\{t_0, t_1, ..., t_N\}$, where $N+1$ is the number of nodes of the discretization. Following \cref{eq:diffeq} and via the use of an integration scheme, like Runge-Kutta 7-8, one obtains the states of the two spacecraft at node $i+1$, which depend on the state and the control at node $i$:
\begin{subequations}
\begin{align}
    \quad & \vec{x}_{p,i+1} = \vec{f}_{p,i}(t_i,\vec{x}_{p,i},\vec{u}_i,\vec{p}_p) \quad & i \in \{0,... \hspace{2pt}, N-1\}, \quad &  \label{eq:dynamicsa} \\
    \quad & \vec{x}_{s,i+1} = \vec{f}_{s,i}(t_i,\vec{x}_{s,i},\vec{p}_s) \quad & i \in \{0,... \hspace{2pt}, N-1\} \quad & s \in \{1,... \hspace{2pt} ,n_{conj}\},  
\end{align}
\label{eq:dynamics} 
\end{subequations}
where $\vec{f}_{p,i}(\cdot): \mathbb{R}_{[t_0,t_f]}\times\mathbb{R}^6\times\mathbb{R}^3\times\mathbb{R}^{m_p}\rightarrow\mathbb{R}^6$  and $\vec{f}_{s,i}(\cdot): \mathbb{R}_{[t_0,t_f]}\times\mathbb{R}^6\times\mathbb{R}^{m_s}\rightarrow\mathbb{R}^6$  are the functions that describe the dynamics at node $i$, $\vec{x}_{p,i} = \vec{x}_p(t_i)$, $\vec{x}_{s,i} = \vec{x}_s(t_i)$ and $\vec{u}_{i} = \vec{u}(t_i)$.
The \gls{da} tool is used to introduce perturbations on the primary state ($\vec{x}_{p,i}+\delta\vec{x}_{p,i}$)
and acceleration ($\vec{u}_i+\delta\vec{u}_i$) at each node, effectively expressing \cref{eq:dynamicsa} through Taylor polynomials:
\begin{equation}
\begin{aligned}
\quad & \vec{x}_{p,i+1} = \mathcal{T}^q_{\vec{x}_{p,i+1}}(\vec{x}_{p,i},\vec{u}_i) \quad & i \in \{0,... \hspace{2pt}, N-1\},
\end{aligned}
\label{eq:perturbations}
\end{equation}
where in general the expression $\mathcal{T}^q_{y}(x)$ indicates the $q^\text{th}$-order Taylor expansion of the variable $y$ as a function of $x$,  around the expansion point $\tilde{x}$ in which the polynomial is computed. The reader can find a detailed explanation of the use of \gls{da} in \cite{Armellin2010}.

\subsubsection{Linearization of the Dynamics}
\label{sec:linearization}
The discretized dynamics are linearized around a reference trajectory, so that a linear constraint can be set in the \gls{socp}. This is an iterative process which can require numerous runs, referred to as major iterations and denoted by index  $j$.
Once the optimization problem of a major iteration $j-1$ is solved, the solution ($\mathbb{x}^{j-1}$, $\mathbb{u}^{j-1}$) becomes available, a column vector comprising the state and control at each node.
Employing \gls{da}, this solution serves as an expansion point for constructing linear dynamics maps for iteration $j$.
The continuity condition is enforced by requiring that the state after the propagation of node $i$ is equal to the state before the propagation of node $i+1$. 
The state and acceleration of each node are expanded around the reference points, which are the output of the previous major iteration $\tilde{\Vec{x}}_i^j=\vec{x}_i^{j-1}$ and $\tilde{\Vec{u}}_i^j=\vec{u}_i^{j-1}$
\begin{equation}
\begin{aligned}
   \quad & \vec{x}_{i}^j = \tilde{\Vec{x}}_i^j + \delta \vec{x}_i^j, \quad & \vec{u}_{i}^j = \tilde{\Vec{u}}_i^j + \delta \vec{u}_i^j \quad & i \in \{0,... \hspace{2pt}, N-1\},
\end{aligned}
\label{eq:exBefore}
\end{equation} 
The state at the subsequent node is obtained through the propagation of the first-order dynamics. The linear continuity constraint imposes that the linearly propagated state be equal to the state at the subsequent node

\begin{equation}
\begin{aligned}
    \quad & \vec{x}_{i+1}^j = \boldsymbol{A}_{i+1}^j\vec{x}_i^j  + \boldsymbol{B}_{i+1}^j
    \vec{u}_i^j + \vec{c}_i^j   \quad & i\in \{0,... \hspace{2pt}, N-1\},
\end{aligned}
\label{eq:forcedDyn}
\end{equation}
where $\boldsymbol{A}_{i+1}^j\in\mathbb{R}^{6\times6}$ is the \acrlong{stm}, $\boldsymbol{B}_{i+1}^j\in\mathbb{R}^{6\times3}$ is the control-state transition matrix, $\vec{c}_i^j = \bar{\vec{x}}_{i+1} - \boldsymbol{A}_{i+1}^j\tilde{\vec{x}}_i^j  - \boldsymbol{B}_{i+1}^j
\tilde{\vec{u}}_i^j$ is the residual of the linearization, and $\bar{\vec{x}}_{i+1} = \vec{f}_{p,i}(t_i,\tilde{\vec{x}}_i,\tilde{\vec{u}}_i)$ is the constant part of the propagation of the state.
The initial condition is fixed because the maneuver cannot alter it 
\begin{equation}
    \vec{x}_0^j = \vec{x}_0^0.
    \label{eq:initBound}
\end{equation}

\subsubsection{Lossless Relaxation of the Control Magnitude Constraint}
Equation \cref{eq:ocpUEq} is a non-convex equality constraint. Following the work from \cite{Wang2018}, we introduce a lossless relaxation to convexify the constraint: \cref{eq:ocpUEq} is transformed into an inequality constraint, and the control magnitude is added to the optimization vector. In this way, \cref{eq:ocpUEq} becomes a second-order cone constraint. The variable $u_i$ is now allowed to take values higher than the norm of the control that acts on the dynamics
\begin{equation}
     u_i \geq \sqrt{u_{i,1}^2 + u_{i,2}^2 + u_{i,3}^2} \quad  i \in \{0,... \hspace{2pt}, N-1\}.
     \label{eq:slackConstr}
\end{equation}
The discretized forms of \cref{eq:ocpObj} and \cref{eq:ocpUmax} become respectively
\begin{equation}
J = \sum_{i=0}^{N-1} u_i\cdot\Delta t_i, \label{eq:obj}
\end{equation}

\begin{equation}
\begin{aligned}
    0 \leq u_i \leq 1 \quad & i \in \{0,... \hspace{2pt}, N-1\},
    \label{eq:slackbound}
\end{aligned}
\end{equation}
where $\Delta t_i = t_{i+1}-t_i$ is the time for which the control is active from node $i$ to $i+1$.
This linearization of the objective function is lossless, meaning that the optimal solution for the convexified problem is also optimal for the original problem.

\subsubsection{Trust region approach}
Using a trust region constraint in formulating the \gls{socp} is greatly beneficial for the convergence of the algorithm \cite{Malyuta2021Tutorial}. In this work, the approach originally introduced in \cite{Bernardini2023} and adapted in \cite{Pavanello2023} is used, which is reliant on the definition of the \gls{nli} from \cite{Losacco2023}. The interested reader is invited to refer to these works for details on the derivation of the constraint. The gist of the approach is to consider the radius of the trust region of an optimization variable as inversely proportional to its nonlinearity. The resulting constraints are 
\begin{subequations}
    \begin{align}
        \quad &\vec{\xi}_{i} \circ [\vec{x}_{i}\transp; \hspace{2pt} \vec{u}_{i}\transp]\transp \preceq \vec{\xi}_{i} \circ [\tilde{\vec{x}}_{i}\transp; \hspace{2pt} \tilde{\vec{u}}_{i}\transp]\transp + \bar{\nu}\cdot\mathbf{1} 
        \quad & i\in\{1,... \hspace{2pt}, N\}, \\
        \quad &\vec{\xi}_{i} \circ[\vec{x}_{i}\transp; \hspace{2pt} \vec{u}_{i}\transp]\transp \succeq \vec{\xi}_{i}\circ[\tilde{\vec{x}}_{i}\transp; \hspace{2pt} \tilde{\vec{u}}_{i}\transp]\transp - \bar{\nu}\cdot\mathbf{1} 
        \quad & i\in\{1,... \hspace{2pt}, N\},
    \end{align}
    \label{eq:tr}
\end{subequations}
where $\vec{\xi}_i$ is the measure of the nonlinearity of the associated state vector in the second-order \gls{da} propagation of the dynamics \cite{Bernardini2023}; $\bar{\nu}$ is a user-defined value that imposes a limit to the maximum \gls{nli} available to the solution.

The introduction of the trust region constraint can cause artificial infeasibility, so virtual controls are added to the dynamics constraints, \cref{eq:forcedDyn}. The new constraints with virtual controls become
\begin{equation}
\begin{aligned}
     \quad & \vec{x}_{i+1}^j - \boldsymbol{A}_{i+1}\vec{x}_i^j  - \boldsymbol{B}_{i+1}^j \vec{u}_i^j + \Vec{\upsilon}_{i+1}^j - \vec{c}_i^j = 0   \quad & i\in \{0,... \hspace{2pt}, N-1\},
\label{vcConstraints}
\end{aligned}
\end{equation}
where $\Vec{\upsilon}_{i}^j\in \mathbb{R}^6$ is the virtual control vector for the node $i$ at major iteration $j$.

In the objective function, a term is added which is proportional to the entity of the virtual controls in order to minimize them. This term is multiplied by a weight $\kappa_{vc}\in\mathbb{R}$, which must be sufficiently high. 

\begin{subequations}
    \begin{align}
    \upsilon_i \geq ||(\upsilon_{ca,i}^j)|| \quad & i \in \{1,... \hspace{2pt}, N\}, \label{eq:vcCone}\\
     J_{vc} = \kappa_{vc}\sum_{i=0}^N \upsilon_i,  \label{eq:vcCost}
    \end{align}
\end{subequations}

\subsection{Convexifcation of the PoC Constraint for Short-Term Encounters}
\label{sec:caShort}
Two alternative approaches are proposed to convexify the \gls{poc} constraint. The first one is based on the linearization of the \gls{smd}, and it consists of setting a separate constraint for each conjunction; the second is based on a direct linearization of the \gls{tpoc} formula of \cref{eq:pcTot}. 

\subsubsection{Linearized SMD constraint}
\label{sec:projLin}
The \gls{poc} constraint for a single conjunction can be formulated by leveraging the methodology proposed in references \cite{Mao2021} and \cite{Armellin2021}. The reader is invited to refer to this last reference for an in-depth explanation of the process.
An iterative projection and linearization algorithm is utilized to convexify the nonlinear \gls{poc} constraint.
The iterations are nested inside the major iterations used to linearize the dynamics and take the name of \textit{minor iterations}, denoted by the symbol $k$. For each conjunction node, the projection convex sub-problem aims to find the point on the surface of an ellipse that is closest to the relative position $\Delta\vec{r}_s^{j,k-1}$ from the previous minor iteration (where $s\in\{1, ... \hspace{2pt},n_{conj}\}$). In the formulation of the convex sub-problem, the indices $s$ and $j$ are dropped since one of these problems is solved multiple times inside the same major iteration and for each conjunction.

Let the covariance and relative position be expressed in the B-plane reference frame $\mathcal{B}$.
It is possible to define a transformation matrix $\boldsymbol{V}\in\mathbb{R}^2$ that diagonalizes the covariance matrix. We call the reference frame in which the covariance matrix is diagonal $\mathcal{C}$.
A simple quadratic optimization problem is used to find the point on the ellipse that is closest to the reference trajectory point
\begin{subequations}
\begin{align}
\min_{\hat{z}_\mathcal{C}} \quad &  ||\vec{z}_\mathcal{C}^k - \Delta\vec{r}_{q,\mathcal{C}}^{k-1}|| \label{eq:objsub}\\
\textrm{s.t.} \quad & (\vec{z}_\mathcal{C}^k)\transp (\boldsymbol{P}_\mathcal{C})^{-1}\vec{z}_\mathcal{C}^k \leq \smdlim, \label{eq:consub}
\end{align}
\label{eq:probEllipsoid}
\end{subequations}
The objective \cref{eq:objsub} imposes the minimization of the Euclidean distance between the relative position of the previous iteration and the optimization variable $\vec{z}_\mathcal{C}$. \cref{eq:consub} is a relaxed condition on the optimization variable to be inside the ellipsoid. These relaxed conditions are lossless since the minimization of the objective guarantees that the optimized variable is always positioned on the ellipse's surface, maximizing the distance from its center.

Once $\vec{z}_\mathcal{C}$ is determined, the solution is transformed back into the original reference frame using the equation  
$ \vec{z}_\mathcal{B}^k = \vec{V}\transp\vec{z}_\mathcal{C}^k$.
After obtaining $\vec{z}_s^{j,k}=\vec{z}_\mathcal{B}^k$, a linearization of the \gls{smd} constraint is applied, effectively turning it into a \gls{koz} constraint. Specifically, we ensure that the new optimized relative position lies within the semi-plane defined by the line tangent to the ellipse on $\vec{z}_s^{j,k}$. The equation of the constraint is, then
\begin{equation}
\begin{aligned}
    \nabla(\smd)_s^{j,k}\Big|_{\vec{z}_s^{j,k}}\cdot (\Delta\vec{r}_i^{j,k}-\vec{z}_s^{j,k}) \geq 0 \quad & i \in \{1,... \hspace{2pt}, N\}.
    \label{eq:caConstr}
\end{aligned}
\end{equation}

\subsubsection{Linearized TPoC constraint}
\label{sec:linTpoc}
An alternative to the use of the linearized \gls{smd} constraint is the linearized \gls{tpoc} constraint. The formula of the \gls{tpoc} from \cref{eq:pcTot} is linearized using \gls{da} as a function of the relative position in each of the encounters
\begin{equation}
    \delta\tpc = \vec{\nabla} \delta \mathbb{r},
    \label{eq:deltaPoC}
\end{equation}
where $\vec{\nabla}$ is the gradient of the expression that gives the total \gls{poc} as a function of the relative position at each \gls{tca} and $\mathbb{r}$ is the stacked vector of the positions of the primary in the nodes in which the conjunctions happen, i.e., the \glspl{tca}
\begin{equation}
  \mathbb{r} = \begin{bmatrix}
      \vec{r}_1\transp & ... & \vec{r}_{n_{conj}}\transp.
  \end{bmatrix}\transp
\end{equation}
The \gls{poc} constraint, in this case, is a scalar one
\begin{equation}
    \tpc = \vec{\nabla}\mathbb{r}  + d \leq \tpclim,
\label{eq:pcconstr}
\end{equation}
where $d = \tilde{P}_{TC} - \vec{\nabla}\tilde{\mathbb{r}}$ is the residual of the linearization (the tilde indicates the values computed at the expansion point). Note that using the absolute or the relative position of the primary is equivalent because in \cref{eq:pcconstr}, the position of the secondary would be eliminated. A trust region constraint is used to favor the convergence of the optimization with this constraint.

The trust region method is analogous to the one used to impose a trust region constraint based on the dynamics in \cite{Pavanello2023}. \cref{eq:pcTot} is computed using a second-order \gls{da} expansion, which yields the gradient of the function with a first order-expansion
\begin{equation}
    \vec{\nabla} = \bar{\vec{\nabla}}+\delta\vec{\nabla},
\end{equation}
where $\bar{\vec{\nabla}}\in\mathbb{R}^{3n_{conj}}$ is the constant part of the gradient and $\delta\vec{\nabla}\in\mathbb{R}^{3n_{conj}}$ is its first-order expansion w.r.t. the position vector.
\begin{equation}
\begin{aligned}
    \quad & \nabla_i = \bar{\nabla}_i + \delta\nabla_i =  \bar\nabla_i + \sum_{j=1}^{3n_{conj}} a_{ij} \delta r_j
    \quad & i\in\{1,... \hspace{2pt} ,3 n_{conj}\}
\end{aligned}
\end{equation}
The non-linearity vector for the single relative position variable $r_j$ is computed as the ratio between the norm of the second-order terms and that of the first-order ones
\begin{equation}
    \nu_j = \frac{\sqrt{\sum_{i=1}^{3n_{conj}} a_{ij}^2}}{||\bar{\vec{\nabla}}||} |\delta r_j| = \xi_j|\delta r_j|.
    \label{eq:nli}
\end{equation}
It is possible to write two vector equations for each conjunction, which bound the position of the given node
\begin{subequations}
    \begin{align}
        \vec{\xi}_{k} \circ \Delta\vec{r}_k \preceq \vec{\xi}_{k} \circ \Delta\tilde{\vec{r}}_k + \vec{\nu}
        \quad & k\in\{1,... \hspace{2pt} ,n_{conj}\}, \\
        \vec{\xi}_{k} \circ\Delta\vec{r}_k \succeq \vec{\xi}_{k}\circ\Delta\tilde{\vec{r}}_k - \vec{\nu}
        \quad & k\in\{1,... \hspace{2pt} ,n_{conj}\},
    \end{align}
    \label{eq:trCa}
\end{subequations}
where $\vec{\xi}_k$ and $\vec{\nu}_k\in\mathbb{R}^3$ include the three components of the single encounter.

\subsubsection{Refinement of the solution}
\label{sec:limAdapt}

The \gls{poc} constraint of the original \gls{ocp} \cref{eq:ocpPc} bound the \gls{tpoc} to below a limit value $\tpclim$. This limit can be split, and a portion of it can be assigned to the single conjunctions

\begin{equation}
    \tpclim = 1-\prod_{s=1}^{n_{conj}}(1-\bar{P}_{C,s}).
    \label{eq:multLims}
\end{equation}
Theoretically, in \cref{eq:multLims}, the single limits $\bar{P}_{C,s}$ can be optimized to build constraints that allow for the minimization of the $\dv$. Under the assumption that a higher $\dv$ leads to a lower \gls{tpoc}, we need to find a method to recover a solution for which the final \gls{tpoc} gets as close as possible to the total limit.

\begin{figure}[tb!]
    \centering
    \input{Diagram}
    \caption{High-level flow of the algorithm.}
    \label{fig:flow}
\end{figure}
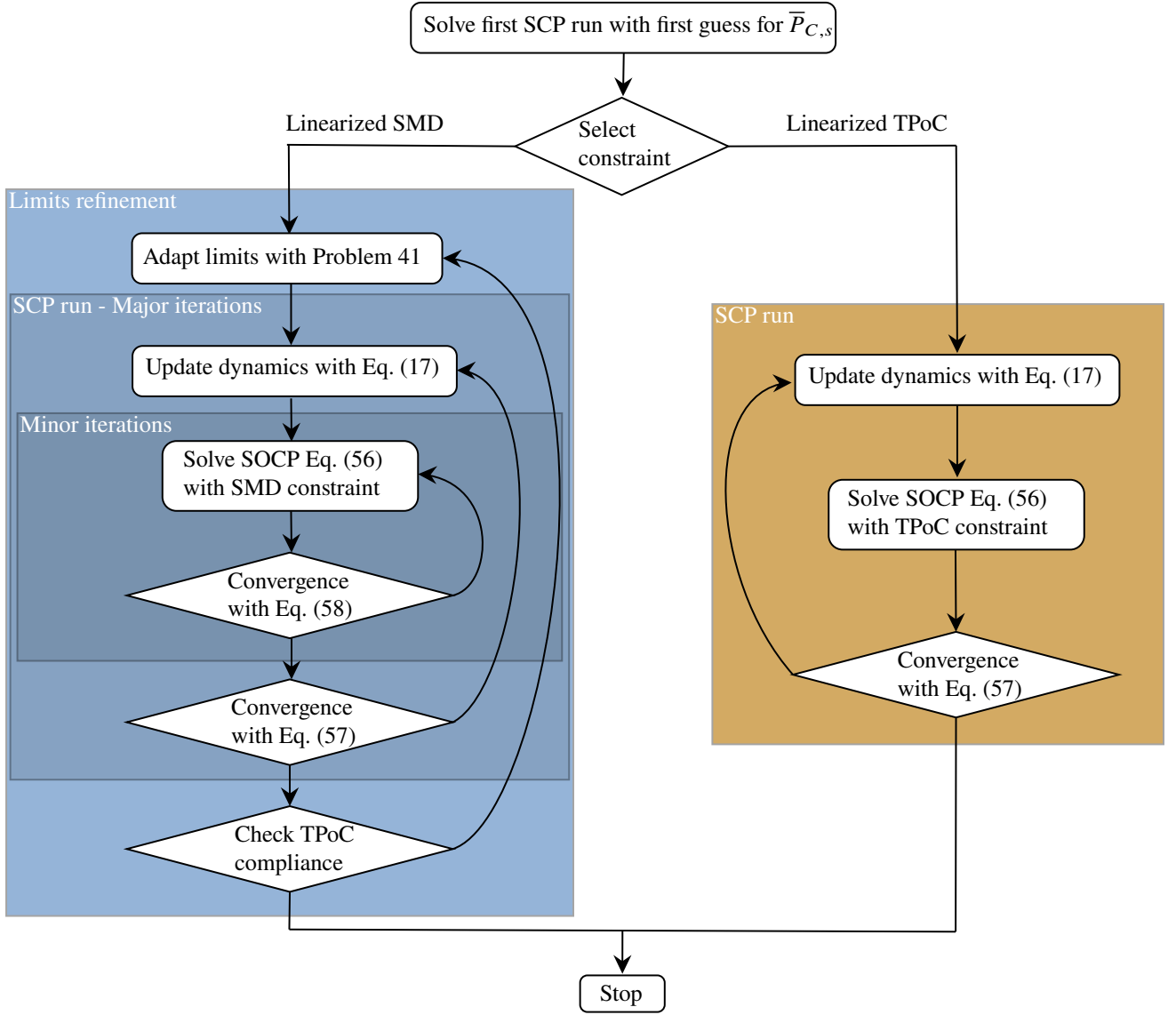

A first \gls{scp} is solved using a first guess for the limits of the single conjunctions. They are all set to the same value by leveraging the approximated \cref{eq:pcTotApprox}

\begin{equation}
    \begin{aligned}
        \quad & \bar{P}_{C,s} = \frac{\tpclim}{n_{conj}} \quad & s\in\{1, ... \hspace{2pt},n_{conj}\},
    \end{aligned}
\end{equation}
so that \cref{eq:multLims} is automatically satisfied.  This first convex run's total \gls{poc} is typically much lower than $\tpclim$, so a refinement is executed. We have two alternatives for the \gls{poc} constraint of this second \gls{scp}, which are reported schematically in \cref{fig:flow}. The first option is to include the linearized \gls{tpoc} constraint on the right side of \cref{fig:flow}. 
In fact, since it is based on a heavy linearization of the \gls{poc} formula, it is unlikely to converge when a good reference solution is not provided, and it can only be used to refine the solution obtained with the first convex run. An alternative approach, reported in the left side of \cref{fig:flow}, consists of adapting the limit \gls{poc} values $\bar{P}_{C,s}$ and running a new \gls{scp} with the new optimized values.

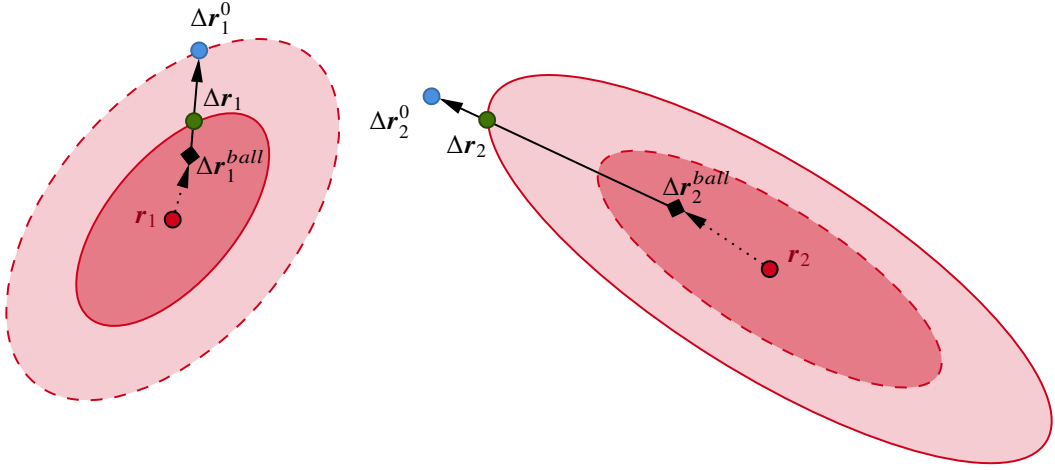
\begin{figure}[tb!]
    \centering
    \input{limAdaptScheme}
    \caption{Scheme of the adaptation of the \gls{poc} limits for two conjunctions.}
    \label{fig:schemeAdapt}
\end{figure}

To tune the limit of each conjunction to optimal values, a \gls{nlp} is solved. In \cref{fig:schemeAdapt}, the purpose of the \gls{nlp} is shown. The \gls{poc} limits imposed in the first \gls{scp} yield \gls{smd} limits that identify the dashed ellipses. Given these two limits, the convex solver may find a solution like the blue dots: the first conjunction dominates the scenario because $\vec{r}_1^0$ is on the border of the ellipse, while $\vec{r}_2^0$ is far from it. After the limits adaptation, more importance is given to the first conjunction. The limit \gls{poc} of the first conjunction is increased so that its \gls{smd} limit decreases, and the solution to the new \gls{scp} might be found closer to the ballistic miss distance. The minimization of the miss distance of the second conjunction is of no concern because, in the first run, the optimized point was found far from the limit. So, the \gls{smd} threshold of the second is allowed to grow (or the \gls{poc} threshold to shrink) until it reaches a value for which the optimized point is on the limit.

The optimization variables of the \gls{nlp} are collected in two vectors of parameters, $\Vec{\alpha}$ and $\Vec{\rho}\in\mathbb{R}^{n_{conj}}$. $\Vec{\alpha}$ is such that for the $s$-th conjunction $\bar{P}_{C,s} = \tpclim\cdot10^{-\alpha_s}$. We choose to introduce this parameter as an exponent to favor the algorithm's convergence since the limit \gls{poc} values can span a wide range of orders of magnitude; the minus sign of the exponent is useful to keep the optimization variables always positive.
Our working hypothesis affirms that getting closer to the ballistic relative position would result in a lower $\dv$. So, the optimized relative position is bound to lie in the same direction as the point resulted from the previous \gls{scp} w.r.t. the ballistic relative position.
\begin{equation}
    \Delta\hat{r}_s = \frac{\Delta\Vec{r}^0_s-\Delta\vec{r}_s^{ball}}{||\Delta\Vec{r}^0_s-\Delta\vec{r}_s^{ball}||},
\end{equation}
where $\Delta\vec{r}_s^0$ and $\Delta\vec{r}_s^{ball}\in \mathbb{R}^3$ are the relative position of the previous \gls{scp} run and of the ballistic trajectory w.r.t. the secondary $s$. In this way, we only need to optimize the magnitude of the difference between the optimized relative position and the ballistic one ($\rho_s$) for each conjunction $s$
\begin{equation}
\begin{aligned}
    \quad & \rho_s = ||\Delta \vec{r}_s-\Delta \vec{r}_s^{ball}|| \quad & s\in\{1, ... \hspace{2pt},n_{conj}\}.
\end{aligned}
\end{equation}
$\Vec{\rho}$ collects the $\rho_s$ for every conjunction.
Now we can write the \gls{nlp} that is used to optimize the limit probabilities. In Problem \ref{prob:nlpWeights}, the variables with superscript $0$ indicate the output of the previous \gls{scp}

\begin{subequations}
\begin{align}
    \min_{\vec{\rho},\vec{\alpha}} \quad & \vec{\rho} \cdot \vec{\beta} \quad & \label{eq:nlpObj}\\
    \text{s.t.} \quad & \Delta\vec{r}_s = \Delta\vec{r}_s^{ball} + \rho_s\Delta\hat{r}_s  \quad & s\in\{1, ... \hspace{2pt},n_{conj}\}, \label{eq:nlpR} \\
    \quad & \bar{P}_{C,s} = \tpclim\cdot 10^{-\alpha_s} \quad & s\in\{1, ... \hspace{2pt},n_{conj}\}, \label{eq:nlpAlpha}\\
    \quad & P_{C,s}(\Delta\vec{r}_s,\boldsymbol{P}_s,HBR_s) \leq \bar{P}_{C,s} \quad & s\in\{1, ... \hspace{2pt},n_{conj}\}, \label{eq:nlpPoC} \\
    \quad & 1 - \tpclim - \prod_{s=1}^{n_{conj}} (1-\bar{P}_{C,s}) = 0, \quad & \label{eq:nlpPoCTot} \\
    \quad & \beta_s = \begin{cases} 1 & \mbox{if } |P^0_{C,s}-\bar{P}^0_{C,s}|/\bar{P}^0_{C,s} \leq \varepsilon  \\ 0 & \mbox{otherwise} \end{cases} \quad & s\in\{1, ... \hspace{2pt},n_{conj}\}, \label{eq:nlpBeta}\\
    \quad & 0 \leq \rho_s \leq \rho_s^0 & s\in\{1, ... \hspace{2pt},n_{conj}\},\label{eq:nlplimsRho} \\
    \quad & 0 \leq \alpha_s \leq 10 \quad & s\in\{1, ... \hspace{2pt},n_{conj}\}. \label{eq:nlplimsAlpha}  
    \end{align}
\label{prob:nlpWeights}
\end{subequations}
The objective function in \cref{eq:nlpObj} is the norm of the scaled distance of the optimized point from the ballistic one.
With \cref{eq:nlpR}, the optimized point is bound to lie in the same direction as the point resulted from the previous \gls{scp} w.r.t. the ballistic relative position.
\cref{eq:nlpAlpha} defines the \gls{poc} limits as functions of the optimization variable $\vec{\alpha}$. \cref{eq:nlpPoC} sets the collision avoidance constraint for each secondary. \cref{eq:nlpPoCTot} imposes that the combination of the \gls{poc} limits must give the \gls{tpoc} limit.
In \cref{eq:nlpBeta}, the vector $\vec{\beta}\in\mathbb{R}^{n_{conj}}$ is used to ensure that the conjunctions that affect the objective function are the ones for which the \gls{poc} contributes in a non-negligible manner to the \gls{tpoc}; $\varepsilon\in\mathbb{R}$ is a tolerance value, e.g., $10^{-2}$.
The last two constraints, \cref{eq:nlplimsRho} and \cref{eq:nlplimsAlpha}, are bound constraints on the optimization variables.

As a solution of Problem \ref{prob:nlpWeights}, one obtains the weights in $\vec{\alpha}$ that are used to find the limit \gls{poc} for each conjunction. Using this objective function requires the underlying hypothesis that a solution to the convex problem that is closer to the original ballistic trajectory would require a lower $\dv$, which might not always be the case.
Lastly, it is important to feed the \gls{nlp} a suitable first guess in terms of $\vec{\alpha}_0$ and $\Vec{\rho}_0$: the value of $\Vec{\rho}_0$ is taken from the result of the previous \gls{scp}, and $\vec{\alpha}_0$ is derived from the values of the \gls{poc} associated with the single conjunctions

\begin{equation}
    \begin{aligned}
        \quad & \alpha_{0,s} = -\log\left(\frac{P_{C,s}^0}{\tpclim}\right) \quad & s\in\{1, ... \hspace{2pt},n_{conj}\}.
    \end{aligned}
\end{equation}



\subsection{Conjunction representation with GMM}
\label{sec:gmm}
\begin{figure}[tb!]
    \centering
    \subfloat[Head-to-head encounter.]{    \input{HeadToHeadEnc}\label{fig:head2head}}
    \subfloat[Perpendicular encounter.]{\input{PerpendicularEnc}\label{fig:perpendicular}}
    \caption{Different geometric configurations of the conjunctions depending on the direction of relative motion.} 
    \label{fig:gmmConf}
\end{figure}
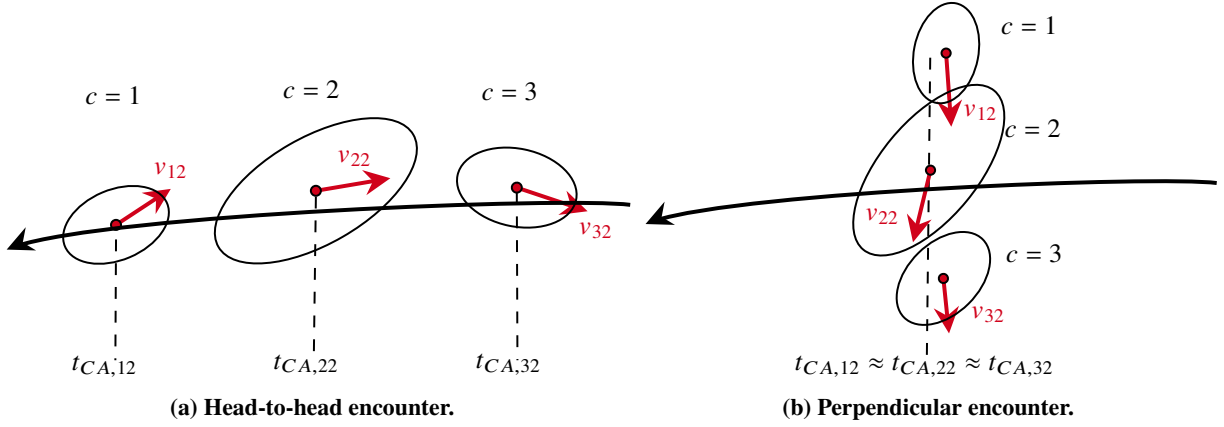

\begin{figure}[tb!]\centering
\def\svgwidth{\textwidth}
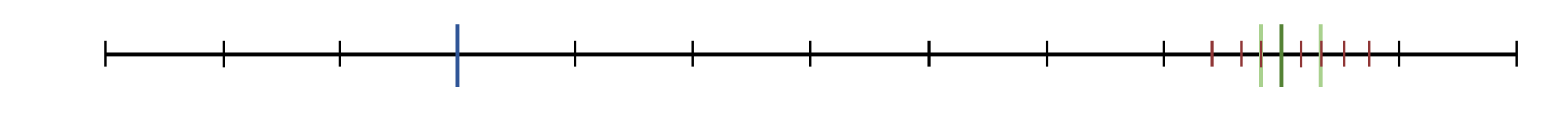
\caption{Grid adjustment to account for different \glspl{tca} of the mixands.}
\label{fig:grid}
\end{figure}

The \gls{gmm} representation of the uncertainty can be used in the framework of the multiple encounters with just a few adjustments. It is easy to visualize every mixand as a separate secondary, thus treating the conjunctions separately, like in the scenario where multiple encounters happen with different objects. As shown in \cref{fig:gmmConf}, depending on the direction of the split and on the direction of relative motion between the two bodies, the \glspl{tca} of the different mixands can be very close, if not identical. This is equivalent to the primary being required to avoid the collision with more than two very close objects at times that are close to one another. The issue arises of the previously established time grid (on which the dynamics have been discretized) not being able to catch the exact \glspl{tca}. This calls for refining the grid near the nominal \gls{tca}. \cref{fig:grid} depicts the idea of the grid refinement for a case with three mixands: the original equally spaced time grid is represented in black. While for the first conjunction, in this particular case, the three mixands all arrive at the close approach at the same time, for the second conjunction, they have had the time to drift from each other; thus, the close encounter with the primary happens at the nominal time for the central mixand (dark green line) and at slightly different times for the other two (light green lines). A grid thickening algorithm is used to be able to catch the offset \glspl{tca}, yielding new nodes like the ones represented in red.

First of all, the real \gls{tca} must be computed. One could do this by propagating the trajectory from node $t_{CA,s}$ previous to the \gls{tca} one and the one subsequent to it and finding the time for which the minimum separation is reached. This method, though, is highly inefficient. We leverage \gls{da} to obtain the same result. For a conjunction $s\in\{1, ... \hspace{2pt},n_{conj}\}$, let us take the nominal states of the two spacecraft at \gls{tca} and expand the secondary state
\begin{equation}
\vec{x}_{s,CA}^j = \Tilde{\vec{x}}_{s,CA}^j + \delta \vec{x}_{s,CA}^j,
\end{equation}
We propagate the states of the two satellites by a \gls{da} variable $\delta t_s$, obtaining, with a Taylor polynomial notation
\begin{equation}
    \begin{aligned}
     \quad & \vec{x}_{s,CA}^{j*} = \mathcal{T}^1_{\vec{x}_{s,CA}^j}(\delta t_s,\delta\vec{x}_{s,CA}^j), \quad & \vec{x}_{p,CA}^{j*} = \mathcal{T}^1_{\vec{x}_{p,CA}^{j*}}(\delta t_s),
    \end{aligned}
\end{equation}
where the stars indicate the variables at the real \gls{tca}. The relative state is a new Taylor polynomial dependent on $\vec{x}_{s,CA}^j$, $\tilde{\vec{x}}_{p,CA}^j$ and $\delta t_s$: $\Delta\vec{x}_{s,CA}^{j*} = \vec{x}_{p,CA}^{j*} - \vec{x}_{s,CA}^{j*}$. It can be split into the position and velocity parts to write the equation of the minimum distance. As mentioned in \cref{sec:caShort}, the relative position is minimum when the relative velocity is orthogonal to it

\begin{equation}
    \Delta\vec{r}_{CA}^*\cdot\Delta\vec{v}_{CA}^* = \mathcal{T}^1_{\Delta\vec{r}_{CA}^*\cdot\Delta\vec{v}_{CA}^*}(\delta t_s,\delta\vec{x}_{s,CA}^j) = 0.
\end{equation}
This parametric implicit equation can be solved for $\delta t_s$ using a polynomial partial inversion technique, as the one used in \cite{Armellin2010} and \cite{He2018}. The equation for $\delta t_s$ reads

\begin{equation}
        \delta t_s = \mathcal{T}^1_{\delta t_s}(\Delta\vec{r}_{CA}^*\cdot\Delta\vec{v}_{CA}^*,\delta\vec{x}_{s,CA}^j).
\end{equation}
The exact \gls{tca} of the single mixand, then, is

\begin{equation}
    t_{CA,s}^*= t_{CA,s} + \delta t_s.
\end{equation}

\noindent According to the new \glspl{tca}, nodes are added to the time grid, and the states of the primary and the relevant mixands are propagated to such times.

If a \gls{gmm} model is used to propagate the uncertainty, we can treat each mixand as an independent secondary object. This way, we can navigate back to the algorithm described in the previous section. If we want to adapt the limits of the mixands to refine the solution, we can use Problem \ref{prob:nlpWeights}, with two slight differences. Firstly, in the first guess, we use the approximated version of \cref{eq:pcGmm} and the limits are set proportional to the weights 

\begin{equation}
    \begin{aligned}
    \quad & \bar{P}_{C,cs} = \frac{\gamma_c}{n_{conj}} \tpclim\quad & c\in\{1, ... \hspace{2pt},n_{mix}\} \quad & s\in\{1, ... \hspace{2pt},n_{conj}\},
    \end{aligned}
    \label{eq:pclimgmm}
\end{equation}
where $n_{mix}$ is the number of \gls{gmm} mixands, $s$ indicates the conjunction and $c$ indicates the mixand; this is not a necessary condition since any combination of thresholds is valid as long as it respects the total probability threshold, but it is found to favor convergence in practice.\footnote{For example, the equation $\bar{P}_{C,cs} = \tpclim/(n_{mix}\cdot n_{conj})$ would be equally valid.}
Secondly, in \cref{eq:nlpPoC}, the \gls{poc} of the conjunction obtained with Chan's formulas must be multiplied by its weight
\begin{equation}
    \begin{aligned}
    \quad & \gamma_c P_{C,cs}(\vec{r}_{cs},\vec{C}_{cs},HBR_s)  \leq \bar{P}_{C,cs} \quad & c\in\{1, ... \hspace{2pt},n_{mix}\} \quad & s\in\{1, ... \hspace{2pt},n_{conj}\},
    \end{aligned}
\end{equation}
and in \cref{eq:nlpPoCTot}, the nested product must be included to account for the \gls{gmm} components in each conjunction

\begin{equation}
    1 - \tpclim - \prod_{s=1}^{n_{conj}}\prod_{c=1}^{n_{mix}} (1-\bar{P}_{C,cs}) = 0.
    \label{eq:psd}
\end{equation}
If \cref{eq:pclimgmm} is substituted into \cref{eq:psd} we see that the total limit is respected by the initial arbitrary split.

One last remark about the optimization with the \gls{gmm} method regards the computation of the \gls{smd} limit. To obtain a realistic \gls{smd} limit from the \gls{poc} limit with the inversion of Chan's formula, the mixand weight must be taken into account

\begin{equation}
    \begin{aligned}
    \quad & \bar{d}_{m,cs}^2 = f\left(\frac{\bar{P}_{C,cs}}{\gamma_c}\right) \quad & c\in\{1, ... \hspace{2pt},n_{mix}\} \quad & s\in\{1, ... \hspace{2pt},n_{conj}\}.
    \end{aligned}
\end{equation}

\subsection{Convexifcation of the PoC Constraint for Long-Term Encounters}
\label{sec:caLong}
When dealing with the long-term encounter problem, the adopted approach is derived from the one described in the previous \cref{sec:caShort}, with some substantial differences presented in this section.
A detailed description of the projection and linearization algorithm for the single conjunction in long-term encounters is given in reference \cite{Pavanello2023}. Here, it suffices to say that it is analogous to the one presented in \cref{sec:projLin}, but instead of working in the B-plane, we set the problem in the 3D \gls{eci} reference frame. There is no particular need to use an \gls{eci} frame, but it is convenient because it does not require transformations with respect to the frame used for the propagation of the dynamics.

\subsubsection{Linearized TIPoC Constraint}
\label{sec:tipcLin}
A similar approach to \cref{sec:linTpoc} is used for the \gls{tipc} constraint in the case of long-term encounters. In this case, the constraint has to be applied indistinctly at every node of the optimization, so the stack of positions changes node-wise $\mathbb{r}\in\mathbb{R}^{3n_{conj}\times N}$ and so does the Jacobian $\boldsymbol{J}\in\mathbb{R}^{3n_{conj}\times N}$. The node-wise linearized \gls{tipc} constraint reads
\begin{equation}
\begin{aligned}
    \quad & P_{TIC,i} = \vec{\nabla}_i\mathbb{r}_i  + d_i \leq \tipclim, \quad & i=\{1,... \hspace{2pt}, N\}
\end{aligned}
\label{eq:ipcconstr}
\end{equation}
where $\vec{\nabla}_i\in\mathbb{R}^{3n_{conj}}$ is the gradient of the expression that gives the \gls{tipc} as a function of the relative position w.r.t. each secondary and $d = \tilde{P}_{TIC,i} - \vec{\nabla}_i\tilde{\mathbb{r}}_i$ is the residual of the linearization.
The same trust region method from \cref{sec:linTpoc} is employed, bearing in mind that, in this case, the trust region must be applied to the relative position variables at each node of the optimization.

\subsubsection{Refinement of the solution}
When dealing with multiple secondaries or propagating \gls{gmm} components, also in long-term encounters, it is important to refine the solution with one of the two methods introduced in \cref{sec:limAdapt}. The first one, based on the linearization of the \gls{tipc}, was shown in \cref{sec:tipcLin}. The second one consists of adapting the \gls{ipc} threshold of each component to obtain the minimum $\dv$ maneuver. An analogous \gls{nlp} to the one presented in Problem \ref{prob:nlpWeights} is used to tune the weights of the limits after the problem is solved with a first guess for the limits. The first guess, as in the short-term scenario, is for the \gls{ipc} limits to be proportional to the mixand weights (using the approximated version of \cref{eq:ipcGmm})
\begin{equation}
\begin{aligned}
    \quad & \bar{P}_{IC,cs} = \frac{\gamma_c}{n_{conj}}\tipclim \quad & c\in\{1, ... \hspace{2pt},n_{mix}\} \quad & s\in\{1, ... \hspace{2pt},n_{conj}\}.
\end{aligned}
\end{equation}
Only the optimization node where the \gls{tipc} is maximum is used in the \gls{nlp} because the maximum point in the profile is expected to stay the same or be shifted by a short time. This single value of \gls{tipc} is indicated as $\oldhat{P}_{TIC} = \max_i(P_{TIC_i})$, where $i\in[0,N]$ and the correspondent values of \gls{ipc} for the single secondaries are $\oldhat{P}_{IC_s}$

\begin{subequations}
\begin{align}
    \min_{\rho,\vec{\alpha}} \quad & |\rho| \quad &  \label{eq:nlpObjLT}\\
    \text{s.t.} \quad & \vec{r} = \vec{r}^{ball} + \rho\hat{r},  \quad & \label{eq:nlpRLT} \\
    \quad & \bar{P}_{IC,cs} = \tipclim\cdot 10^{-\alpha_{cs}} \quad & c\in\{1, ... \hspace{2pt},n_{mix}\} \quad & s\in\{1, ... \hspace{2pt},n_{conj}\}, \label{eq:nlpAlphaLT}\\
    \quad & \gamma_c \oldhat{P}_{IC,cs}(\Delta\vec{r}_{cs},\boldsymbol{P}_{cs},HBR_s)\leq \bar{P}_{IC,cs} \quad & c\in\{1, ... \hspace{2pt},n_{mix}\} \quad & s\in\{1, ... \hspace{2pt},n_{conj}\}, \label{eq:nlpIpc} \\
    \quad & 1 - \tipclim -  \prod_{s=1}^{n_{conj}} \prod_{c=1}^{n_{mix}} (1-\bar{P}_{IC,cs}) = 0, \quad &  \quad & \label{eq:nlpPoCTotLT} \\
    \quad & 0 \leq \rho \leq \rho^0, & \quad &\label{eq:nlplimsRhoLT} \\
    \quad & 0 \leq \alpha_{cs} \leq 10 \quad & c\in\{1, ... \hspace{2pt},n_{mix}\} \quad & s\in\{1, ... \hspace{2pt},n_{conj}\}, \label{eq:nlplimsAlphaLT}  
    \end{align}
\label{prob:nlpWeightsLT}
\end{subequations}
where the position variables $\vec{r},\Vec{r}^{ball},\Delta\vec{r}_{cs}\in\mathbb{R}^3$ are all computed at the time step of maximum \gls{tipc}. 
In this case, the objective function in \cref{eq:nlpObjLT} is the minimization of the scalar variable $\rho$. The remaining equations from \cref{eq:nlpRLT} to \cref{eq:nlplimsAlphaLT} can be interpreted in the same way as those in Problem \ref{prob:nlpWeights}. The same considerations from \cref{sec:gmm} are also valid for the long-term encounters when propagating with a \gls{gmm}.

\subsection{Final Form of the Problem}
This section summarizes the problem we face and provides some insights into the optimization process. 
First, we finalize the low-level optimization problem, the \gls{socp}.
The complete objective function is given by the sum of \cref{eq:obj} and \cref{eq:vcCost}:
\begin{equation}
    J = \kappa_{vc}\sum_{i=0}^N \upsilon_i + \sum_{i=0}^{N-1} u_i\cdot\Delta t_i,
    \label{eq:totObj}
\end{equation}
where the node-wise control action is weighted to account for variable time steps.
The final \gls{socp} is reported in Problem \eqref{eq:optFinal}.
\begin{equation}
\begin{aligned}
\min_{\mathbb{x},\mathbb{u}} \quad & \text{\cref{eq:totObj}}\\
\mbox{s.t.} 
\quad & \text{\cref{eq:dynamics}, \cref{eq:caConstr} or \cref{eq:pcconstr} with \cref{eq:trCa}, \cref{eq:initBound,eq:slackConstr,eq:slackbound,eq:tr,eq:vcCone}}
\end{aligned}
\label{eq:optFinal}
\end{equation}


The major iterations' convergence criterion is computed as the difference between the maximum control acceleration of two consecutive iterations
\begin{equation}
    e_M = \max_i(||\vec{u}_i^j-\vec{u}_i^{j-1}||).
    \label{eq:majConv}
\end{equation}
Likewise, the minor iterations' criterion (for the linearized \gls{smd} method) is the difference between the relative position vectors of two consecutive minor iterations
\begin{equation}
    e_m = \max_i(||\vec{r}_i^j-\vec{r}_i^{j-1}||).
    \label{eq:minConv}
\end{equation}

It is relevant to remember that the successive linearization process is performed on the output of the previous optimization problem, so there is no online check on the accuracy of the forward-propagated solution. Such a check is performed after the optimization has reached convergence in all its levels. A validation error is computed as the maximum difference between the states at the last output of the optimization and the states obtained by forward-propagating the trajectory using the converged control profile
\begin{equation}
    e_{validation} = \max_i(||\vec{x}_i^{j=end}-\vec{x}_i^{val}||).
    \label{eq:validationError}
\end{equation}

To conclude the discussion, the high-level work-flow of the algorithm, as it has been presented in the previous sections, is shown in \cref{fig:flow}.

\section{Results}
\label{sec:results}
Realistic test cases for multiple encounters in \gls{leo} are analyzed. Three scenarios are considered. In the first one, the primary must deal with multiple short-term conjunctions with different secondaries; in the second one, the primary encounters the secondary more than once, and the uncertainties are propagated using \gls{gmm}; in the last one, a long term encounter with a single secondary is studied using \glspl{gmm}. Results from both the refinement methods from \cref{sec:caShort} and \cref{sec:caLong} are compared.
In all the scenarios, the discretization time step is such that the orbit of the primary is divided into $60$ nodes. The threshold \gls{tpoc} and \gls{tipc} are both set to $10^{-6}$. The tolerance on the major iterations is set to $10^{-3}$, and the one on the minor iterations for the linearized \gls{smd} method is set to $10^{-6}$.
The simulations are run with MATLAB r2022b on AMD Ryzen 9 6900HS @ 3.3GHz. The optimization is performed using MOSEK 10.0.24, which implements a state-of-the-art primal-dual interior-point solver. 
Scaled variables are used both during the integration of the equations of motion and in the optimization to favor numerical stability and convergence: the scaling constants are $L_{sc} = a_p$, $V_{sc} = \sqrt{\mu_E/a_p}$, $T_{sc} = \sqrt{a_P^3/\mu_E}$, and $A_{sc} = \mu_E/a_P^2$, where $a_P$ and $\mu_E$ are the semi-major axis of the primary's orbit and Earth's gravitational constant. In all the considered scenarios, the trajectories are discretized using a constant time step into 60 nodes per orbit.

\begin{table}[b!]
\centering
\caption{Case 1: conjunctions' parameters.}
\label{tab:conj1}
    \begin{tabular}{l|lllllll}
    \Xhline{4\arrayrulewidth}
     & $t_{CA}$ & $\Delta\vec{r}$ [\si{m}]  & $\Delta\vec{v}$ [\si{km/s}] & $d_{miss}$ [\si{m}] & $\pc$ [-] & $\gamma_v$ [deg] \\
    \Xhline{3\arrayrulewidth}
    \textbf{CDM 1}  &  1 \si{h} 34  \si{min} 3  \si{s}& [-11.45;\hspace{2pt} -19.44 \hspace{2pt}  7.44]  & [ 0.37; \hspace{2pt} -2.85; \hspace{2pt}  4.32] & 23.8  & 0.0017 & 40   \\ 
    \textbf{CDM 2}  &  2  \si{h} 4  \si{min} 20 \si{s}& [ 99.90;\hspace{2pt}  14.60 \hspace{2pt}  2.78]  & [-3.94; \hspace{2pt} -5.85; \hspace{2pt}  3.36] & 101.0 & 0.0018 & 62  \\ 
    \textbf{CDM 3}  &  3  \si{h} 43 \si{min} 11 \si{s}& [ 45.70;\hspace{2pt} -3.73  \hspace{2pt}  28.29] & [-0.15; \hspace{2pt} -1.30; \hspace{2pt}  0.90] & 53.9  & 0.0017 & 12  \\
    \textbf{CDM 4}  &  4  \si{h} 46 \si{min} 55 \si{s}& [ 8.90  \hspace{2pt}  49.27 \hspace{2pt} -50.88] & [ 0.93; \hspace{2pt} -0.63; \hspace{2pt} 11.56] & 71.4  & 0.0025 & 100 \\ 
    \textbf{CDM 5}  &  5  \si{h} 53 \si{min} 54 \si{s}& [-8.10  \hspace{2pt} -3.02  \hspace{2pt} -20.62] & [ 13.86;\hspace{2pt} -3.71; \hspace{2pt} -4.92] & 22.4  & 0.0138 & 180 \\ 
    \textbf{CDM 6}  &  9 \si{h} 22 \si{min} 43 \si{s}& [-12.64 \hspace{2pt} -14.81 \hspace{2pt}  13.76] & [ 5.64; \hspace{2pt} -3.00  \hspace{2pt}  8.62] & 23.8  & 0.0023 & 90  \\
    \textbf{CDM 7}  &  11 \si{h} 47  \si{min} 47 \si{s}& [ 2.09  \hspace{2pt}  13.74 \hspace{2pt}  1.89]  & [ 7.93; \hspace{2pt} -7.06; \hspace{2pt}  3.21] & 14.0  & 0.0022 & 94  \\ 
    \textbf{CDM 8}  &  13 \si{h} 17 \si{min} 3  \si{s}& [-4.18  \hspace{2pt} -14.19 \hspace{2pt}  14.49] & [ 0.37; \hspace{2pt} -2.23; \hspace{2pt}  1.35] & 20.7  & 0.0020 & 20  \\
    \textbf{CDM 9}  &  14 \si{h} 46 \si{min} 20   \si{s}& [ 8.92  \hspace{2pt} -1.02  \hspace{2pt} -12.78] & [ 1.43; \hspace{2pt}  4.52; \hspace{2pt} 13.58] & 15.6  & 0.0050 & 143 \\
    \textbf{CDM 10} &  15 \si{h} 56 \si{min} 28 \si{s}& [-19.88 \hspace{2pt}  9.10  \hspace{2pt} -4.18]  & [ 0.63; \hspace{2pt} -9.52; \hspace{2pt} -0.31] & 22.3  & 0.0012 & 78  \\ 
    \Xhline{4\arrayrulewidth}
    \end{tabular}
\end{table}

\begin{table}[tb!]
\centering
\caption{Case 1: elements of the relative position covariance matrix of each conjunction.}
\label{tab:covs1}
    \begin{tabular}{l|llllll}
    \Xhline{4\arrayrulewidth}
     & $P_{rr}$ [\si{km^2}] & $P_{tt}$ [\si{km^2}]  & $P_{nn}$ [\si{km^2}] & $P_{rt}$ [\si{km^2}] & $P_{tn}$ [\si{km^2}] & $P_{nr}$ [\si{km^2}] \\
    \Xhline{3\arrayrulewidth}
    \textbf{CDM 1}  & 0.0078  & 0.2887 & 0.0145 &  0.0466 &  0.0103 &  0.0644 \\ 
    \textbf{CDM 2}  & 0.0662  & 0.1279 & 0.1170 & -0.0918 &  0.0879 & -0.1220 \\ 
    \textbf{CDM 3}  & 0.2473  & 0.0022 & 0.0616 &  0.0215 &  0.1232 &  0.0109 \\
    \textbf{CDM 4}  &  0.0004 & 0.1457 & 0.1649 & -0.0055 &  0.0057 & -0.1548 \\ 
    \textbf{CDM 5}  &  0.2579 & 0.0199 & 0.0333 & -0.0710 & -0.0924 &  0.0256 \\
    \textbf{CDM 6}  &  0.0001 & 0.2000 & 0.1110 & -0.0007 &  0.0005 & -0.1486 \\
    \textbf{CDM 7}  &  0.0015 & 0.2027 & 0.1069 & -0.0162 &  0.0120 & -0.1469 \\ 
    \textbf{CDM 8}  &  0.1814 & 0.0012 & 0.1284 & -0.0142 & -0.1523 &  0.0119 \\
    \textbf{CDM 9}  &  0.0028 & 0.0002 & 0.3081 & -0.0003 &  0.0278 & -0.0017 \\
    \textbf{CDM 10} &  0.0004 & 0.1327 & 0.1779 & -0.0036 &  0.0040 & -0.1535 \\
    \Xhline{4\arrayrulewidth}
    \end{tabular}
\end{table}

\subsection{Case 1: Multiple Short-Term Encounters in LEO}
\label{sec:case1}

Most impacts in LEO are head-to-head \cite{ESA2023}, so they involve very high relative velocities.
The test case presented simulates a series of endogenous conjunctions in the Starlink constellation under $J2$ dynamics. A total of ten conjunctions happening in the course of ten orbits of the primary is synthetically generated. The primary is on a circular orbit with an inclination of $53$ \si{deg}, a radius of $6928$ \si{km}, a RAAN, and an argument of periapsis equal to $0$ \si{deg}; at the starting time, $t_0$, its true anomaly is $0$ \si{deg}. The parameters that define each conjunction are reported in \cref{tab:conj1} and \cref{tab:covs1}: the relative states are expressed in the \gls{eci} reference frame, whereas the covariances are in the \gls{rtn} of the secondary. The covariance of the primary is assumed to be already absorbed by that of the secondary, so the number of parameters needed to set the scenario is reduced. $d_{miss} = ||\Delta\vec{r}_{s,CA}||$ is the miss distance of the un-maneuvered case, and $\gamma_v=\mathrm{acos}\left(\Vec{v}_{p,CA}\cdot\vec{v}_{s,CA}\right)$ is the angle between the velocity vectors of the two spacecraft; since the encounters happen between satellites on circular orbits at the same altitude, the magnitude of the velocity of the two spacecraft is equal. The \gls{hbr} is the same for each Starlink satellite and is conservatively assumed to be $3$ \si{m}. The propagation is performed using the harmonics up to the second degree of the gravitational field. The spacecraft is supposed to mount a low-thruster with a maximum thrust of $5.2$ \si{mN}, which, considering the Starlink's mass of 230 kg, corresponds to a maximum acceleration of $0.02$ \si{mm/s^2}.
\begin{figure}[tb!]
    \centering
    \subfloat[2 conjunctions.]{\includegraphics[width = 0.48\textwidth]{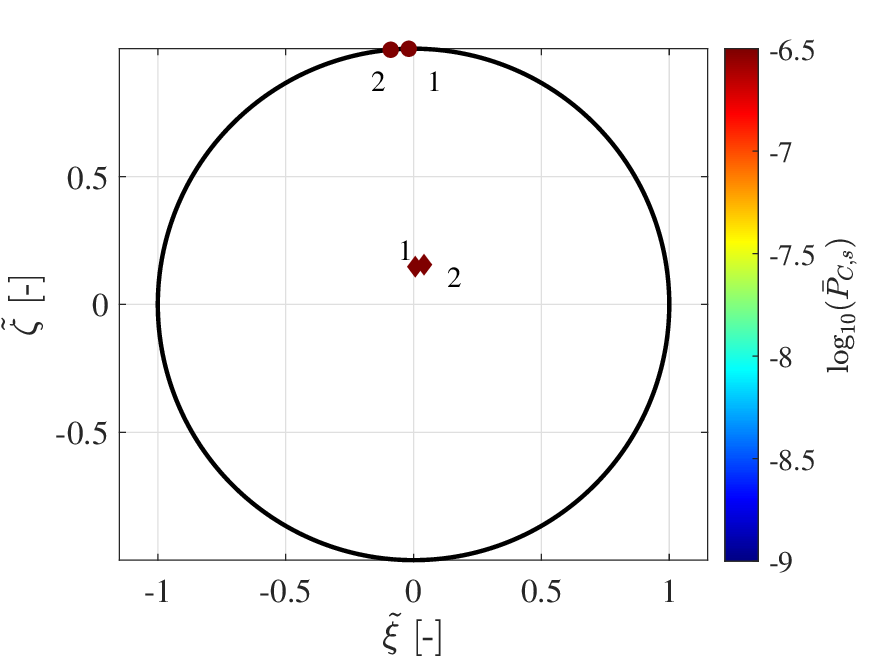}\label{fig:Bplane_1}}
    \subfloat[5 conjunctions.]{\includegraphics[width = 0.48\textwidth]{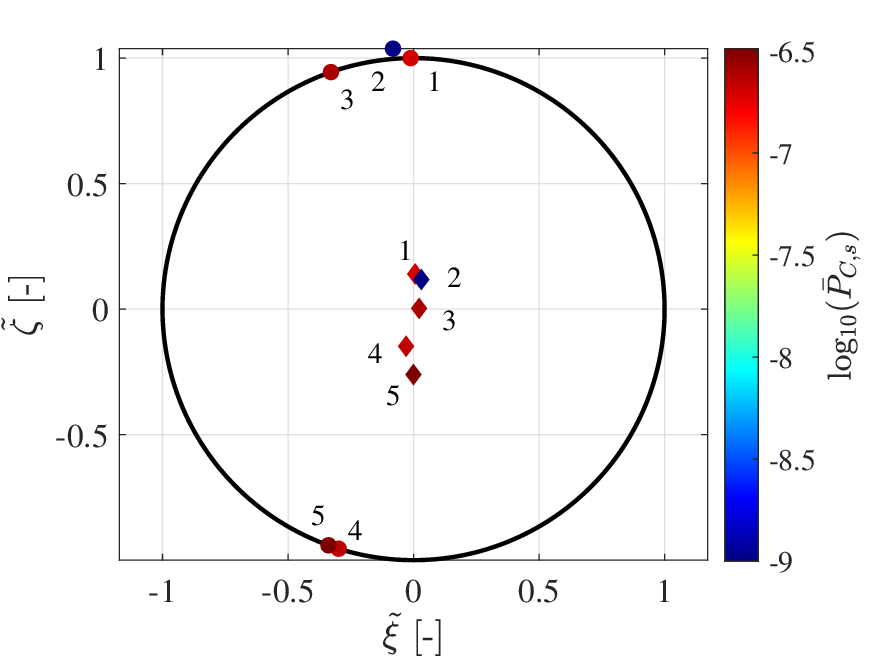}\label{fig:Bplane_2}} \\
    \subfloat[7 conjunctions.]{\includegraphics[width = 0.48\textwidth]{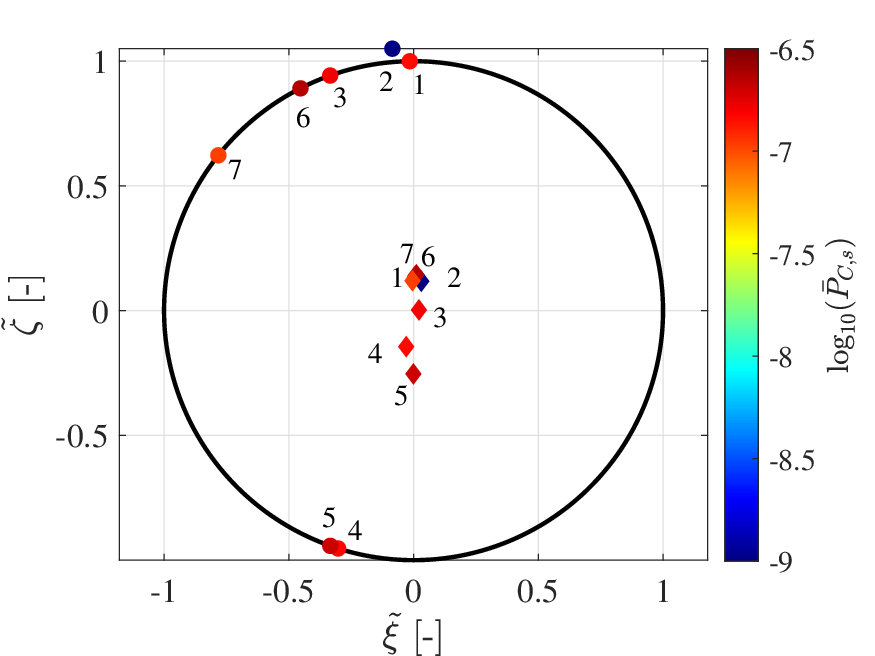}\label{fig:Bplane_3}}
    \subfloat[10 conjunctions.]{\includegraphics[width = 0.48\textwidth]{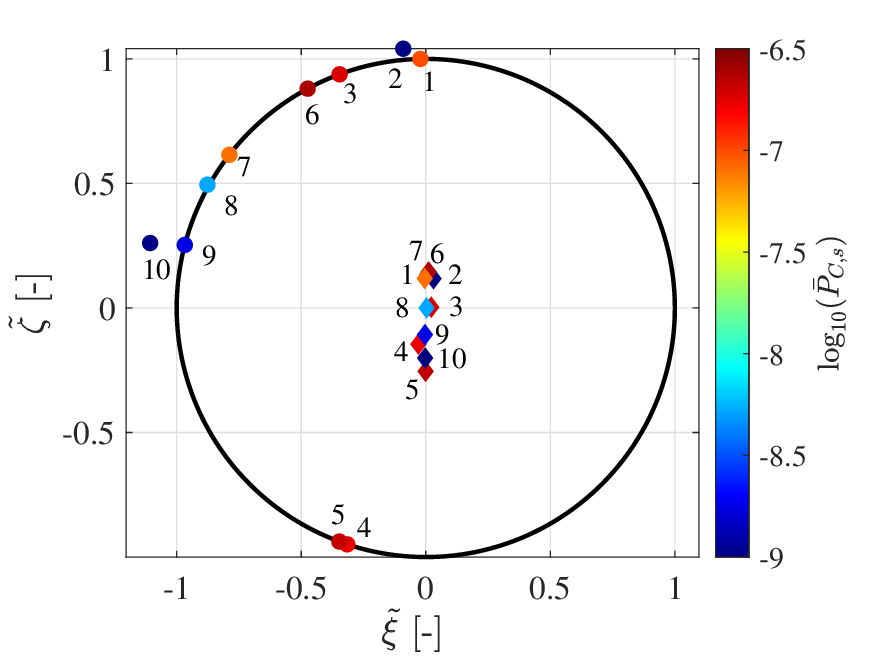}\label{fig:Bplane_4}}
    \caption{Case 1: equivalent B-plane representation of the conjunctions.} 
    \label{fig:bplane_Case1}
\end{figure}
The algorithm is run with both the refinements methods introduced in \cref{sec:caShort}, and the results are compared.

\begin{figure}[tb!]
    \centering
    \subfloat[2 conjunctions.]{\includegraphics[width = 0.48\textwidth]{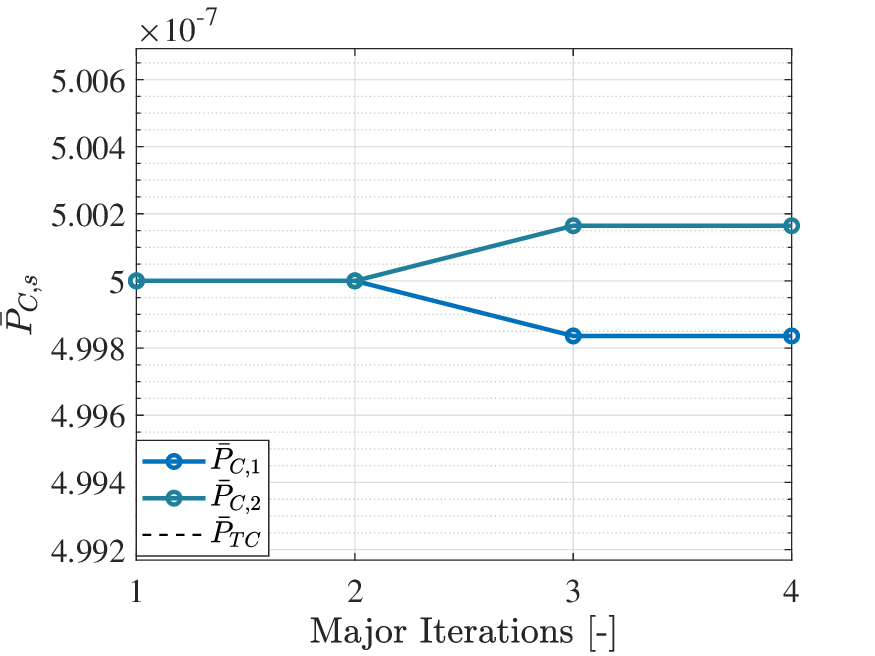}\label{fig:limAdapt_1}}
    \subfloat[5 conjunctions.]{\includegraphics[width = 0.48\textwidth]{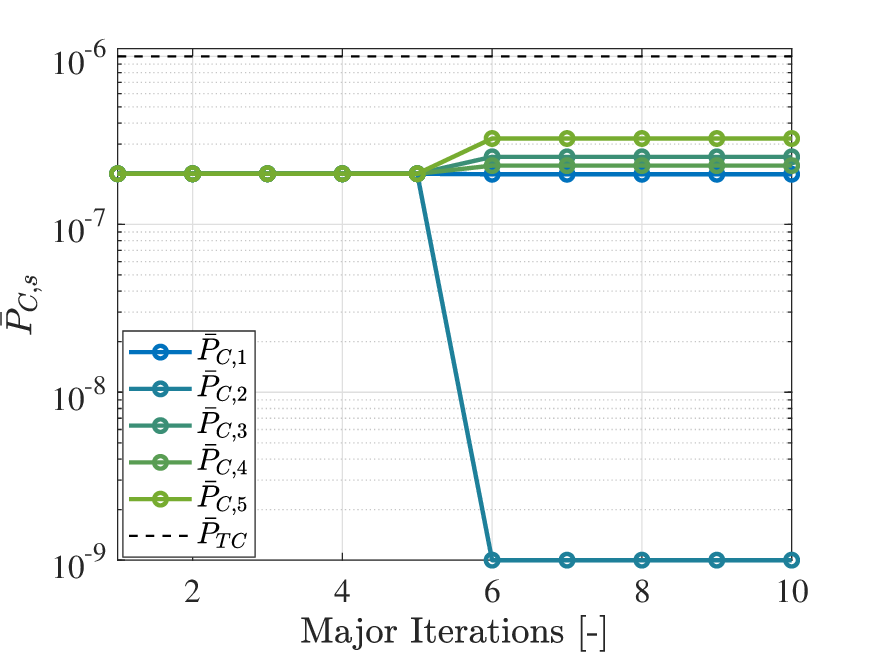}\label{fig:limAdapt_2}} \\
    \subfloat[7 conjunctions.]{\includegraphics[width = 0.48\textwidth]{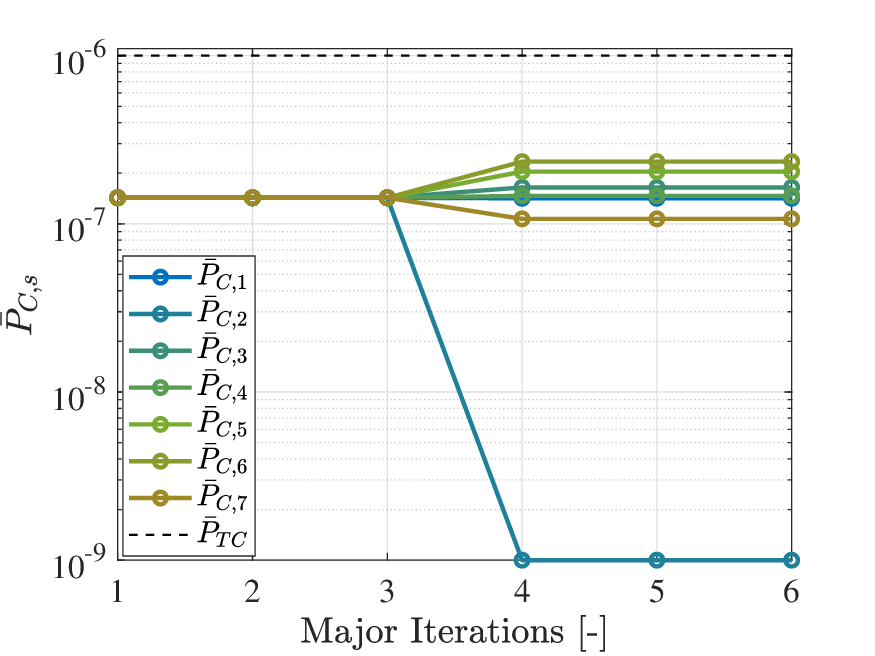}\label{fig:limAdapt_3}}
    \subfloat[10 conjunctions.]{\includegraphics[width = 0.48\textwidth]{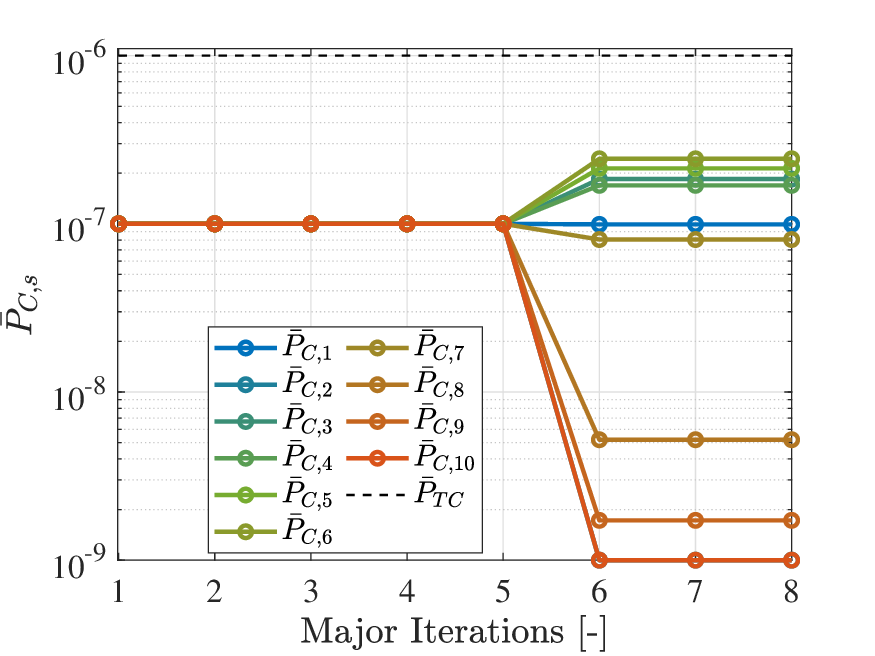}\label{fig:limAdapt_4}}
    \caption{Case 1: evolution of the $\bar{P}_{C,s}$ of the single conjunctions.} 
    \label{fig:limitAdapt_Case1}
\end{figure}

\begin{figure}[tb!]
    \centering
    \subfloat[2 conjunctions.]{\includegraphics[width = 0.48\textwidth]{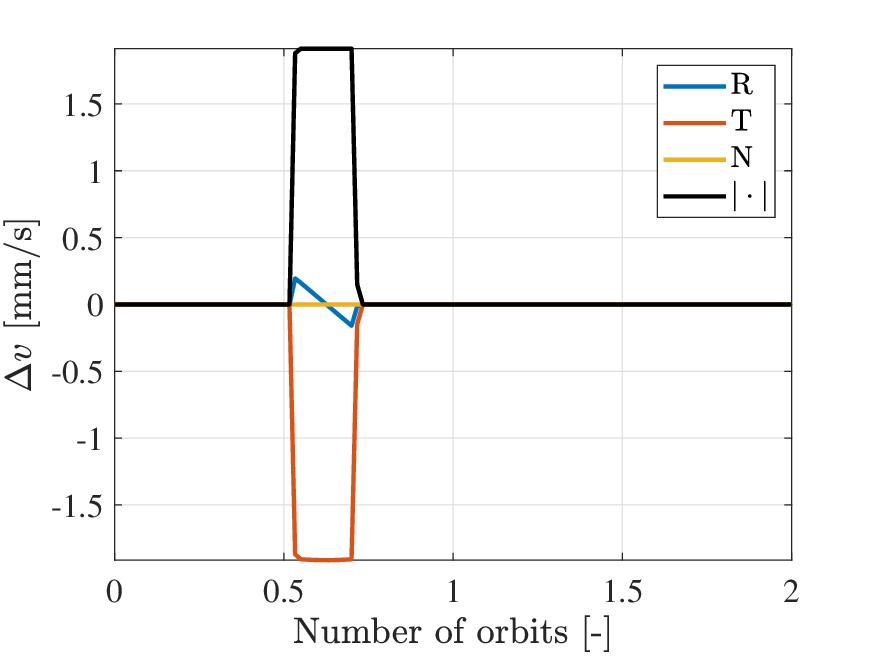}\label{fig:dv_1}}
    \subfloat[5 conjunctions.]{\includegraphics[width = 0.48\textwidth]{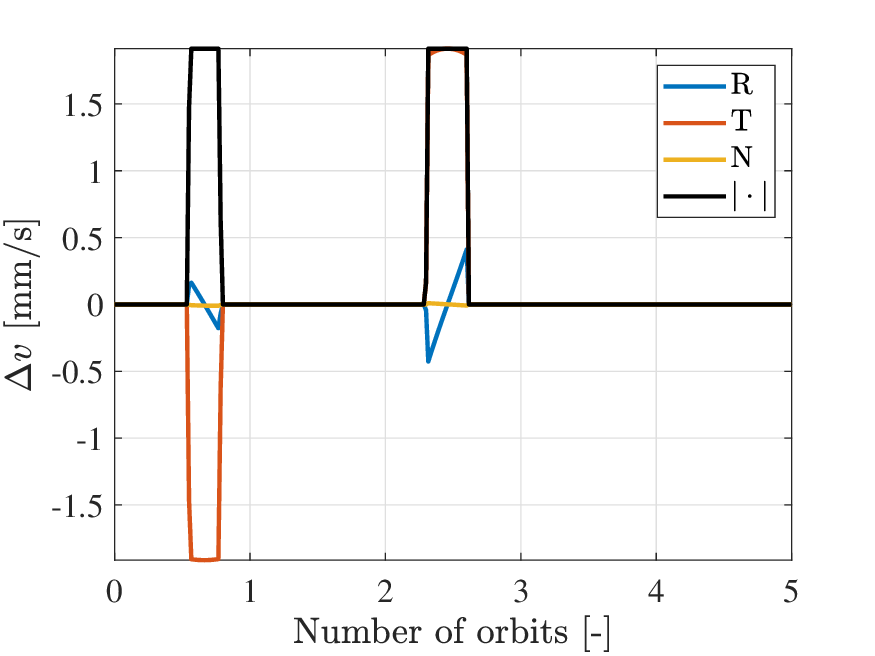}\label{fig:dv_2}} \\
    \subfloat[7 conjunctions.]{\includegraphics[width = 0.48\textwidth]{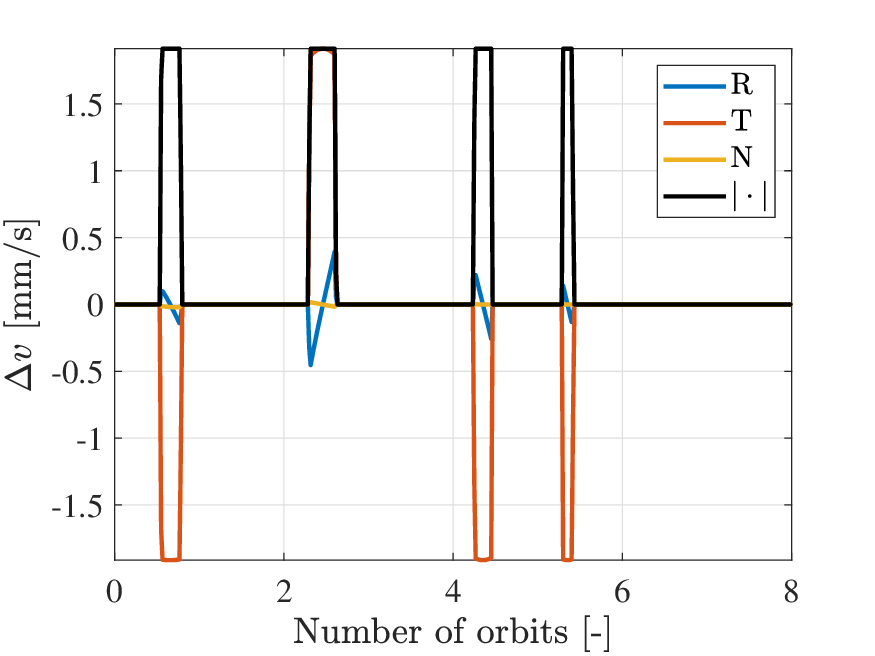}\label{fig:dv_3}}
    \subfloat[10 conjunctions.]{\includegraphics[width = 0.48\textwidth]{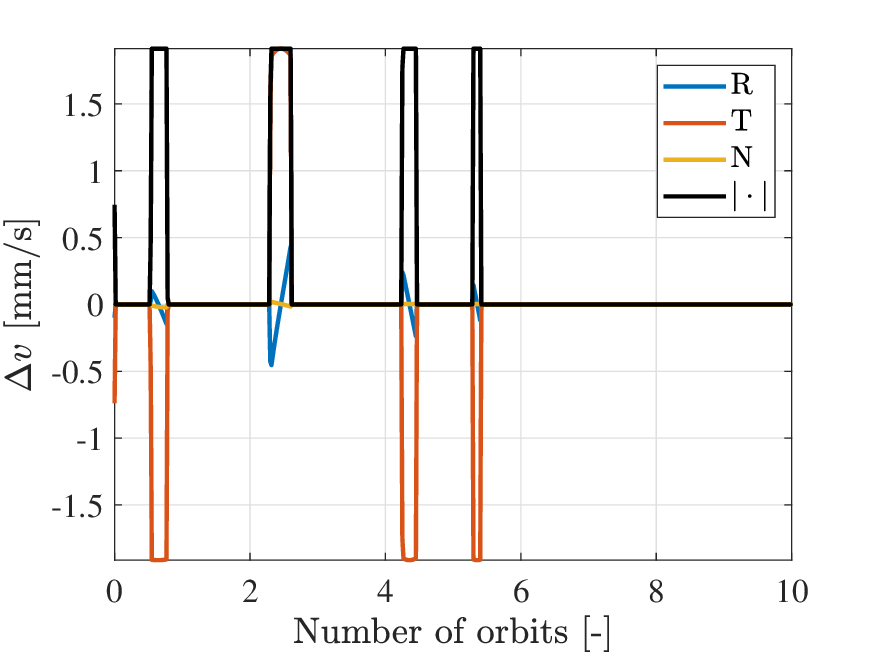}\label{fig:dv_4}}
    \caption{Case 1: $\dv$ to perform the \gls{cam}.} 
    \label{fig:dv_Case1}
\end{figure}

In \cref{fig:bplane_Case1}, the equivalent B-planes of the conjunctions for each case considered are shown. To represent all the conjunctions in the same plot, a coordinate transformation is performed on the \gls{koz} ellipse and on the relative position points. The transformation comprises a rotation to diagonalize the covariance matrix and a stretch to circularize it. Considering a generic vector $\vec{w} = [\xi \hspace{4pt} \zeta]\transp$ expressed in the B-plane reference frame of a conjunction $s$, the transformed vector $\Tilde{\vec{w}}= [\tilde{\xi} \hspace{4pt} \tilde{\zeta}]\transp$ is
\begin{equation}
  \tilde{\vec{w}} = \mathrm{diag}\left([1/a; \hspace{4pt} 1/b]\right) \Vec{V}\vec{w},
\end{equation}
where $a,b\in\mathbb{R}_\mathrm{+}$ are the semiaxes of the ellipse associated with the optimized probability threshold, and $\Vec{V}\in\mathbb{R}^{2\times2}$ is the rotation matrix that diagonalizes the covariance. The ballistic trajectory points are shown as diamonds, and the optimized points obtained using the linearized \gls{smd} constraint with limit adaptation are shown as circles. The conjunction number is indicated for each marker, and the color of the marker is associated with the value of the limit \gls{poc} of the corresponding conjunction after the limit adaptation. The result of the same limit adaptation is reported in \cref{fig:limitAdapt_Case1}: after the first \gls{scp} is solved, in all the considered cases, a single \gls{scp} is needed to find a solution that grants a \gls{tpoc} value sufficiently close to the threshold. We notice from \cref{fig:limitAdapt_Case1} that depending on the conjunctions configuration and their number, the algorithm automatically finds the most important ones and assigns them a higher weight. For example, the second conjunction is equally important to the first one in the first scenario, but it becomes irrelevant in the next ones. As a second example, the last three conjunctions that differentiate the last scenario from the third have no effect on the final result since their $\bar{P}_C$ value drops significantly after the adaptation. This is why the computed maneuvers for the last two scenarios are almost identical (see \cref{fig:dv_3} and \cref{fig:dv_4}). 

Through an analysis of this type, it is also possible to identify clusters of conjunctions that cause a maneuver. Indeed, the time, direction, and magnitude of the first two burns remain almost unchanged across the four cases. In general, the maneuvers are all mainly in the tangential direction. This confirms the idea proposed by other researchers \cite{Bombardelli2015,Hernando-Ayuso2020} that raising or lowering the orbit is the most effective way of performing the \gls{cam} if enough warning time is given.

In \cref{tab:res1}, the results in terms of final \gls{tpoc} and $\dv$ are shown for the four combinations of conjunctions considered. Convergence results are also compared in the number of iterations, convergence time, and maximum positional validation error. The two refinement methods do not differ greatly in the resulting $\dv$ profile or in the final equivalent B-plane configuration, as is shown in \cref{tab:res1}. Nonetheless, the method that relies on the linearized \gls{smd} constraint with limit adaptation yields more accurate results and can target the desired \gls{tpoc} threshold with higher precision. The \gls{poc} value at \gls{tca} for the first method is always kept below the threshold value with a margin of $.002\%$, while the linearized \gls{tpoc} method reaches \gls{tpoc} values up to $42\%$ below the threshold. This means that the \gls{cam} optimization does not find the optimal result with the second method, leading to higher $\dv$ impulses than the first method. The run times of the two methods are similar when a low number of conjunctions is considered, but the linearized \gls{smd} method is definitely faster when we start to consider a high $n_{conj}$. In general, we can infer that the linearized \gls{smd} constraint is more appropriate when dealing with more than two conjunctions because each \gls{poc} constraint is linearized separately, so the coarseness of approximation of the nonlinearity does not scale with the complexity of the problem. 

\begin{table}[tb!]
\centering
\caption{Case 1: compared results for the \gls{smd} and \gls{tpoc} linearized constraints.}
\label{tab:res1}
\begin{tabular}{l|llll|llll}
\Xhline{4\arrayrulewidth}
\textbf{Constraint}           & \multicolumn{4}{c|}{\textbf{Linearized \gls{smd}}}  & \multicolumn{4}{c}{\textbf{Linearized \gls{tpoc}}}     \\ \hline
$n_{conj}$                 & 2            & 5            & 7             & 10           & 2            & 5           & 7             & 10             \\ \Xhline{3\arrayrulewidth}
$\dv$ [\si{mm/s}]          & $21.16$      & $60.05$      & $102.04$      & $102.14$     & $21.16$      & $62.32$     & $103.18$      & $105.58$       \\ \hline
$P_{TC}\times 10^{6}$ [-]  & $1.001$      & $1.001$      & $1.001$       & $0.998$      & $1.002$      & $0.801$     & $0.824$       & $0.583$        \\ \hline
$n_{major}$ [-]            & 4            & 6            & 6             & 8            & 5            & 7           & 10            & 20             \\ \hline
$n_{minor}$ [-]            & 5            & 9            & 10            & 13           & -            & -          & -            & -             \\ \hline
$e_{validation}$ [\si{mm}] & $0.92$       & $1.88$       & $0.72$        & $5.12$       & $37.61$      & $5.48$      & $43.40$       & $1940.80$      \\ \hline
$T_{sim}$ [\si{s}]         & $4.53$       & $8.34$       & $14.52$       & $28.03$      & $4.47$       & $8.82$     & $20.66$        & $62.15$       \\ \Xhline{4\arrayrulewidth}
\end{tabular}
\end{table}

\subsection{Case 2: Repeating Short-Term Encounter in LEO with GMM}
Let us consider a scenario in which a secondary spacecraft's orbit intersects the primary's twice, thus generating a repeating encounter. The first encounter is described by a \gls{cdm}, whereas the second one is detected by propagating the dynamics. The covariance of the secondary at the first \gls{tca} is split into a number of \gls{gmm} mixands, which are propagated separately until the second \gls{tca}. Here, we assume that the covariance of the primary is negligible w.r.t. to the secondary's. The state of the primary is described at \gls{tca} by $a = 6800$ \si{km}, $e=0$, $\omega = \Omega = \theta = 0$ \si{deg}, and $i = 91.67$ \si{deg}; the orbit of the secondary is elliptical, has opposite inclination w.r.t. the primary and its period is a multiple of the period of the primary. The full secondary's covariance matrix at the first \gls{tca} is defined in \gls{rtn} coordinates in \cref{tab:covcase2}. The primary covariance matrix is assumed to be null.

\begin{table}[tb!]
\centering
\caption{Case 2: elements of the relative position covariance matrix at the first \gls{tca}.}
\label{tab:covcase2}
    \begin{tabular}{l|llllll}
    \Xhline{4\arrayrulewidth}
    \textbf{Spacecraft} & $P_{rr}$ [\si{m^2}] & $P_{tt}$ [\si{m^2}]  & $P_{nn}$ [\si{m^2}] & $P_{\dot{r}\dot{r}}$ [\si{mm^2/s^2}] & $P_{\dot{t}\dot{t}}$ [\si{mm^2/s^2}] & $P_{\dot{n}\dot{n}}$ [\si{mm^2/s^2}] \\
    \Xhline{3\arrayrulewidth}
    Primary  & 0  & 0 & 0 &  0 &  0 &  0 \\
    \hline
    Secondary  & 2.025  & 10.000 & 2.500 & 0.225 &  5.625 & 0.625 \\ 
    \Xhline{4\arrayrulewidth}
    \end{tabular}
\end{table}

According to reference \cite{Yanez2019}, when propagating Cartesian coordinates in \gls{leo}, for typical initial uncertainty, after five orbits ca., the normality of the distribution is lost. So, we construct the test such that after six orbits of the primary, the head-to-head encounter happens again. The eccentricity of the secondary orbit is $0.6971$; its semi-major axis is $22453.1$ \si{km}, the inclination is $181.67$ \si{deg}, the argument of perigee is $0$ \si{deg}, the RAAN is $1.72\cdot 10^{-3}$ \si{deg} and the true anomaly at the first conjunction is $0.11$ \si{deg}. The orbits are propagated backward for one orbital period of the primary and forward up to the second conjunction. In this way, the first conjunction is treated as a standard one, i.e., without splitting, and the second one is treated with the propagated \gls{gmm} components. In this case, we are considering a chemical propulsion system with a maximum instantaneous acceleration of $1$ \si{mm/s^2}. Different simulations with different $n_{mix}$ are compared to one where no split is performed, i.e., the second conjunction is treated with a linearly propagated covariance. Moreover, the solutions obtained with the two constraints are compared in \cref{tab:res2}.

\begin{figure}[tb!]
    \centering
    \subfloat[$n_{mix}=1$.]{\includegraphics[width = 0.48\textwidth]{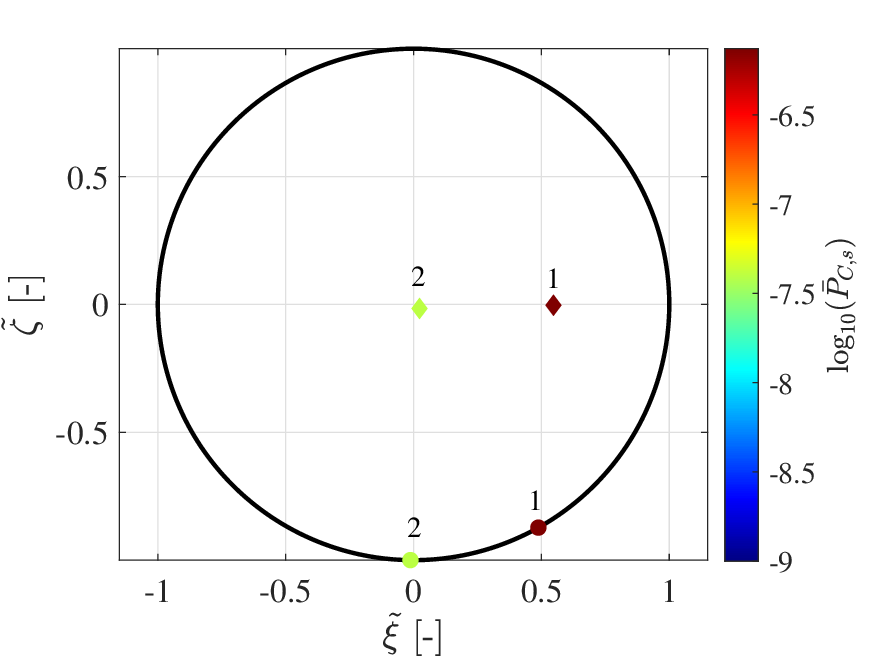}\label{fig:bplane_gmm_short_1}}
    \subfloat[$n_{mix}=3$.]{\includegraphics[width = 0.48\textwidth]{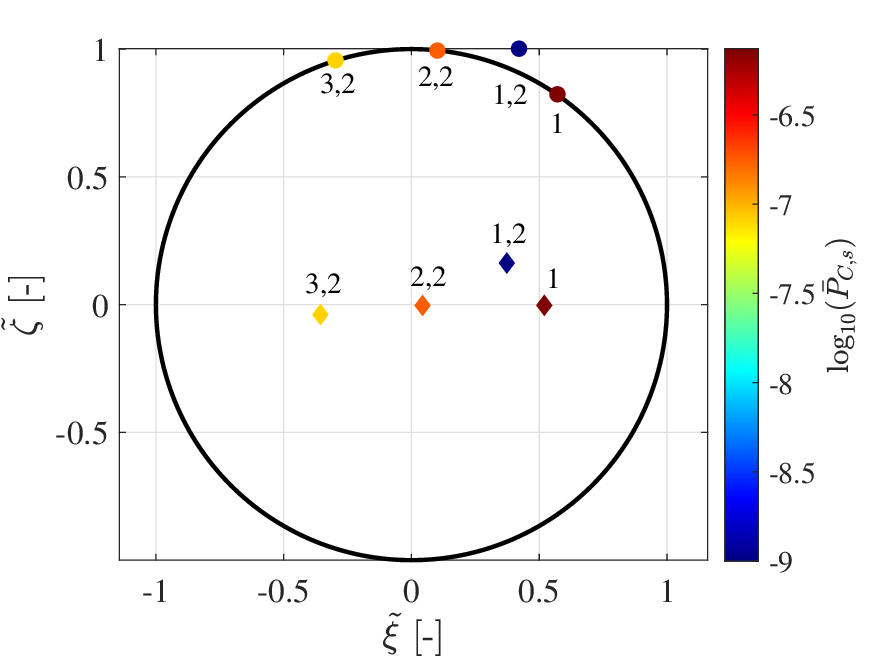}\label{fig:bplane_gmm_short_2}} \\
    \subfloat[$n_{mix}=5$.]{\includegraphics[width = 0.48\textwidth]{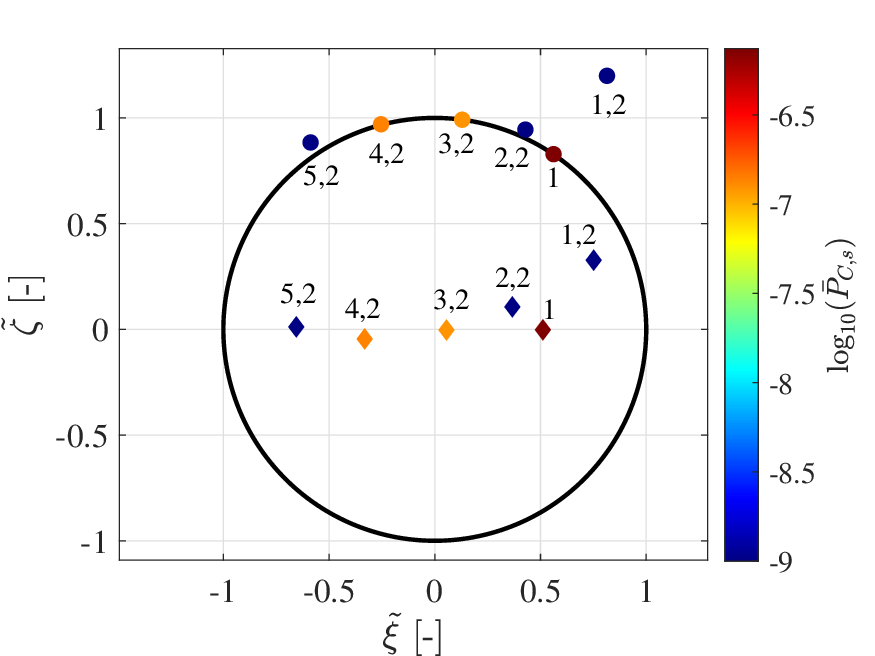}\label{fig:bplane_gmm_short_3}}
    \subfloat[$n_{mix}=7$.]{\includegraphics[width = 0.48\textwidth]{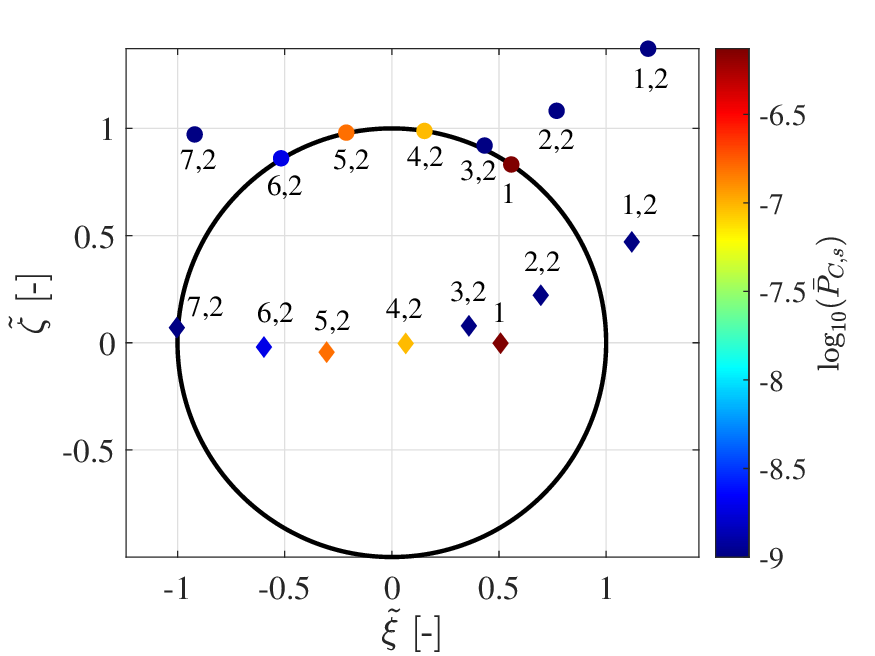}\label{fig:bplane_gmm_short_4}}
    \caption{Case 2: equivalent B-plane representation of the conjunctions.} 
    \label{fig:bplane_Case2}
\end{figure}

In \cref{tab:res2}, we notice that the use of the \gls{gmm} permits better to represent the \gls{koz} and reduce the required $\dv$ by at least $4$ \si{mm/s}. \cref{fig:bplane_Case2} shows the equivalent B-planes of the two conjunctions for different values of $n_{mix}$. The first conjunction is always dominant: since it appears first avoiding it means that the satellite will most likely avoid the subsequent one as well. Nonetheless, in the other plots, it is useful to see which mixands play a role in defining the limits: in particular, in \cref{fig:bplane_gmm_short_2} mixands 3 and 2, in \cref{fig:bplane_gmm_short_3} mixands 4 and 3 and in \cref{fig:bplane_gmm_short_4} mixands 5 and 4 are the ones that are not led to the lowest threshold. This means that in the second (extended) conjunction, the primary passes closest to the central mixand and the one right after it compared to all the others. For all the other mixands, the threshold \gls{poc} is brought down to the minimum value of $10^{-9}$, meaning they do not play a role in defining the maneuver.  
\begin{figure}[tb!]
    \centering
    \subfloat[$n_{mix}=1$.]{\includegraphics[width = 0.48\textwidth]{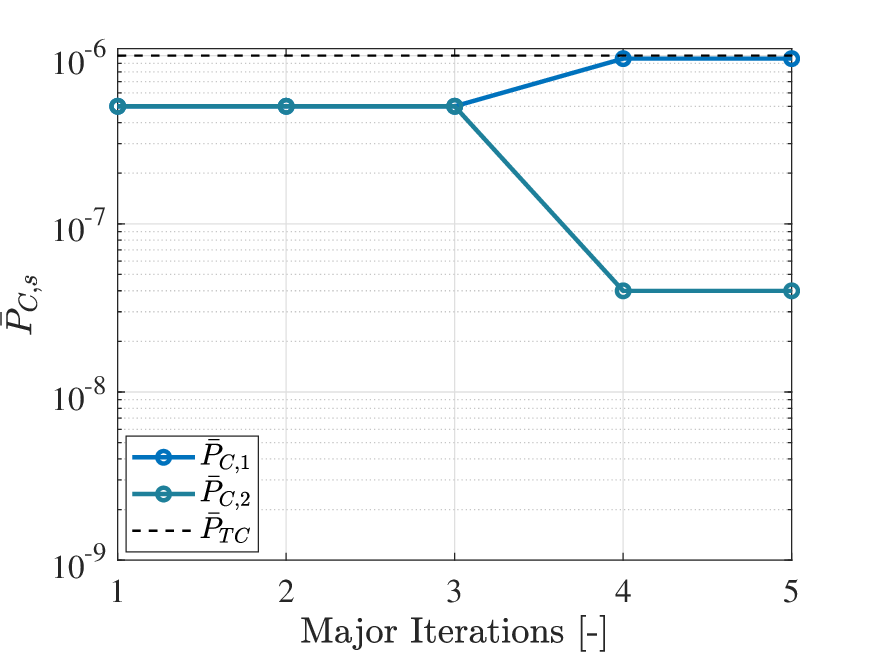}\label{fig:limAdapt_gmm_1}}
    \subfloat[$n_{mix}=7$.]{\includegraphics[width = 0.48\textwidth]{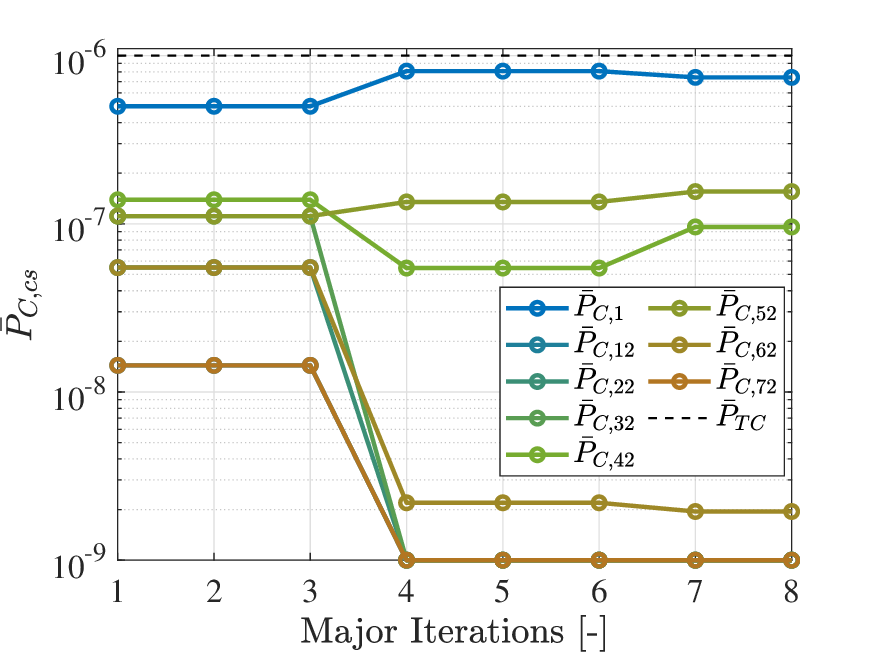}\label{fig:limAdapt_gmm_2}} 
    \caption{Case 2: evolution of the $\bar{P}_{C,s}$ of the single conjunctions and mixands.} 
    \label{fig:limAdapt_Case2}
\end{figure}

Depending on the number of mixands used, the \glspl{tca} of the single mixands are between $1.32$ and $2.91$ \si{s} apart from each other, which is accounted for in the selection of the time grid for the discretization. This happens because the encounter is head-to-head, and the split is performed very close to the direction of the velocity of the secondary so that the primary runs into each mixand consecutively, as \cref{fig:head2head} illustrates. In \cref{fig:limAdapt_Case2}, the adaptation of the limits for two representative numbers of mixand is shown: with no \gls{gmm} split, the first conjunction dominates the final \gls{tpoc} computation, while in the case with $n_{mix}=7$, the second conjunction is considered slightly more important. 
In \cref{tab:dv_Case2}, one sees that, in accordance with most published research on single encounters \cite{Hernando-Ayuso2020}, in all the considered cases, the $\dv$ is almost purely tangential. Nonetheless, in the case where no split is performed, the direction of thrust is opposite to that of the other cases. In \cref{fig:bplane_gmm_short_1}, the reader can see that the tangential maneuver moves the relative position along the $\tilde{\xi}$ axis, so there is little to no difference in performing a maneuver on the positive or the negative tangential direction and the optimizer prefers one over the other. In all cases with a split, the preferred maneuver is in the positive tangential direction. The number of mixands considered does not affect the maneuvering time. The maneuver is always performed half an orbit before the first \gls{tca}, and an almost negligible second maneuver is performed half an orbit before the second \gls{tca}.

Regarding the convergence properties, the two methods behave similarly to the case in \cref{sec:case1}, with the first method being slightly faster when more mixands are considered. In this case, though, it is important to point out that because of the coarseness of the linearization, the linearized \gls{tpoc} method always reaches a solution that is around $5\%$ over the \gls{tpoc} threshold. This allows it to find a lower $\dv$ at the cost of precision. The validation error for both methods is always acceptable. 


\begin{table}[tb!]
\centering
\caption{Case 2: compared results for the \gls{smd} and \gls{tpoc} linearized constraints.}
\label{tab:res2}
\begin{tabular}{l|llll|llll}
\Xhline{4\arrayrulewidth}
\textbf{Constraint}           & \multicolumn{4}{c|}{\textbf{Linearized \gls{smd}}}  & \multicolumn{4}{c}{\textbf{Linearized \gls{tpoc}}}     \\ \hline
$n_{mix}$                  & 1            & 3         & 5        & 7         & 1        & 3       & 5         & 7         \\ \Xhline{3\arrayrulewidth}       
$\dv$ [\si{mm/s}]          & $57.06$      & $52.64$   & $52.98$  & $53.37$   & $56.52$  & $52.23$ & $52.55$   & $52.15$   \\ \hline
$P_{TC}\times 10^{6}$ [-]  & $1.002$      & $0.968$   & $0.936$  & $0.910$   & $1.053$  & $1.039$ & $1.038$   & $1.048$   \\ \hline
$n_{major}$ [-]            & 5            & 8         & 8        & 8         & 7        & 7       & 8         & 9          \\ \hline
$n_{minor}$ [-]            & 12           & 16        & 16       & 19        & -       & -      & -        & -            \\ \hline
$e_{validation}$ [\si{mm}] & $9.92$       & $3.64$    & $3.21$   & $2.98$    & $0.46$   & $3.81$  & $4.81$    & $10.61$   \\ \hline
$T_{sim}$ [\si{s}]         & $6.57$       & $19.76$   & $23.01$  & $27.68$   & $9.15$   & $16.98$ & $24.58$   & $30.17$   \\ \Xhline{3\arrayrulewidth}
\end{tabular}
\end{table}

\begin{table}
\centering
\caption{Case 2: $\dv$ to perform the \gls{cam}.}
\label{tab:dv_Case2}
\begin{tabular}{l|llll|llll}
\Xhline{4\arrayrulewidth}
$n_{mix}$                & $t_1$ [s] & $\mathrm{R}_1$ [\si{mm/s}]            & $\mathrm{T}_1$ [\si{mm/s}]            & $\mathrm{N}_1$ [\si{mm/s}]  & $t_2$ [s] & $\mathrm{R}_2$ [\si{mm/s}]            & $\mathrm{T}_2$ [\si{mm/s}]            & $\mathrm{N}_2$ [\si{mm/s}]  \\ \Xhline{3\arrayrulewidth}
$1$ & $2790.26$  &  $-0.69$  &  $-56.63$    &  $-0.04$  & $25112.31$  & $0$        & $-0.04$  & $0$ \\ \hline
$3$ & $2790.26$  &  $-0.68$   &  $52.60$    &  $0.01$   & $25112.31$  & $0$        & $0.03$   & $0$ \\ \hline
$5$ & $2790.26$  &  $-0.67$   &  $53.04$    &  $0.01$   & $25112.31$  &  $0$   $0.08$   &   & $0$ \\ \hline
$7$ & $2790.26$  &  $-0.01$   &  $53.37$    &  $0.25$   & -           & -          & -        & -   \\
\Xhline{4\arrayrulewidth}\end{tabular}
\end{table}

\subsection{Case 3: Single Long-Term Encounter in LEO with GMM}
In this scenario, a single long-term conjunction in \gls{leo} is studied with \gls{gmm}. The problem becomes equivalent to having a long-term conjunction with multiple secondaries. 
The scenario is taken from the \gls{leo} test case of reference \cite{Pavanello2023}: the primary spacecraft is on a circular \gls{leo} with $a=6800$ \si{km}, and $\omega=\Omega=i=\theta=0$ \si{deg}; the relative state of the secondary as expressed in the primary's \gls{rtn} at \gls{tca}, i.e., at the starting time is 
\begin{equation*}
    \vec{x}_{rel} = [-1.004814 \hspace{5pt} 0 \hspace{5pt} 0 \hspace{5pt} 0.112790 \hspace{5pt} 3.388089 \hspace{5pt} 0.011228]\transp,
\end{equation*}
where the first three elements are in [\si{km}] and the last three in [\si{m/s}]. The covariance matrices of the two spacecraft are reported in \cref{tab:covcase3}. The combined \gls{hbr} is $32$ \si{m}, and the maximum instantaneous acceleration achievable by the primary is 0.2 \si{mm/s^2}, corresponding to a maximum thrust of $40$ \si{mN} in a $200$ \si{kg} satellite. The dynamics include J2 perturbation, lunisolar attraction, atmospheric drag, and \gls{srp}. The models used to compute the perturbation forces are computed with \gls{aida} \cite{Morselli2014High}. The physical parameters of the spacecraft are reported in \cref{tab:paramsLeo}.

\begin{table}[tb!]
\centering
\caption{Case 3: elements of the relative position covariance matrix of the two spacecraft.}
\label{tab:covcase3}
    \begin{tabular}{lllllll}
    \Xhline{4\arrayrulewidth}
   \textbf{Spacecraft} & $C_{rr}$ [\si{m^2}] & $C_{tt}$ [\si{m^2}] & $C_{nn}$ [\si{m^2}] & $C_{\dot{r}\dot{r}}$ [\si{m^2/s^2}] & $C_{\dot{t}\dot{t}}$ [\si{m^2/s^2}] & $C_{\dot{n}\dot{n}}$ [\si{m^2/s^2}] \\
    \Xhline{3\arrayrulewidth}
    Primary    & $0.625$ & $10$ & $3.025$ & $0.00625$ & $0.05625$ & $0.00225$ \\ \hline 
    Secondary  & $5.625$ &$ 90$ & $27.225$ & $0.05625$ & $0.50625$ & $0.02025$ \\ \hline 
    \Xhline{4\arrayrulewidth}
    \end{tabular}
\end{table}

\begin{table}[b!]
\centering
\caption{\gls{leo} scenario: physical properties of the spacecraft.}
\label{tab:paramsLeo}
    \begin{tabular}{lllllll}
    \Xhline{4\arrayrulewidth}
    \textbf{Spacecraft} & $m$ [\si{kg}] & $A_{drag}$ [\si{m^2}]  & $C_D$ [-] & $A_{SRP}$ [\si{m^2}] & $C_r$ [-] & \gls{hbr} [\si{m}]  \\
    \Xhline{3\arrayrulewidth}
    Primary   & $200$ & $1$     & $2.2$ & $1$    & $1.31$ & $25$\\ \hline 
    Secondary & $50$  & $0.05$  & $2$   & $0.05$ & $1.31$ & $7$\\
    \Xhline{4\arrayrulewidth}
    \end{tabular}
\end{table}

\begin{figure}[tb!]
    \centering
    \subfloat[Before the maneuver.]{\includegraphics[width = \textwidth]{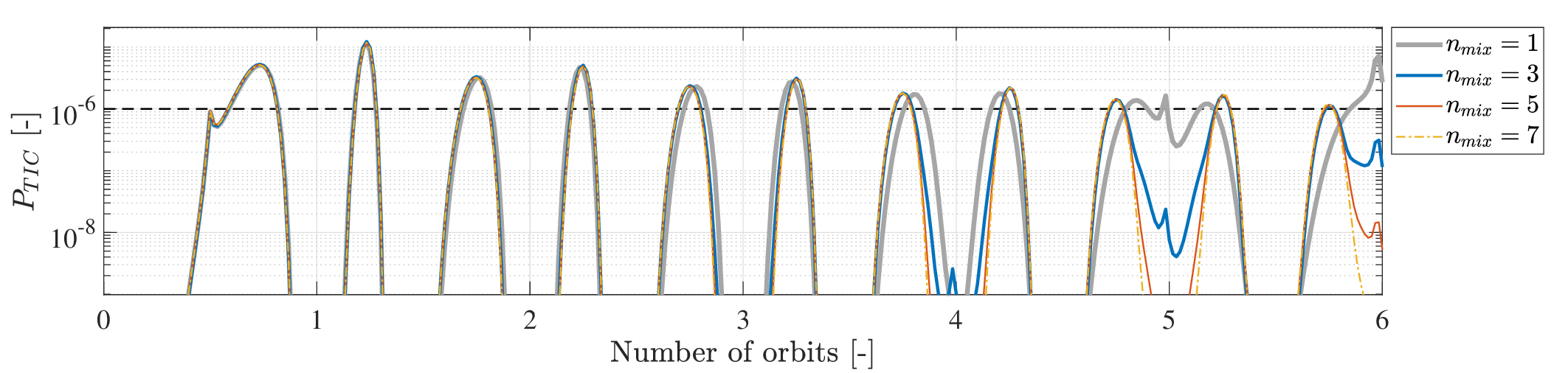}\label{fig:origIpc}} \\
    \subfloat[After the maneuver, \gls{smd} constraint.]{\includegraphics[width = \textwidth]{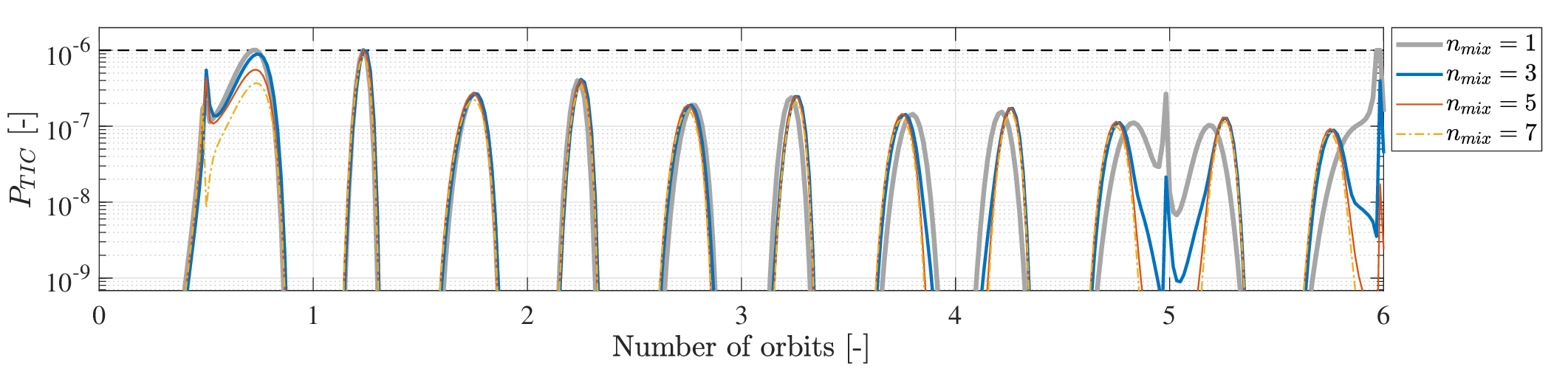}\label{fig:optIpc}} \\
    \subfloat[After the maneuver, \gls{tipc} constraint.]{\includegraphics[width = \textwidth]{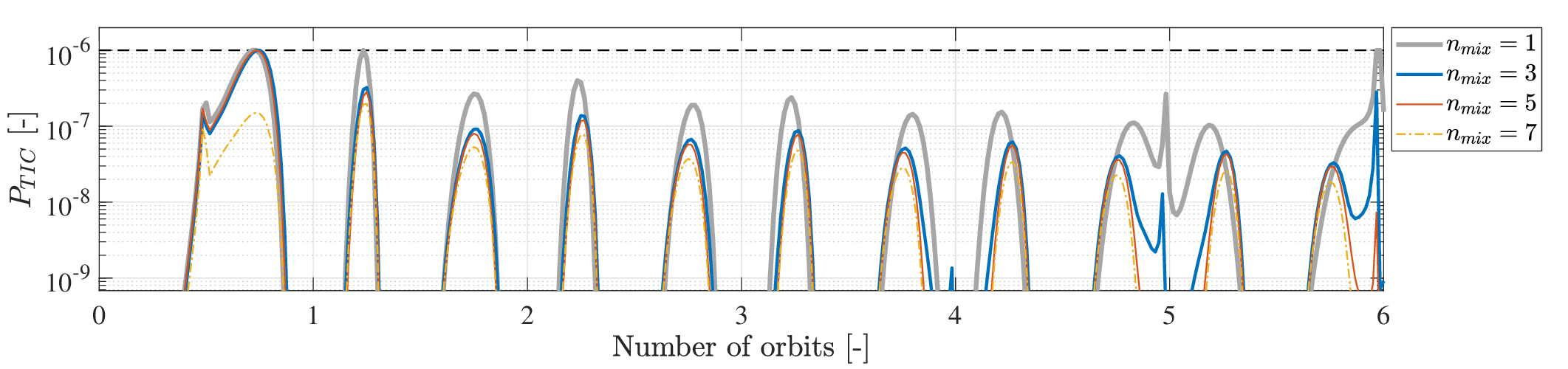}\label{fig:optIpc_ipcConstr}} 
    \caption{Case 3: \gls{tipc} profile for different $n_{mix}$ in the long-term encounter scenario.}
    \label{fig:tipc}
\end{figure}

\begin{figure}[tb!]
    \centering
    \subfloat[$n_{mix}=3$.]{\includegraphics[width=0.48\textwidth]{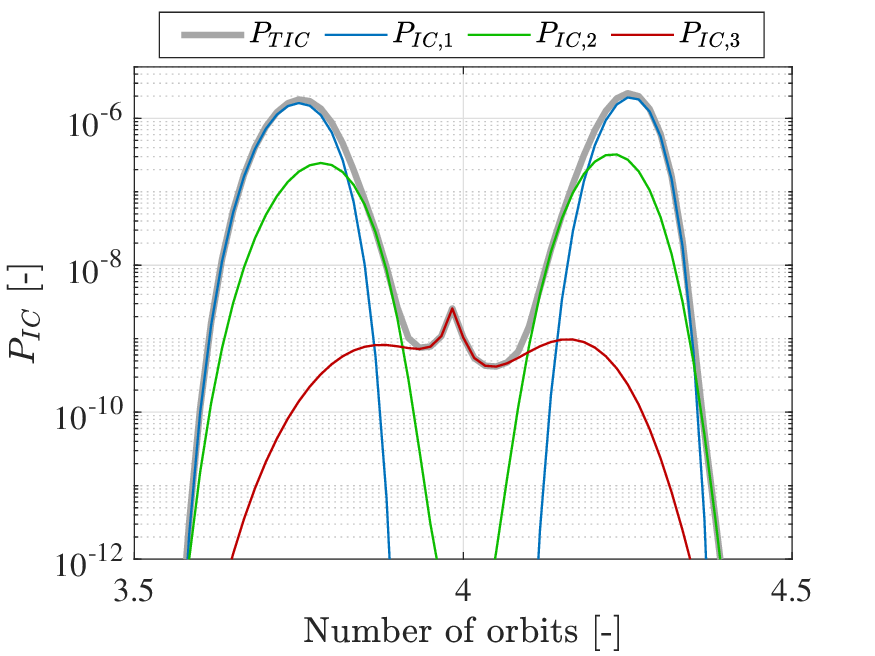}\label{fig:ipc_comb3}}\subfloat[$n_{mix}=7$: the contribution of mixands 6 and 7 is negligible.]{\includegraphics[width=0.48\textwidth]{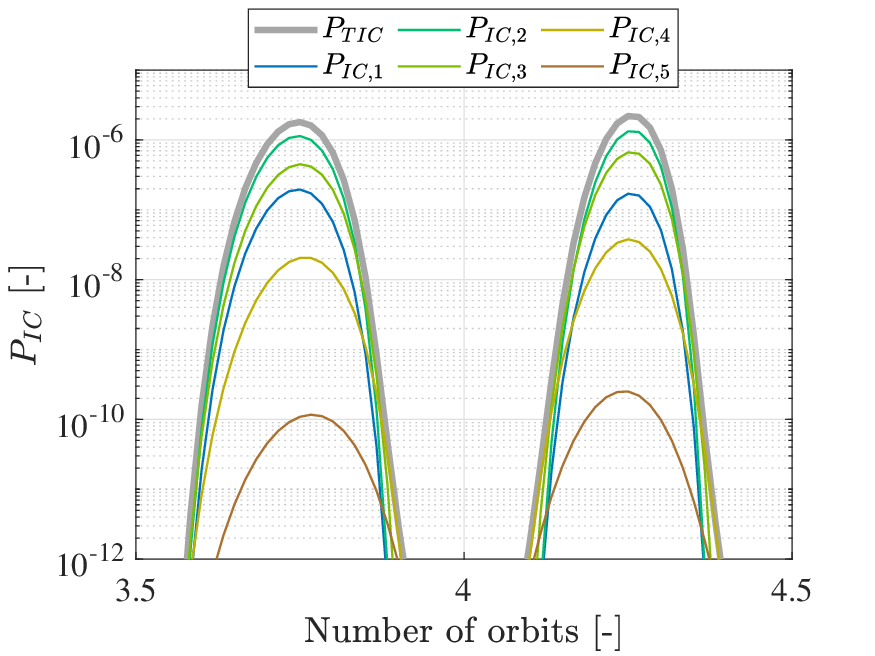}\label{fig:ipc_comb7}} 
    \caption{Case 3: Breakdown of the contribution of each mixand to the \gls{tipc}. Zoom in a central orbit.}
    \label{fig:ipc_comb}
\end{figure}

The propagation spans six orbits of the primary from \gls{tca}. Different orders of the \gls{gmm} are compared in terms of $\dv$, \gls{tipc} profile, and convergence properties in \cref{tab:res3} and \cref{fig:tipc}. It is interesting to notice from \cref{fig:origIpc} that, when a low $n_{mix}$ is used, after four orbits, the profile tends to lose its periodicity property, which is instead preserved for a high $n_{mix}$. This suggests that the 3D \gls{koz} starts to become too wide when no \gls{gmm} is used, while it is represented more accurately when $n_{mix}$ is high enough. In \cref{fig:ipc_comb}, the contribution of each mixand to the ballistic \gls{tipc} is shown for three (\cref{fig:ipc_comb3}) and seven (\cref{fig:ipc_comb7}) mixands. It is clear that having more mixands allows us to refine where the majority of the probability resides. In \cref{fig:ipc_comb3}, the total probability is mostly given by the first mixand. Instead, in \cref{fig:ipc_comb7}, we see that most of the mixands do not contribute at all, being below $10^{-9}$, and the highest contributions are given by mixands 2, 3, and 1. Moreover, using seven mixands allows us to preserve the profile's periodicity, which starts to fade when three mixands are used because $P_{IC,3}$ has a non-negligible effect given by the inflation of the corresponding covariance.

\begin{figure}[tb!]
    \centering
    \subfloat[$n_{mix} = 1$.]{\includegraphics[width = 0.48\textwidth]{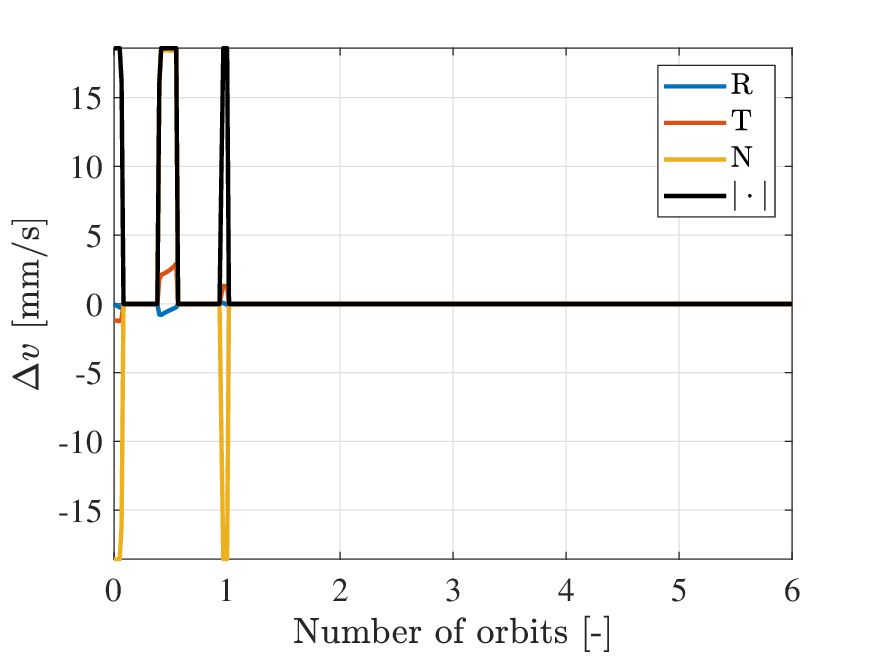}\label{fig:dvLong_1}}
    \subfloat[$n_{mix} = 3$.]{\includegraphics[width = 0.48\textwidth]{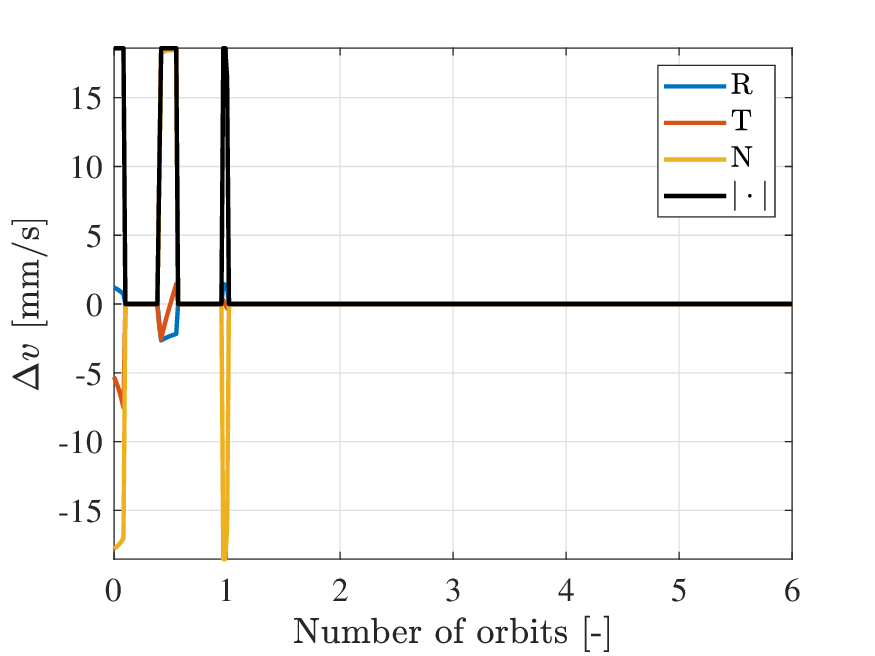}\label{fig:dvLong_2}} \\
    \subfloat[$n_{mix} = 5$.]{\includegraphics[width = 0.48\textwidth]{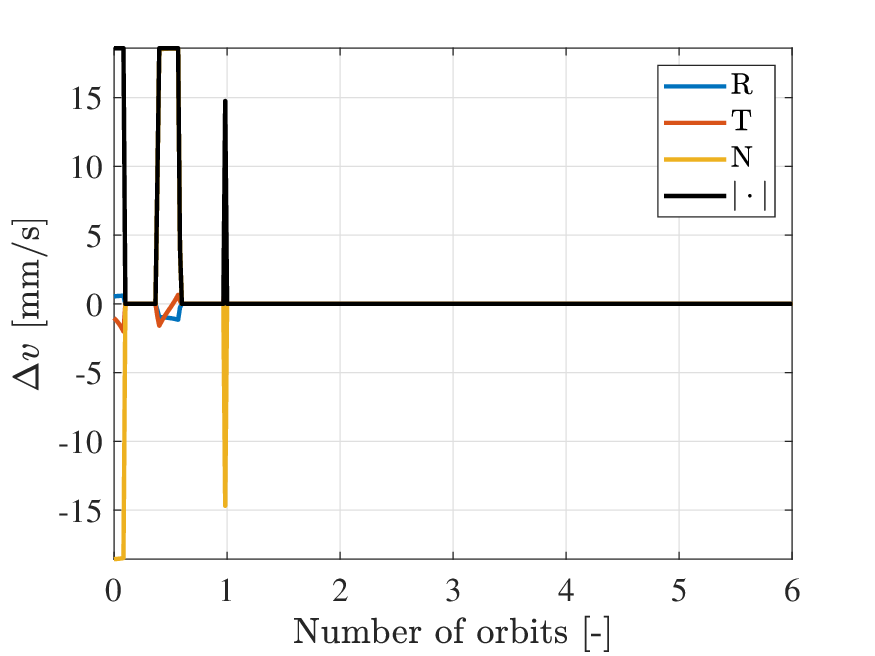}\label{fig:dvLong_3}}
    \subfloat[$n_{mix} = 7$.]{\includegraphics[width = 0.48\textwidth]{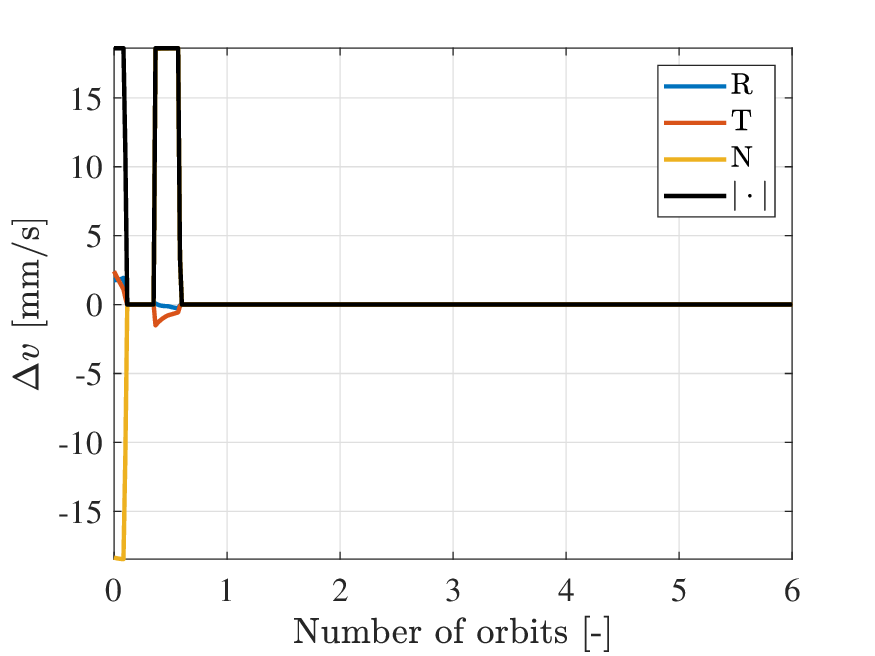}\label{fig:dvLong_4}}
    \caption{Case 3: $\dv$ to perform the \gls{cam}.} 
    \label{fig:dv_Case3}
\end{figure}

In \cref{fig:dv_Case3}, the $\dv$ profiles for the four considered cases are reported. The computed maneuver is very similar, and the normal component is always predominant.
From \cref{tab:res3}, we see that the linearized \gls{smd} approach is much more efficient in this test case: the solution requires fewer major iterations, and the $\dv$ is always lower. Since there is no guarantee for either method to find the global optimum, we observe that the solution found in the two cases falls into two different local minimum wells. In fact, the \gls{tipc} profile in the case with $n_{mix}=1$ had three equivalent points of global maxima. The linearized \gls{smd} method's solution for any $n_mix$ preserves the second maximum (around $1.3$ orbital periods), while the linearized \gls{tipc} one preserves the first maximum (around $0.73$ orbital periods). The solution found by the first method is closer to the global optimum. Moreover, most importantly, the second method with $n_{mix}=7$ finds a solution that is far from the local optimum achieved by $n_{mix}=3$ and $n_{mix}=5$, so the linearization starts to be too coarse when more than five mixands are considered in this test case.

\begin{table}[tb!]
\centering
\caption{Case 3: compared results for the \gls{smd} and \gls{tipc} linearized constraints.}
\label{tab:res3}
    \begin{tabular}{l|l|lll|lll}
    \Xhline{4\arrayrulewidth}
     \textbf{Constraint}  & \multicolumn{4}{c|}{\textbf{Linearized \gls{smd}}} & \multicolumn{3}{c}{\textbf{Linearized \gls{tipc}}}     \\ \hline
     $n_{mix}$ [-]                         & $1$          & $3$          & $5$           & $7$          & $3$          & $5$         & $7$       \\ \Xhline{3\arrayrulewidth}
     $\dv$ [\si{mm/s}]                     & $339.42$     & $342.30$    & $344.12$      & $367.77$     & $401.76$    & $420.78$    & $469.08$  \\ \hline
     $\max_i(P_{TIC,i})\times 10^{-6}$ [-] & $1.006$      & $9.93$       & $0.989$       & $0.829$      & $1.000$      & $1.000$     & $0.196$   \\ \hline
     $n_{major}$ [-]                       & 4            & 9            & 6             & 7            & 7            & 7           & 20        \\ \hline
     $n_{minor}$ [-]                       & 15           & 21           & 30            & 15           & -            & -           & -         \\ \hline
     $e_{validation}$ [\si{mm}]            & $12.73$      & $47.06$      & $1.63$        & $30.79$      & $53.43$      & $44.39$     & $20.32$   \\ \hline
     $T_{sim}$ [\si{s}]                    & $52.85$      & $134.75$     & $164.58$      & $142.02$     & $100.63$     & $114.00$    & $272.64$  \\ \Xhline{4\arrayrulewidth}
    \end{tabular}
\end{table}
\pagebreak
\section{Conclusions}
This work presented a \acrfull{scp} method to compute collision avoidance maneuvers (\glspl{cam}) in scenarios involving multiple encounters between satellites.
Both short-term and long-term encounters were considered. In the first case, the \acrfull{tpoc} was used as a collision metric, whereas in the second, we used the \acrfull{tipc}. Moreover, for the scenarios in which uncertainty propagation is needed for long periods, we employed a \acrfull{gmm} method to propagate the uncertainty of the secondary spacecraft and define multiple keep-out zones.

The fuel-optimal \gls{cam} optimization problem was cast as a \acrfull{socp}, where different techniques were used to convexify the non-convex constraints. Two methods were proposed to refine the initial \gls{scp} solution. The first one relies on a \gls{da}-based linearization of the \gls{tpoc} or \gls{tipc} constraint; the second one updates the \acrfull{poc} or \acrlong{ipc} limits associated with the single conjunctions by means of a \acrlong{nlp}.

The proposed methods were tested on three scenarios. The first one comprised multiple consecutive short-term conjunctions: Chan's \gls{poc} method was used to assess the risk, and the results for the two refinements algorithms were compared. The second scenario was a double encounter with the same secondary where the covariance was propagated with a \gls{gmm}. The last one was a single long-term encounter studied with the aid of \glspl{gmm} with different numbers of mixands. The algorithms can always find an optimal solution for the fuel-optimal problem, but global optimality is not guaranteed. The \gls{smd} linearization refinement method was deemed more suitable because it behaved much better in the presence of significant non-linearities of the original problem.

The proposed method is general, being able to deal with very diverse settings, including (i) both short- and long-term encounters, (ii) arbitrary dynamics models, (iii) arbitrary number of encounters, (iv) nonlinear propagation of the uncertainty of the secondary. The primary's uncertainty is propagated linearly when necessary. Since we assume that accurate orbit determination is performed on the maneuverable satellite, the initial uncertainty of the primary is almost negligible with respect to the secondary. If the uncertainty of the primary is larger, a nonlinear propagation model must also be used for the uncertainty of the primary.

In an ever more cluttered orbital environment, we believe that the complexity of the \gls{ca} scenarios for spacecraft will increase dramatically. Given its light computational burden and proven efficacy and reliability in dealing with multiple conjunctions, the proposed algorithm can be an enabling tool for the safe use of Earth's space in the future.

\section*{Funding Sources}
This material is based upon work supported by the Air Force Office of Scientific Research under award number FA2386-21-1-4115.

\bibliography{references}

\end{document}

%% file: Preamble.tex
\usepackage[english]{babel}

\usepackage{textcomp}
\usepackage[super]{nth}
\usepackage{float}
\usepackage{subfloat}
\usepackage{subcaption}
\usepackage{multirow}
\usepackage{makecell}
\usepackage{longtable,tabularx}
\usepackage{lipsum}
\usepackage{pgfplots}
\usepackage{tikz}
\usetikzlibrary{positioning,external}
\usepackage{algorithm}
\usepackage{algpseudocode}
\usepackage{amsmath}
\usepackage{amstext}
\usepackage{amsfonts}
\usepackage{bm}
\usepackage{mathtools}
\usepackage[bb=boondox]{mathalfa}

\usepackage[nopostdot,acronym,nogroupskip,nonumberlist]{glossaries}
\usepackage[capitalize]{cleveref}
\usepackage{apptools}
\usepackage{graphicx,import}
\usepackage{lastpage}
\usepackage{fancyhdr}
\usepackage{ifthen}
\usepackage[]{siunitx}
\usepackage{algorithm}
\usepackage{algpseudocode}
\usepackage{bm}

\let\oldhat\hat
\renewcommand{\vec}[1]{\boldsymbol{#1}}
\renewcommand{\hat}[1]{\oldhat{\boldsymbol{#1}}}

\newcommand*{\doublerule}{\hrule width \hsize height 1pt \kern 0.5mm \hrule width \hsize height 2pt}
\newcommand\doublerulefill{\leavevmode\leaders\vbox{\hrule width .1pt\kern1pt\hrule}\hfill\kern0pt }
\usetikzlibrary{plotmarks}
\pgfplotsset{compat=newest,
every plot/.append style={color=black, mark=none},
every axis/.append style={
label style={font=\fontsize{9pt}{1em}\color{white!15!black}\selectfont},
tick style={font=\fontsize{9pt}{1em}\selectfont, color=black, line cap=round},
ticklabel style= {font=\fontsize{9pt}{1em}\selectfont},
legend style={legend cell align=left, font=\fontsize{9pt}{1em}\selectfont, align=left, draw=white!15!black},
/pgf/number format/1000 sep={},
yticklabel style={
        /pgf/number format/fixed,
        /pgf/number format/precision=5
},
scaled y ticks=false
}}
\newcommand{\smd}{d_m^2}

\newcommand{\tipc}{P_{TIC}}
\newcommand{\pc}{P_C}
\newcommand{\tpc}{P_{TC}}

\newcommand{\smdlim}{\bar{d}_m^2}

\newcommand{\tipclim}{\bar{P}_{TIC}}

\newcommand{\tpclim}{\bar{P}_{TC}}

\renewcommand{\transp}{^\mathrm{T}}
\newcommand{\dv}{\Delta v}
\newcommand\blfootnote[1]{%
  \begingroup
  \renewcommand\thefootnote{}\footnote{#1}%
  \addtocounter{footnote}{-1}%
  \endgroup
}
\input{glossary.tex}

%% file: glossary.tex
\newacronym{leo}{LEO}{low Earth orbit}
\newacronym{geo}{GEO}{geostationary Earth orbit}
\newacronym{smd}{SMD}{squared Mahalanobis distance}
\newacronym{ssa}{SSA}{Space Situational Awareness}
\newacronym{adr}{ADR}{Active Debris Removal}
\newacronym{sk}{SK}{station keeping}
\newacronym{ca}{CA}{collision avoidance}
\newacronym{ipc}{IPoC}{instantaneous probability of collision}
\newacronym{tipc}{TIPoC}{total instantaneous probability of collision}
\newacronym{cpc}{CPoC}{cumulative probability of collision}
\newacronym{poc}{PoC}{probability of collision}
\newacronym{tpoc}{TPoC}{total probability of collision}
\newacronym{cam}{CAM}{collision avoidance maneuver}
\newacronym{scvx}{SCvx}{successive convexification}
\newacronym{erf}{erf}{error function}
\newacronym{da}{DA}{differential algebra}
\newacronym{eci}{ECI}{Earth Centered Inertial}
\newacronym{ecef}{ECEF}{Earth Centered Earth Fixed}
\newacronym{lvlh}{LVLH}{Local Vertical Local Horizon}
\newacronym{hbr}{HBR}{hard body radius}
\newacronym{pdf}{PDF}{Probability Density Function}
\newacronym{nlp}{NLP}{non-linear program}
\newacronym{srp}{SRP}{solar radiation pressure}
\newacronym{mpbvp}{MPBVP}{multi-point boundary value problem}
\newacronym{cw}{CW}{Clohessy-Wiltshire}
\newacronym{ns}{N-S}{North-South}
\newacronym{sdp}{SDP}{Semidefinite program}
\newacronym{qcqp}{QCQP}{Quadratically constrained quadratic program}
\newacronym{qp}{QP}{quadratic program}
\newacronym{zoh}{ZOH}{zeroth-order hold}
\newacronym{ocp}{OCP}{optimal control problem}
\newacronym{ew}{E-W}{East-West}
\newacronym{stm}{STM}{state transition matrix}
\newacronym{rtn}{RTN}{radial, along-track, cross-track}
\newacronym{iss}{ISS}{International Space Station}
\newacronym{lp}{LP}{linear program}
\newacronym{gnc}{GNC}{guidance, navigation, and control}
\newacronym{koz}{KOZ}{keep-out-zone}
\newacronym{socp}{SOCP}{second-order cone program}
\newacronym{cp}{CP}{convex program}
\newacronym{mrv}{MRV}{multivariate random variable}
\newacronym{gmm}{GMM}{Gaussian mixture model}
\newacronym{rv}{RV}{random variable}
\newacronym{scp}{SCP}{sequential convex program}
\newacronym{rso}{RSO}{resident space object}
\newacronym{mc}{MC}{Monte Carlo}
\newacronym{mss}{MSS}{Mahalanobis Shell Sampling}
\newacronym{nli}{NLI}{nonlinearity index}
\newacronym{rk}{RK}{Runge-Kutta}
\newacronym{trr}{TRR}{trust region radius}
\newacronym{tca}{TCA}{time of closest approach}
\newacronym{gto}{GTO}{Geostationary Transfer Orbit}
\newacronym{aida}{AIDA}{Accurate Integrator for Debris Analysis}
\newacronym{ava}{AvA}{active vs active}
\newacronym{cdm}{CDM}{Conjunction Data Message}

%% file: Diagram.tex
\tikzset{every picture/.style={line width=0.75pt}} 

\begin{tikzpicture}[x=0.75pt,y=0.75pt,yscale=-1,xscale=1]

\draw  [color={rgb, 255:red, 0; green, 0; blue, 0 }  ,draw opacity=0.36 ][fill={rgb, 255:red, 211; green, 170; blue, 99 }  ,fill opacity=1 ] (401.4,181.22) -- (657.4,181.22) -- (657.4,431.22) -- (401.4,431.22) -- cycle ;
\draw  [color={rgb, 255:red, 0; green, 0; blue, 0 }  ,draw opacity=0.36 ][fill={rgb, 255:red, 149; green, 180; blue, 216 }  ,fill opacity=1 ] (0.78,116.01) -- (322.52,116.01) -- (322.52,528.69) -- (0.78,528.69) -- cycle ;
\draw  [color={rgb, 255:red, 0; green, 0; blue, 0 }  ,draw opacity=0.36 ][fill={rgb, 255:red, 0; green, 0; blue, 0 }  ,fill opacity=0.1 ] (3.1,175.7) -- (320.1,175.7) -- (320.1,451.3) -- (3.1,451.3) -- cycle ;
\draw  [color={rgb, 255:red, 0; green, 0; blue, 0 }  ,draw opacity=0.36 ][fill={rgb, 255:red, 0; green, 0; blue, 0 }  ,fill opacity=0.1 ] (7.1,243.2) -- (316.1,243.2) -- (316.1,383.8) -- (7.1,383.8) -- cycle ;
\draw   (230.37,15.72) .. controls (230.37,12.56) and (232.93,10) .. (236.09,10) -- (465.65,10) .. controls (468.81,10) and (471.37,12.56) .. (471.37,15.72) -- (471.37,32.88) .. controls (471.37,36.04) and (468.81,38.6) .. (465.65,38.6) -- (236.09,38.6) .. controls (232.93,38.6) and (230.37,36.04) .. (230.37,32.88) -- cycle ;
\draw   (350,64) -- (410.48,91.63) -- (350,119.25) -- (289.52,91.63) -- cycle ;
\draw    (351.15,38.6) -- (350.8,60.2) ;
\draw [shift={(350.75,63.2)}, rotate = 270.93] [fill={rgb, 255:red, 0; green, 0; blue, 0 }  ][line width=0.08]  [draw opacity=0] (10.72,-5.15) -- (0,0) -- (10.72,5.15) -- (7.12,0) -- cycle    ;
\draw    (289.52,91.63) -- (160.05,90.75) ;
\draw    (162.37,169.7) -- (162.4,201.6) ;
\draw [shift={(162.4,204.6)}, rotate = 269.95] [fill={rgb, 255:red, 0; green, 0; blue, 0 }  ][line width=0.08]  [draw opacity=0] (10.72,-5.15) -- (0,0) -- (10.72,5.15) -- (7.12,0) -- cycle    ;
\draw    (162.4,234.8) -- (162.62,256.35) ;
\draw [shift={(162.65,259.35)}, rotate = 269.42] [fill={rgb, 255:red, 0; green, 0; blue, 0 }  ][line width=0.08]  [draw opacity=0] (10.72,-5.15) -- (0,0) -- (10.72,5.15) -- (7.12,0) -- cycle    ;
\draw    (539.75,91.2) -- (540.35,206.03) ;
\draw [shift={(540.37,209.03)}, rotate = 269.7] [fill={rgb, 255:red, 0; green, 0; blue, 0 }  ][line width=0.08]  [draw opacity=0] (10.72,-5.15) -- (0,0) -- (10.72,5.15) -- (7.12,0) -- cycle    ;
\draw    (539.37,319.7) -- (539.63,364.35) ;
\draw [shift={(539.65,367.35)}, rotate = 269.66] [fill={rgb, 255:red, 0; green, 0; blue, 0 }  ][line width=0.08]  [draw opacity=0] (10.72,-5.15) -- (0,0) -- (10.72,5.15) -- (7.12,0) -- cycle    ;
\draw    (539.75,91.2) -- (410.48,91.63) ;
\draw    (539.45,537.65) -- (539.65,416) ;
\draw    (161.37,536.35) -- (161.65,515) ;
\draw    (539.45,537.65) -- (161.37,536.35) ;
\draw   (326.37,566.61) .. controls (326.37,564.31) and (328.23,562.45) .. (330.52,562.45) -- (370.21,562.45) .. controls (372.51,562.45) and (374.37,564.31) .. (374.37,566.61) -- (374.37,579.08) .. controls (374.37,581.37) and (372.51,583.23) .. (370.21,583.23) -- (330.52,583.23) .. controls (328.23,583.23) and (326.37,581.37) .. (326.37,579.08) -- cycle ;
\draw    (350.41,537) -- (350.37,558.55) ;
\draw [shift={(350.37,561.55)}, rotate = 270.1] [fill={rgb, 255:red, 0; green, 0; blue, 0 }  ][line width=0.08]  [draw opacity=0] (10.72,-5.15) -- (0,0) -- (10.72,5.15) -- (7.12,0) -- cycle    ;
\draw    (540.37,238.03) -- (540.37,276.7) ;
\draw [shift={(540.37,279.7)}, rotate = 270] [fill={rgb, 255:red, 0; green, 0; blue, 0 }  ][line width=0.08]  [draw opacity=0] (10.72,-5.15) -- (0,0) -- (10.72,5.15) -- (7.12,0) -- cycle    ;
\draw    (161.05,90.75) -- (161.35,138.7) ;
\draw [shift={(161.37,141.7)}, rotate = 269.64] [fill={rgb, 255:red, 0; green, 0; blue, 0 }  ][line width=0.08]  [draw opacity=0] (10.72,-5.15) -- (0,0) -- (10.72,5.15) -- (7.12,0) -- cycle    ;
\draw    (162.37,298.9) -- (162.61,319.35) ;
\draw [shift={(162.65,322.35)}, rotate = 269.31] [fill={rgb, 255:red, 0; green, 0; blue, 0 }  ][line width=0.08]  [draw opacity=0] (10.72,-5.15) -- (0,0) -- (10.72,5.15) -- (7.12,0) -- cycle    ;
\draw    (162.65,371) -- (162.65,391.35) ;
\draw [shift={(162.65,394.35)}, rotate = 270] [fill={rgb, 255:red, 0; green, 0; blue, 0 }  ][line width=0.08]  [draw opacity=0] (10.72,-5.15) -- (0,0) -- (10.72,5.15) -- (7.12,0) -- cycle    ;
\draw  [fill={rgb, 255:red, 255; green, 255; blue, 255 }  ,fill opacity=1 ] (161.65,466.35) -- (254.15,490.68) -- (161.65,515) -- (69.15,490.68) -- cycle ;
\draw    (161.65,443) -- (161.65,463.35) ;
\draw [shift={(161.65,466.35)}, rotate = 270] [fill={rgb, 255:red, 0; green, 0; blue, 0 }  ][line width=0.08]  [draw opacity=0] (10.72,-5.15) -- (0,0) -- (10.72,5.15) -- (7.12,0) -- cycle    ;
\draw  [fill={rgb, 255:red, 255; green, 255; blue, 255 }  ,fill opacity=1 ] (72.4,210.72) .. controls (72.4,207.56) and (74.96,205) .. (78.12,205) -- (250.68,205) .. controls (253.84,205) and (256.4,207.56) .. (256.4,210.72) -- (256.4,227.88) .. controls (256.4,231.04) and (253.84,233.6) .. (250.68,233.6) -- (78.12,233.6) .. controls (74.96,233.6) and (72.4,231.04) .. (72.4,227.88) -- cycle ;
\draw  [fill={rgb, 255:red, 255; green, 255; blue, 255 }  ,fill opacity=1 ] (72.15,146.72) .. controls (72.15,143.56) and (74.71,141) .. (77.87,141) -- (242.43,141) .. controls (245.59,141) and (248.15,143.56) .. (248.15,146.72) -- (248.15,163.88) .. controls (248.15,167.04) and (245.59,169.6) .. (242.43,169.6) -- (77.87,169.6) .. controls (74.71,169.6) and (72.15,167.04) .. (72.15,163.88) -- cycle ;
\draw  [fill={rgb, 255:red, 255; green, 255; blue, 255 }  ,fill opacity=1 ] (161.65,394.35) -- (254.15,418.68) -- (161.65,443) -- (69.15,418.68) -- cycle ;
\draw  [fill={rgb, 255:red, 255; green, 255; blue, 255 }  ,fill opacity=1 ] (161.65,322.35) -- (254.15,346.68) -- (161.65,371) -- (69.15,346.68) -- cycle ;
\draw    (254.15,490.68) .. controls (312.08,471.21) and (356.95,161.92) .. (249.99,155.19) ;
\draw [shift={(248.37,155.12)}, rotate = 1.93] [fill={rgb, 255:red, 0; green, 0; blue, 0 }  ][line width=0.08]  [draw opacity=0] (10.72,-5.15) -- (0,0) -- (10.72,5.15) -- (7.12,0) -- cycle    ;
\draw    (254.15,418.68) .. controls (293.76,398.43) and (312.82,225.11) .. (258.91,218.96) ;
\draw [shift={(256.4,218.8)}, rotate = 1.01] [fill={rgb, 255:red, 0; green, 0; blue, 0 }  ][line width=0.08]  [draw opacity=0] (10.72,-5.15) -- (0,0) -- (10.72,5.15) -- (7.12,0) -- cycle    ;
\draw    (254.15,346.68) .. controls (277.88,341.23) and (280.31,282.52) .. (237.1,278.01) ;
\draw [shift={(234.4,277.8)}, rotate = 2.86] [fill={rgb, 255:red, 0; green, 0; blue, 0 }  ][line width=0.08]  [draw opacity=0] (10.72,-5.15) -- (0,0) -- (10.72,5.15) -- (7.12,0) -- cycle    ;
\draw    (447.15,391.68) .. controls (399.13,337.62) and (394.53,231.06) .. (444.1,225.96) ;
\draw [shift={(446.4,225.8)}, rotate = 177.8] [fill={rgb, 255:red, 0; green, 0; blue, 0 }  ][line width=0.08]  [draw opacity=0] (10.72,-5.15) -- (0,0) -- (10.72,5.15) -- (7.12,0) -- cycle    ;
\draw  [fill={rgb, 255:red, 255; green, 255; blue, 255 }  ,fill opacity=1 ] (467.4,288.93) .. controls (467.4,284.55) and (470.95,281) .. (475.33,281) -- (604.47,281) .. controls (608.85,281) and (612.4,284.55) .. (612.4,288.93) -- (612.4,312.72) .. controls (612.4,317.1) and (608.85,320.65) .. (604.47,320.65) -- (475.33,320.65) .. controls (470.95,320.65) and (467.4,317.1) .. (467.4,312.72) -- cycle ;
\draw  [fill={rgb, 255:red, 255; green, 255; blue, 255 }  ,fill opacity=1 ] (89.4,266.93) .. controls (89.4,262.55) and (92.95,259) .. (97.33,259) -- (226.47,259) .. controls (230.85,259) and (234.4,262.55) .. (234.4,266.93) -- (234.4,290.72) .. controls (234.4,295.1) and (230.85,298.65) .. (226.47,298.65) -- (97.33,298.65) .. controls (92.95,298.65) and (89.4,295.1) .. (89.4,290.72) -- cycle ;
\draw  [fill={rgb, 255:red, 255; green, 255; blue, 255 }  ,fill opacity=1 ] (448.4,215.72) .. controls (448.4,212.56) and (450.96,210) .. (454.12,210) -- (626.68,210) .. controls (629.84,210) and (632.4,212.56) .. (632.4,215.72) -- (632.4,232.88) .. controls (632.4,236.04) and (629.84,238.6) .. (626.68,238.6) -- (454.12,238.6) .. controls (450.96,238.6) and (448.4,236.04) .. (448.4,232.88) -- cycle ;
\draw  [fill={rgb, 255:red, 255; green, 255; blue, 255 }  ,fill opacity=1 ] (539.65,367.35) -- (632.15,391.68) -- (539.65,416) -- (447.15,391.68) -- cycle ;

\draw (236.09,13) node [anchor=north west][inner sep=0.75pt]   [align=left] {Solve first SCP run with first guess for $\displaystyle \overline{P}_{C,s}$ };
\draw (324,76) node [anchor=north west][inner sep=0.75pt]  [font=\normalsize] [align=left] {Select \\constraint};
\draw (76.87,148) node [anchor=north west][inner sep=0.75pt]   [align=left] {Adapt limits with Problem \ref{prob:nlpWeights}};
\draw (336.09,566.45) node [anchor=north west][inner sep=0.75pt]   [align=left] {Stop};
\draw (78.4,210.72) node [anchor=north west][inner sep=0.75pt]   [align=left] {Update dynamics with \cref{eq:dynamics}};
\draw (3.1,175.7) node [anchor=north west][inner sep=0.75pt]  [color={rgb, 255:red, 255; green, 255; blue, 255 }  ,opacity=1 ] [align=left] {SCP run - Major iterations};
\draw (401.4,181.22) node [anchor=north west][inner sep=0.75pt]  [color={rgb, 255:red, 255; green, 255; blue, 255 }  ,opacity=1 ] [align=left] {SCP run};
\draw (127,404) node [anchor=north west][inner sep=0.75pt]   [align=left] {Convergence\\with \cref{eq:majConv}};
\draw (125,332) node [anchor=north west][inner sep=0.75pt]   [align=left] {Convergence\\with \cref{eq:minConv}};
\draw (158,72) node [anchor=north west][inner sep=0.75pt]   [align=left] {Linearized SMD};
\draw (442,72) node [anchor=north west][inner sep=0.75pt]   [align=left] {Linearized TPoC};
\draw (7.1,243.2) node [anchor=north west][inner sep=0.75pt]  [color={rgb, 255:red, 255; green, 255; blue, 255 }  ,opacity=1 ] [align=left] {Minor iterations};
\draw (129,476) node [anchor=north west][inner sep=0.75pt]   [align=left] {Check TPoC\\compliance};
\draw (0.78,116.01) node [anchor=north west][inner sep=0.75pt]  [color={rgb, 255:red, 255; green, 255; blue, 255 }  ,opacity=1 ] [align=left] {Limits refinement};
\draw (476.65,285.17) node [anchor=north west][inner sep=0.75pt]   [align=left] {Solve SOCP \cref{eq:optFinal} \\with TPoC constraint};
\draw (99.65,263.17) node [anchor=north west][inner sep=0.75pt]   [align=left] {Solve SOCP \cref{eq:optFinal} \\with SMD constraint};
\draw (454.4,215.72) node [anchor=north west][inner sep=0.75pt]   [align=left] {Update dynamics with \cref{eq:dynamics}};
\draw (505,377) node [anchor=north west][inner sep=0.75pt]   [align=left] {Convergence\\with \cref{eq:majConv}};

\end{tikzpicture}

%% file: limAdaptScheme.tex
\tikzset{every picture/.style={line width=0.75pt}} 

\begin{tikzpicture}[x=0.75pt,y=0.75pt,yscale=-1,xscale=1]

\draw  [color={rgb, 255:red, 208; green, 2; blue, 27 }  ,draw opacity=1 ][fill={rgb, 255:red, 208; green, 2; blue, 27 }  ,fill opacity=0.2 ] (316.66,91.25) .. controls (332.09,66.56) and (407.02,85.57) .. (484.02,133.71) .. controls (561.02,181.84) and (610.93,240.87) .. (595.5,265.56) .. controls (580.07,290.24) and (505.14,271.23) .. (428.14,223.1) .. controls (351.14,174.96) and (301.23,115.93) .. (316.66,91.25) -- cycle ;
\draw  [color={rgb, 255:red, 208; green, 2; blue, 27 }  ,draw opacity=1 ][fill={rgb, 255:red, 208; green, 2; blue, 27 }  ,fill opacity=0.2 ][dash pattern={on 4.5pt off 4.5pt}] (85.06,235.54) .. controls (59.45,213.73) and (70,159.28) .. (108.62,113.92) .. controls (147.24,68.55) and (199.31,49.45) .. (224.93,71.26) .. controls (250.54,93.06) and (239.99,147.51) .. (201.37,192.88) .. controls (162.75,238.24) and (110.67,257.34) .. (85.06,235.54) -- cycle ;
\draw  [color={rgb, 255:red, 208; green, 2; blue, 27 }  ,draw opacity=1 ][fill={rgb, 255:red, 208; green, 2; blue, 27 }  ,fill opacity=0.4 ][dash pattern={on 4.5pt off 4.5pt}] (371.04,125.24) .. controls (380.45,110.19) and (426.16,121.78) .. (473.12,151.14) .. controls (520.09,180.5) and (550.53,216.5) .. (541.12,231.56) .. controls (531.71,246.62) and (486,235.02) .. (439.04,205.66) .. controls (392.07,176.3) and (361.63,140.3) .. (371.04,125.24) -- cycle ;
\draw  [color={rgb, 255:red, 208; green, 2; blue, 27 }  ,draw opacity=1 ][fill={rgb, 255:red, 208; green, 2; blue, 27 }  ,fill opacity=0.4 ] (112.81,202.94) .. controls (99.44,191.57) and (107.5,160.16) .. (130.79,132.79) .. controls (154.09,105.43) and (183.81,92.47) .. (197.18,103.85) .. controls (210.54,115.23) and (202.49,146.64) .. (179.2,174) .. controls (155.9,201.37) and (126.18,214.32) .. (112.81,202.94) -- cycle ;
\draw  [color={rgb, 255:red, 0; green, 0; blue, 0 }  ,draw opacity=1 ][fill={rgb, 255:red, 0; green, 0; blue, 0 }  ,fill opacity=1 ] (164.37,116.8) -- (168.1,121.7) -- (163.2,125.44) -- (159.46,120.53) -- cycle ;
\draw  [color={rgb, 255:red, 0; green, 0; blue, 0 }  ,draw opacity=1 ][fill={rgb, 255:red, 0; green, 0; blue, 0 }  ,fill opacity=1 ] (410.1,143.54) -- (412.75,149.1) -- (407.18,151.75) -- (404.54,146.18) -- cycle ;
\draw  [color={rgb, 255:red, 57; green, 109; blue, 165 }  ,draw opacity=1 ][fill={rgb, 255:red, 74; green, 144; blue, 226 }  ,fill opacity=1 ] (281.67,90.14) .. controls (282.24,88.01) and (284.43,86.74) .. (286.56,87.31) .. controls (288.7,87.87) and (289.97,90.07) .. (289.4,92.2) .. controls (288.83,94.33) and (286.64,95.6) .. (284.51,95.04) .. controls (282.37,94.47) and (281.1,92.28) .. (281.67,90.14) -- cycle ;
\draw    (164.37,116.8) -- (168.04,74) ;
\draw [shift={(168.21,72.01)}, rotate = 94.91] [fill={rgb, 255:red, 0; green, 0; blue, 0 }  ][line width=0.08]  [draw opacity=0] (12,-3) -- (0,0) -- (12,3) -- cycle    ;
\draw    (404.54,146.18) -- (291.21,93.05) ;
\draw [shift={(289.4,92.2)}, rotate = 25.12] [fill={rgb, 255:red, 0; green, 0; blue, 0 }  ][line width=0.08]  [draw opacity=0] (12,-3) -- (0,0) -- (12,3) -- cycle    ;
\draw  [color={rgb, 255:red, 39; green, 71; blue, 3 }  ,draw opacity=1 ][fill={rgb, 255:red, 65; green, 117; blue, 5 }  ,fill opacity=1 ] (161.74,103.56) .. controls (161.83,101.35) and (163.69,99.63) .. (165.89,99.72) .. controls (168.1,99.81) and (169.82,101.67) .. (169.74,103.87) .. controls (169.65,106.08) and (167.79,107.8) .. (165.58,107.71) .. controls (163.38,107.63) and (161.66,105.77) .. (161.74,103.56) -- cycle ;
\draw  [color={rgb, 255:red, 39; green, 71; blue, 3 }  ,draw opacity=1 ][fill={rgb, 255:red, 65; green, 117; blue, 5 }  ,fill opacity=1 ] (309.73,101.91) .. controls (310.37,99.8) and (312.59,98.59) .. (314.71,99.22) .. controls (316.83,99.85) and (318.03,102.08) .. (317.4,104.2) .. controls (316.77,106.32) and (314.54,107.52) .. (312.42,106.89) .. controls (310.31,106.26) and (309.1,104.03) .. (309.73,101.91) -- cycle ;
\draw  [color={rgb, 255:red, 0; green, 0; blue, 0 }  ,draw opacity=1 ][fill={rgb, 255:red, 208; green, 2; blue, 27 }  ,fill opacity=1 ] (151.16,152.25) .. controls (151.79,150.14) and (154.02,148.93) .. (156.14,149.56) .. controls (158.25,150.19) and (159.46,152.42) .. (158.83,154.54) .. controls (158.2,156.66) and (155.97,157.86) .. (153.85,157.23) .. controls (151.73,156.6) and (150.53,154.37) .. (151.16,152.25) -- cycle ;
\draw  [color={rgb, 255:red, 0; green, 0; blue, 0 }  ,draw opacity=1 ][fill={rgb, 255:red, 208; green, 2; blue, 27 }  ,fill opacity=1 ] (452.89,175.98) .. controls (454.23,174.22) and (456.74,173.88) .. (458.5,175.22) .. controls (460.26,176.55) and (460.6,179.06) .. (459.27,180.82) .. controls (457.93,182.58) and (455.42,182.92) .. (453.66,181.59) .. controls (451.9,180.25) and (451.56,177.74) .. (452.89,175.98) -- cycle ;
\draw  [dash pattern={on 0.84pt off 2.51pt}]  (452.89,175.98) -- (414.41,150.22) ;
\draw [shift={(412.75,149.1)}, rotate = 33.81] [fill={rgb, 255:red, 0; green, 0; blue, 0 }  ][line width=0.08]  [draw opacity=0] (12,-3) -- (0,0) -- (12,3) -- cycle    ;
\draw  [dash pattern={on 0.84pt off 2.51pt}]  (156.14,149.56) -- (162.64,127.36) ;
\draw [shift={(163.2,125.44)}, rotate = 106.31] [fill={rgb, 255:red, 0; green, 0; blue, 0 }  ][line width=0.08]  [draw opacity=0] (12,-3) -- (0,0) -- (12,3) -- cycle    ;
\draw  [color={rgb, 255:red, 57; green, 109; blue, 165 }  ,draw opacity=1 ][fill={rgb, 255:red, 74; green, 144; blue, 226 }  ,fill opacity=1 ] (164.37,67.85) .. controls (164.46,65.65) and (166.31,63.93) .. (168.52,64.01) .. controls (170.73,64.1) and (172.45,65.96) .. (172.36,68.16) .. controls (172.28,70.37) and (170.42,72.09) .. (168.21,72.01) .. controls (166,71.92) and (164.28,70.06) .. (164.37,67.85) -- cycle ;

\draw (162.52,42.01) node [anchor=north west][inner sep=0.75pt]   [align=left] {$\displaystyle \Delta \vec{r}_{1}^{0}$};
\draw (165.37,116.8) node [anchor=north west][inner sep=0.75pt]   [align=left] {$\displaystyle \Delta \vec{r}_{1}^{ball}$};
\draw (168.89,85.72) node [anchor=north west][inner sep=0.75pt]   [align=left] {$\displaystyle \Delta \vec{r}_{1}$};
\draw (400.54,127.18) node [anchor=north west][inner sep=0.75pt]   [align=left] {$\displaystyle \Delta \vec{r}_{2}^{ball}$};
\draw (291.66,109.25) node [anchor=north west][inner sep=0.75pt]   [align=left] {$\displaystyle \Delta \vec{r}_{2}$};
\draw (252.66,94.25) node [anchor=north west][inner sep=0.75pt]   [align=left] {$\displaystyle \Delta \vec{r}_{2}^{0}$};
\draw (134.37,147.8) node [anchor=north west][inner sep=0.75pt]  [color={rgb, 255:red, 126; green, 6; blue, 47 }  ,opacity=1 ] [align=left] {$\displaystyle \textcolor[rgb]{0.55,0,0.07}{\vec{r}_{1}}$};
\draw (463.5,167.22) node [anchor=north west][inner sep=0.75pt]  [color={rgb, 255:red, 126; green, 6; blue, 47 }  ,opacity=1 ] [align=left] {$\displaystyle \textcolor[rgb]{0.55,0,0.07}{\vec{r}}\textcolor[rgb]{0.55,0,0.07}{_{2}}$};

\end{tikzpicture}

%% file: HeadToHeadEnc.tex
\tikzset{every picture/.style={line width=0.75pt}} 

\begin{tikzpicture}[x=0.75pt,y=0.75pt,yscale=-1,xscale=1]

\draw [color={rgb, 255:red, 208; green, 2; blue, 27 }  ,draw opacity=1 ][line width=1.5]    (145.85,124.82) -- (121.56,140.79) ;
\draw [shift={(149.19,122.62)}, rotate = 146.67] [fill={rgb, 255:red, 208; green, 2; blue, 27 }  ,fill opacity=1 ][line width=0.08]  [draw opacity=0] (13.4,-6.43) -- (0,0) -- (13.4,6.44) -- (8.9,0) -- cycle    ;
\draw [color={rgb, 255:red, 208; green, 2; blue, 27 }  ,draw opacity=1 ][line width=1.5]    (255.99,117.91) -- (222.36,123.59) ;
\draw [shift={(259.93,117.24)}, rotate = 170.42] [fill={rgb, 255:red, 208; green, 2; blue, 27 }  ,fill opacity=1 ][line width=0.08]  [draw opacity=0] (13.4,-6.43) -- (0,0) -- (13.4,6.44) -- (8.9,0) -- cycle    ;
\draw [color={rgb, 255:red, 208; green, 2; blue, 27 }  ,draw opacity=1 ][line width=1.5]    (355.6,132.64) -- (323.69,121.97) ;
\draw [shift={(359.39,133.91)}, rotate = 198.49] [fill={rgb, 255:red, 208; green, 2; blue, 27 }  ,fill opacity=1 ][line width=0.08]  [draw opacity=0] (13.4,-6.43) -- (0,0) -- (13.4,6.44) -- (8.9,0) -- cycle    ;
\draw   (114.64,124.45) .. controls (128.85,118.43) and (143.48,120.86) .. (147.3,129.88) .. controls (151.12,138.91) and (142.69,151.1) .. (128.48,157.12) .. controls (114.26,163.15) and (99.64,160.72) .. (95.82,151.69) .. controls (91.99,142.67) and (100.42,130.48) .. (114.64,124.45) -- cycle ;
\draw   (208.93,100.3) .. controls (236.07,84.66) and (264.08,82.4) .. (271.49,95.26) .. controls (278.91,108.12) and (262.92,131.23) .. (235.78,146.87) .. controls (208.64,162.52) and (180.63,164.77) .. (173.22,151.91) .. controls (165.81,139.05) and (181.79,115.94) .. (208.93,100.3) -- cycle ;
\draw   (327.66,102.74) .. controls (344.18,106.15) and (355.78,117.52) .. (353.59,128.14) .. controls (351.4,138.77) and (336.24,144.62) .. (319.72,141.21) .. controls (303.21,137.8) and (291.6,126.43) .. (293.8,115.8) .. controls (295.99,105.18) and (311.15,99.33) .. (327.66,102.74) -- cycle ;
\draw  [fill={rgb, 255:red, 208; green, 2; blue, 27 }  ,fill opacity=1 ] (119.03,140.79) .. controls (119.03,139.39) and (120.16,138.26) .. (121.56,138.26) .. controls (122.95,138.26) and (124.08,139.39) .. (124.08,140.79) .. controls (124.08,142.18) and (122.95,143.32) .. (121.56,143.32) .. controls (120.16,143.32) and (119.03,142.18) .. (119.03,140.79) -- cycle ;
\draw  [fill={rgb, 255:red, 208; green, 2; blue, 27 }  ,fill opacity=1 ] (219.83,123.59) .. controls (219.83,122.19) and (220.96,121.06) .. (222.36,121.06) .. controls (223.75,121.06) and (224.88,122.19) .. (224.88,123.59) .. controls (224.88,124.98) and (223.75,126.11) .. (222.36,126.11) .. controls (220.96,126.11) and (219.83,124.98) .. (219.83,123.59) -- cycle ;
\draw  [fill={rgb, 255:red, 208; green, 2; blue, 27 }  ,fill opacity=1 ] (321.32,121.11) .. controls (321.8,119.8) and (323.25,119.12) .. (324.56,119.6) .. controls (325.87,120.08) and (326.54,121.53) .. (326.07,122.84) .. controls (325.59,124.15) and (324.14,124.83) .. (322.83,124.35) .. controls (321.52,123.87) and (320.84,122.42) .. (321.32,121.11) -- cycle ;
\draw  [dash pattern={on 4.5pt off 4.5pt}]  (121.56,148.13) -- (121.23,209.11) ;
\draw  [dash pattern={on 4.5pt off 4.5pt}]  (222.36,126.53) -- (221.31,208.24) ;
\draw  [dash pattern={on 4.5pt off 4.5pt}]  (323.69,124.5) -- (324.01,209.11) ;
\draw [line width=1.5]    (70.42,151.73) .. controls (119.43,135.49) and (348.04,126.51) .. (380.89,130.68) ;
\draw [shift={(66.4,153.26)}, rotate = 336.04] [fill={rgb, 255:red, 0; green, 0; blue, 0 }  ][line width=0.08]  [draw opacity=0] (13.4,-6.43) -- (0,0) -- (13.4,6.44) -- (8.9,0) -- cycle    ;

\draw (97.87,205.13) node [anchor=north west][inner sep=0.75pt]    {$t_{CA,12}$};
\draw (199.47,204.59) node [anchor=north west][inner sep=0.75pt]    {$t_{CA,22}$};
\draw (301.62,204.38) node [anchor=north west][inner sep=0.75pt]    {$t_{CA,32}$};
\draw (104.49,69.19) node [anchor=north west][inner sep=0.75pt]    {$c=1$};
\draw (204.25,66.37) node [anchor=north west][inner sep=0.75pt]    {$c=2$};
\draw (304.39,67.3) node [anchor=north west][inner sep=0.75pt]    {$c=3$};
\draw (139.14,107.83) node [anchor=north west][inner sep=0.75pt]    {$\textcolor[rgb]{0.82,0.01,0.11}{v_{12}}$};
\draw (231.81,102) node [anchor=north west][inner sep=0.75pt]    {$\textcolor[rgb]{0.82,0.01,0.11}{v}\textcolor[rgb]{0.82,0.01,0.11}{_{22}}$};
\draw (353.59,138.14) node [anchor=north west][inner sep=0.75pt]    {$\textcolor[rgb]{0.82,0.01,0.11}{v}\textcolor[rgb]{0.82,0.01,0.11}{_{32}}$};

\end{tikzpicture}

%% file: PerpendicularEnc.tex
\tikzset{every picture/.style={line width=0.75pt}} 

\begin{tikzpicture}[x=0.75pt,y=0.75pt,yscale=-1,xscale=1]

\draw  [dash pattern={on 4.5pt off 4.5pt}]  (545.92,55.02) -- (544.19,208.24) ;
\draw [color={rgb, 255:red, 208; green, 2; blue, 27 }  ,draw opacity=1 ][line width=1.5]    (556.26,85.89) -- (553.92,52.7) ;
\draw [shift={(556.54,89.88)}, rotate = 265.97] [fill={rgb, 255:red, 208; green, 2; blue, 27 }  ,fill opacity=1 ][line width=0.08]  [draw opacity=0] (13.4,-6.43) -- (0,0) -- (13.4,6.44) -- (8.9,0) -- cycle    ;
\draw [color={rgb, 255:red, 208; green, 2; blue, 27 }  ,draw opacity=1 ][line width=1.5]    (538.17,144.31) -- (546.06,111.79) ;
\draw [shift={(537.23,148.19)}, rotate = 283.63] [fill={rgb, 255:red, 208; green, 2; blue, 27 }  ,fill opacity=1 ][line width=0.08]  [draw opacity=0] (13.4,-6.43) -- (0,0) -- (13.4,6.44) -- (8.9,0) -- cycle    ;
\draw [color={rgb, 255:red, 208; green, 2; blue, 27 }  ,draw opacity=1 ][line width=1.5]    (554.99,190.22) -- (552.56,166.44) ;
\draw [shift={(555.4,194.2)}, rotate = 264.16] [fill={rgb, 255:red, 208; green, 2; blue, 27 }  ,fill opacity=1 ][line width=0.08]  [draw opacity=0] (13.4,-6.43) -- (0,0) -- (13.4,6.44) -- (8.9,0) -- cycle    ;
\draw   (537.98,49.56) .. controls (540.71,35.69) and (550.07,25.85) .. (558.87,27.58) .. controls (567.68,29.31) and (572.6,41.97) .. (569.87,55.84) .. controls (567.13,69.72) and (557.78,79.56) .. (548.98,77.83) .. controls (540.17,76.09) and (535.25,63.44) .. (537.98,49.56) -- cycle ;
\draw   (526.17,96.15) .. controls (544.4,74) and (567.69,63.04) .. (578.19,71.67) .. controls (588.69,80.31) and (582.43,105.28) .. (564.2,127.44) .. controls (545.97,149.59) and (522.68,160.55) .. (512.18,151.92) .. controls (501.68,143.28) and (507.94,118.31) .. (526.17,96.15) -- cycle ;
\draw   (539.91,153.64) .. controls (550.9,142.79) and (565.47,139.72) .. (572.45,146.78) .. controls (579.43,153.85) and (576.19,168.38) .. (565.21,179.24) .. controls (554.22,190.09) and (539.66,193.16) .. (532.67,186.09) .. controls (525.69,179.03) and (528.93,164.5) .. (539.91,153.64) -- cycle ;
\draw  [fill={rgb, 255:red, 208; green, 2; blue, 27 }  ,fill opacity=1 ] (551.61,52.7) .. controls (551.61,51.43) and (552.65,50.39) .. (553.92,50.39) .. controls (555.2,50.39) and (556.24,51.43) .. (556.24,52.7) .. controls (556.24,53.98) and (555.2,55.02) .. (553.92,55.02) .. controls (552.65,55.02) and (551.61,53.98) .. (551.61,52.7) -- cycle ;
\draw  [fill={rgb, 255:red, 208; green, 2; blue, 27 }  ,fill opacity=1 ] (543.74,111.79) .. controls (543.74,110.52) and (544.78,109.48) .. (546.06,109.48) .. controls (547.33,109.48) and (548.37,110.52) .. (548.37,111.79) .. controls (548.37,113.07) and (547.33,114.11) .. (546.06,114.11) .. controls (544.78,114.11) and (543.74,113.07) .. (543.74,111.79) -- cycle ;
\draw  [fill={rgb, 255:red, 208; green, 2; blue, 27 }  ,fill opacity=1 ] (550.39,165.65) .. controls (550.82,164.45) and (552.15,163.83) .. (553.35,164.26) .. controls (554.55,164.7) and (555.17,166.03) .. (554.74,167.23) .. controls (554.3,168.43) and (552.97,169.05) .. (551.77,168.61) .. controls (550.57,168.18) and (549.95,166.85) .. (550.39,165.65) -- cycle ;
\draw [line width=1.5]    (406.72,137.35) .. controls (453.07,122.61) and (660.46,114.49) .. (690.4,118.29) ;
\draw [shift={(402.33,138.98)}, rotate = 336.04] [fill={rgb, 255:red, 0; green, 0; blue, 0 }  ][line width=0.08]  [draw opacity=0] (13.4,-6.43) -- (0,0) -- (13.4,6.44) -- (8.9,0) -- cycle    ;

\draw (476.25,203.75) node [anchor=north west][inner sep=0.75pt]    {$t_{CA,12} \approx t_{CA,22} \approx t_{CA,32}$};
\draw (580.05,32.34) node [anchor=north west][inner sep=0.75pt]    {$c=1$};
\draw (581.87,84.44) node [anchor=north west][inner sep=0.75pt]    {$c=2$};
\draw (582.62,147.57) node [anchor=north west][inner sep=0.75pt]    {$c=3$};
\draw (561.23,77.29) node [anchor=north west][inner sep=0.75pt]    {$\textcolor[rgb]{0.82,0.01,0.11}{v_{12}}$};
\draw (511.78,129.59) node [anchor=north west][inner sep=0.75pt]    {$\textcolor[rgb]{0.82,0.01,0.11}{v}\textcolor[rgb]{0.82,0.01,0.11}{_{22}}$};
\draw (565.21,179.24) node [anchor=north west][inner sep=0.75pt]    {$\textcolor[rgb]{0.82,0.01,0.11}{v}\textcolor[rgb]{0.82,0.01,0.11}{_{32}}$};

\end{tikzpicture}

%% file: gridGmm.pdf_tex
\begingroup%
  \makeatletter%
  \providecommand\color[2][]{%
    \errmessage{(Inkscape) Color is used for the text in Inkscape, but the package 'color.sty' is not loaded}%
    \renewcommand\color[2][]{}%
  }%
  \providecommand\transparent[1]{%
    \errmessage{(Inkscape) Transparency is used (non-zero) for the text in Inkscape, but the package 'transparent.sty' is not loaded}%
    \renewcommand\transparent[1]{}%
  }%
  \providecommand\rotatebox[2]{#2}%
  \newcommand*\fsize{\dimexpr\f@size pt\relax}%
  \newcommand*\lineheight[1]{\fontsize{\fsize}{#1\fsize}\selectfont}%
  \ifx\svgwidth\undefined%
    \setlength{\unitlength}{960bp}%
    \ifx\svgscale\undefined%
      \relax%
    \else%
      \setlength{\unitlength}{\unitlength * \real{\svgscale}}%
    \fi%
  \else%
    \setlength{\unitlength}{\svgwidth}%
  \fi%
  \global\let\svgwidth\undefined%
  \global\let\svgscale\undefined%
  \makeatother%
  \begin{picture}(1,0.078125)%
    \lineheight{1}%
    \setlength\tabcolsep{0pt}%
    \put(0.06391107,0.01254328){\color[rgb]{0,0,0}\makebox(0,0)[lt]{\lineheight{1.25}\smash{\begin{tabular}[t]{l}$t_0$\end{tabular}}}}%
    \put(0.27659527,0.01040111){\color[rgb]{0,0,0}\makebox(0,0)[lt]{\lineheight{1.25}\smash{\begin{tabular}[t]{l}$t_{CA,1}$\end{tabular}}}}%
    \put(0.79504348,0.01164289){\color[rgb]{0.00392157,0.00392157,0.00392157}\makebox(0,0)[lt]{\lineheight{1.25}\smash{\begin{tabular}[t]{l}$t_{CA,2}$\end{tabular}}}}%
    \put(0,0){\includegraphics[width=\unitlength,page=1]{gridGmm.pdf}}%
  \end{picture}%
\endgroup%